\documentclass[]{siamltex1213}

\usepackage{amsmath}
\usepackage{amsfonts}
\usepackage[dvips]{graphicx}

\DeclareMathOperator{\erfc}{erfc}
\DeclareMathOperator{\Arg}{Arg}
\DeclareMathOperator{\Or}{O}

\title{Orientation-dependent pinning and homoclinic snaking on a planar lattice}

\author{Andrew D. Dean\footnotemark[2]
\and Paul C. Matthews\footnotemark[3]
\and Stephen M. Cox\footnotemark[3]
\and John R. King\footnotemark[3]}

\begin{document}
\maketitle

\newcommand{\od}{\, \mathrm{d}}
\newcommand{\pd}{\partial}
\newcommand{\cc}{\mathrm{c.c.}}
\newcommand{\rhs}{\mathrm{RHS}}
\newcommand{\mb}[1]{\mathbf{#1}}
\newcommand{\lb}{\left}
\newcommand{\rb}{\right}
\newcommand{\fr}[2]{\frac{#1}{#2}}


\renewcommand{\thefootnote}{\fnsymbol{footnote}}
\footnotetext[2]{Department of Biology, University of York, Heslington, York YO10 5DD, UK}
\footnotetext[3]{Mathematical Sciences, University Park, Nottingham NG7 2RD, UK}
\renewcommand{\thefootnote}{\arabic{footnote}}

\pagestyle{myheadings}
\thispagestyle{plain}
\markboth{A.D. Dean, P.C. Matthews, S.M. Cox and J.R. King}{Orientation-dependent pinning and snaking}


\begin{abstract}
We study homoclinic snaking of one-dimensional, localised states on two-dimensional, bistable lattices via the method of exponential asymptotics. Within a narrow region of parameter space, fronts connecting the two stable states are pinned to the underlying lattice. Localised solutions are formed by matching two such stationary fronts back-to-back; depending on the orientation relative to the lattice, the solution branch may `snake' back and forth within the pinning region via successive saddle-node bifurcations. Standard continuum approximations in the weakly nonlinear limit (equivalently, the limit of small mesh size) do not exhibit this behaviour, due to the resultant leading-order reaction-diffusion equation lacking a periodic spatial structure. By including exponentially small effects hidden beyond all algebraic orders in the asymptotic expansion, we find that exponentially small but exponentially growing terms are switched on via error function smoothing near Stokes lines. Eliminating these otherwise unbounded beyond-all-orders terms selects the origin (modulo the mesh size) of the front, and matching two fronts together yields a set of equations describing the snaking bifurcation diagram. This is possible only within an exponentially small region of parameter space---the pinning region.  Moreover, by considering fronts orientated at an arbitrary angle $\psi$ to the $x$-axis, we show that the width of the pinning region is non-zero only if $\tan\psi$ is rational or infinite. The asymptotic results are compared with numerical calculations, with good agreement.
\end{abstract}




\section{Introduction}
\label{Sec_Intro}

The phenomenon known as homoclinic snaking, referring to the existence of a multiplicity of localised solutions within a narrow region of parameter space, has been observed in a wide variety of experimental and theoretical contexts \cite{batiste2006spatially, coombes2003waves, lee1994experimental, richter2005two, umbanhowar1996localized}, and has been the subject of much research over the past decade or so \cite{champneys1998homoclinic, dawes2010emergence, knobloch2008spatially}. While much of the literature is focused on continuous systems, snaking also occurs in discrete problems \cite{taylor2010snaking}. One pertinent physical example which has received much recent attention is nonlinear optics \cite{bortolozzo2004bistability, chong2009multistable, firth2007homoclinic, tlidi1994localized, vladimirov2002two, yulin2010discrete, yulin2011snake}, not least due to the potential use of `cavity solitons' as a basis for purely optical information storage and processing \cite{peschel2003spatial}. Furthermore, snaking has recently been observed in a model of plant hormone distribution \cite{draelants2014localised}; although this is the only such example we have been able to find, the methods therein, as well as those of the present work and, for example, \cite{taylor2010snaking}, are applicable to a wide class of problems, and so we expect there to be many more instances of snaking in problems pertaining to cellular biology waiting to be discovered.

From a purely theoretical perspective, there are two main advantages to studying snaking in a discrete context. First, numerical calculations are much more straightforward, with no need to discretise the system. Secondly, snaking is found in second-order systems, in contrast to the continuous case in which a fourth-order system is necessary (e.g. the Swift-Hohenberg equation \cite{beck2009snakes, burke2006localized, swift1977hydrodynamic}). This is significant: the bifurcation diagrams in both cases are remarkably similar, and so discrete problems provide a relatively simple context in which to study the snaking phenomenon analytically. This will facilitate the analysis of more complicated snaking phenomena than have heretofore been considered, in particular the move from one dimension to two \cite{avitabile2010snake, chong2009multistable, lloyd2008localized, lloyd2009localized, taylor2010snaking}, one of the major challenges in snaking theory \cite{knobloch2008spatially}. In the present work we take a first step towards higher-dimensional snaking by considering fronts (the building blocks of localised solutions) which are oriented at an arbitrary angle $\psi$ relative to the lattice; we refer to such solutions as having been `rotated into the plane'. We will show that the width of the snaking region in parameter space depends discontinuously on the front orientation, vanishing when $\tan\psi$ is irrational. A similar result was recently derived by Kozyreff and Chapman \cite{kozyreff2013analytical} for a wide range of continuous problems exhibiting a Turing instability to a hexagonal pattern. The more general findings of \cite{kozyreff2013analytical} are complemented by the present work, being a much more detailed study of a closely related problem (i.e. a square lattice rather than a hexagonal pattern). In particular, we present a full asymptotic description of the snaking bifurcation diagram, a feature lacking from \cite{kozyreff2013analytical}, along with a more complete formula for the width of the snaking region; our results therefore confirm the general claims made in \cite{kozyreff2013analytical}.

In one dimension, localised solutions (homoclinic connections) in bistable systems are constructed by gluing together stationary back-to-back fronts (heteroclinic connections) between the two stable states \cite{beck2009snakes}. Although in general fronts drift, for a certain range of parameter values they pin to the underlying lattice and are stationary. The periodic structure provides an energy barrier; within the pinning region, fronts lack sufficient energy to de-pin, and, as they are stationary, can be used to construct localised solutions. Outside the pinning region, back-to-back fronts either annihilate one another or grow so the localised patch fills the entire domain, depending on the direction of drift. Pinning of fronts is a well-known feature of discrete problems, observed theoretically in models of nerve cells \cite{carpio2011propagation}, discrete reaction-diffusion systems \cite{fath1998propagation} and elastic crystals \cite{carpio1997dynamics, king2001asymptotics}, and experimentally in coupled chemical reactors \cite{laplante1992propagation}, to name a few examples. To our knowledge, the associated homoclinic snaking has not yet been studied explicitly in such systems; however, it occurs as a direct consequence of the pinning of fronts and therefore will be readily observable.

The localised solution branches bifurcate from the primary, constant solution branch, becoming progressively more localised until they enter the pinning region (where they comprise back-to-back fronts) and begin to snake, turning back and forth between two asymptotes via successive saddle-node bifurcations \cite{beck2009snakes, burke2006localized, burke2007homoclinic}. Typically snaking branches occur in pairs; both have reflection symmetry, but one branch is centred on a lattice point (site-centred) and the other is centred midway between two consecutive points (bond-centred). These are linked by `rungs' of asymmetric solutions. This structure is nearly identical to the continuous case \cite{beck2009snakes, burke2007snakes}, although the symmetries of the snaking solutions are different; note that in continuous systems, pinning is due to slowly-varying fronts locking to underlying fast spatial oscillations \cite{chapman2009exponential, dean2011exponential}, rather than to a lattice.

A typical differential-difference equation on the plane is
\begin{equation}
\label{GenDiffEq}
\fr{\pd u}{\pd t} = \Delta u - \epsilon^2F(u;r),
\end{equation}
where $u\equiv u(x,y,t)$ for $(x,y,t) \in \mathbb{Z}^2\times[0,\infty)$. Here $F(u;r)$ is some nonlinear function of $u$ incorporating a bifurcation parameter $r$, which we assume to be bistable, allowing two stable states, and $\epsilon$ is some scaling which we will take to be small in our asymptotic analysis. The difference operator $\Delta$ comprises the nearest-neighbour stencil
\begin{equation}
\label{Delta}
\Delta u(x,y,t) := u(x+1,y,t) + u(x-1,y,t) + u(x,y+1,t) + u(x,y-1,t) - 4u(x,y,t).
\end{equation}
Although we have here given a specific $\Delta$, our methods are also applicable to a reasonably general class of difference operators. (\ref{GenDiffEq}) is the discrete analogue of the reaction-diffusion equation
\begin{equation}
\label{RDEq}
\fr{\pd u}{\pd T} = \lb( \fr{\pd^2}{\pd X^2} + \fr{\pd^2}{\pd Y^2} \rb) u - F(u;r),
\end{equation}
where we define the slow variables $(X,Y,T) = (\epsilon x,\epsilon y,\epsilon^2 t)\in\mathbb{R}^2\times[0,\infty)$; in other words, (\ref{GenDiffEq}) is the discrete approximation of (\ref{RDEq}) using second-order finite differences with a mesh spacing of $\epsilon$. Note, however, that (\ref{RDEq}) is invariant under arbitrary rotations in the plane, while (\ref{GenDiffEq}) is not. 

\begin{figure}[th!]
\centering
\includegraphics[trim = 0 0.7cm 1cm 0, clip, width=0.333\textwidth]{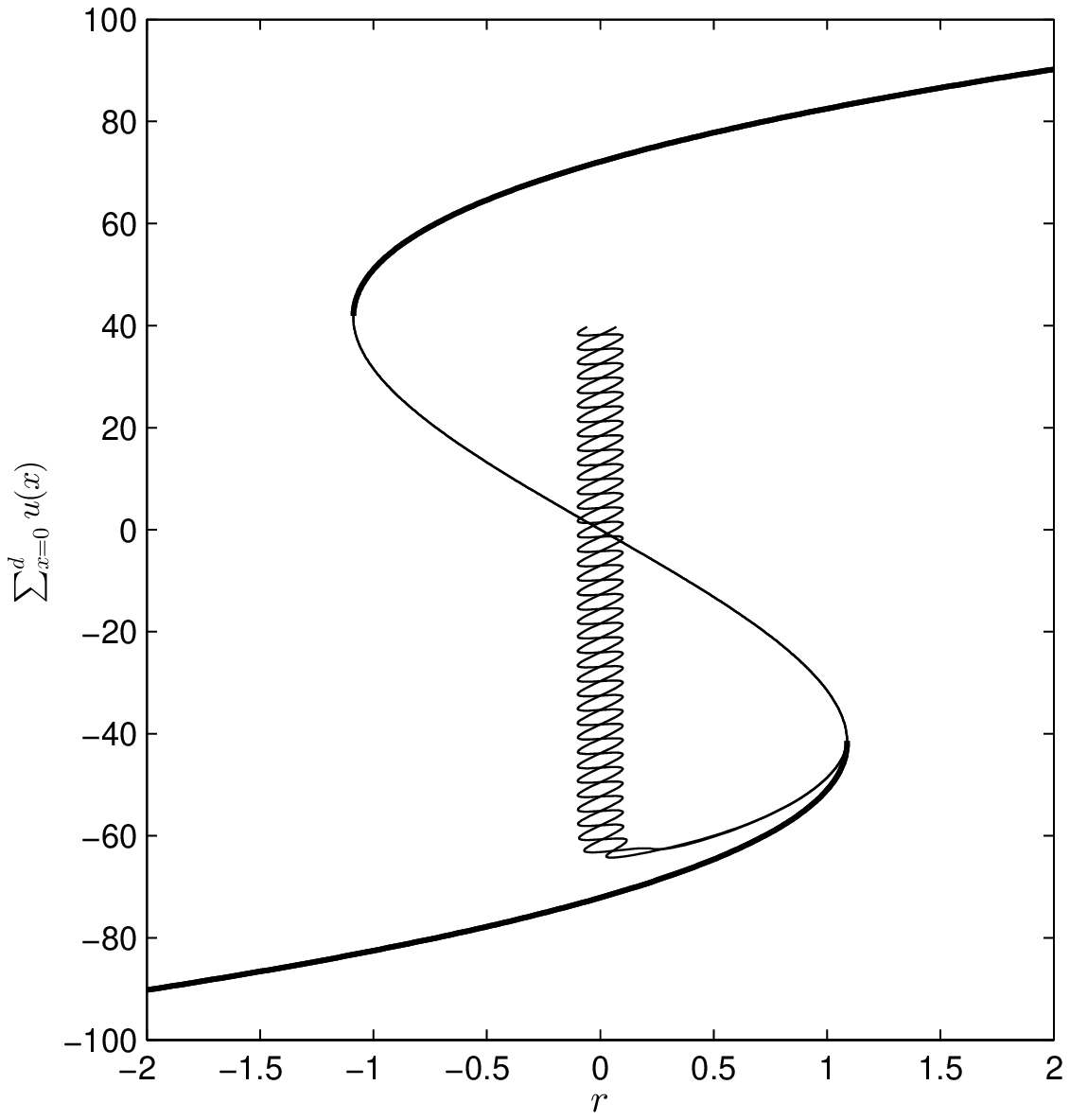}\includegraphics[trim = 0 0.7cm 1cm 0, clip, width=0.333\textwidth]{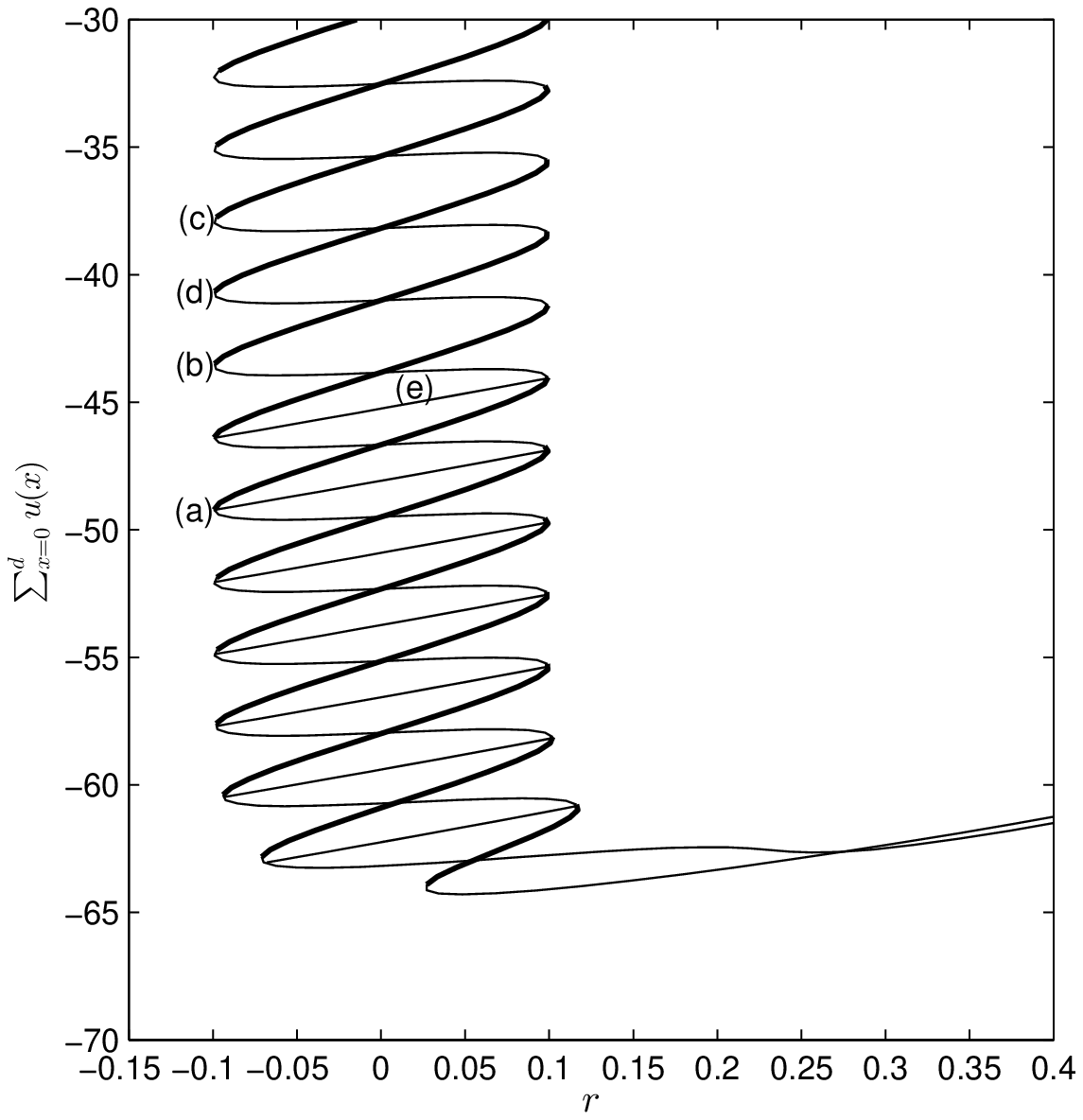}\includegraphics[trim = 0 1.1cm 1cm 0, clip, width=0.334\textwidth, keepaspectratio]{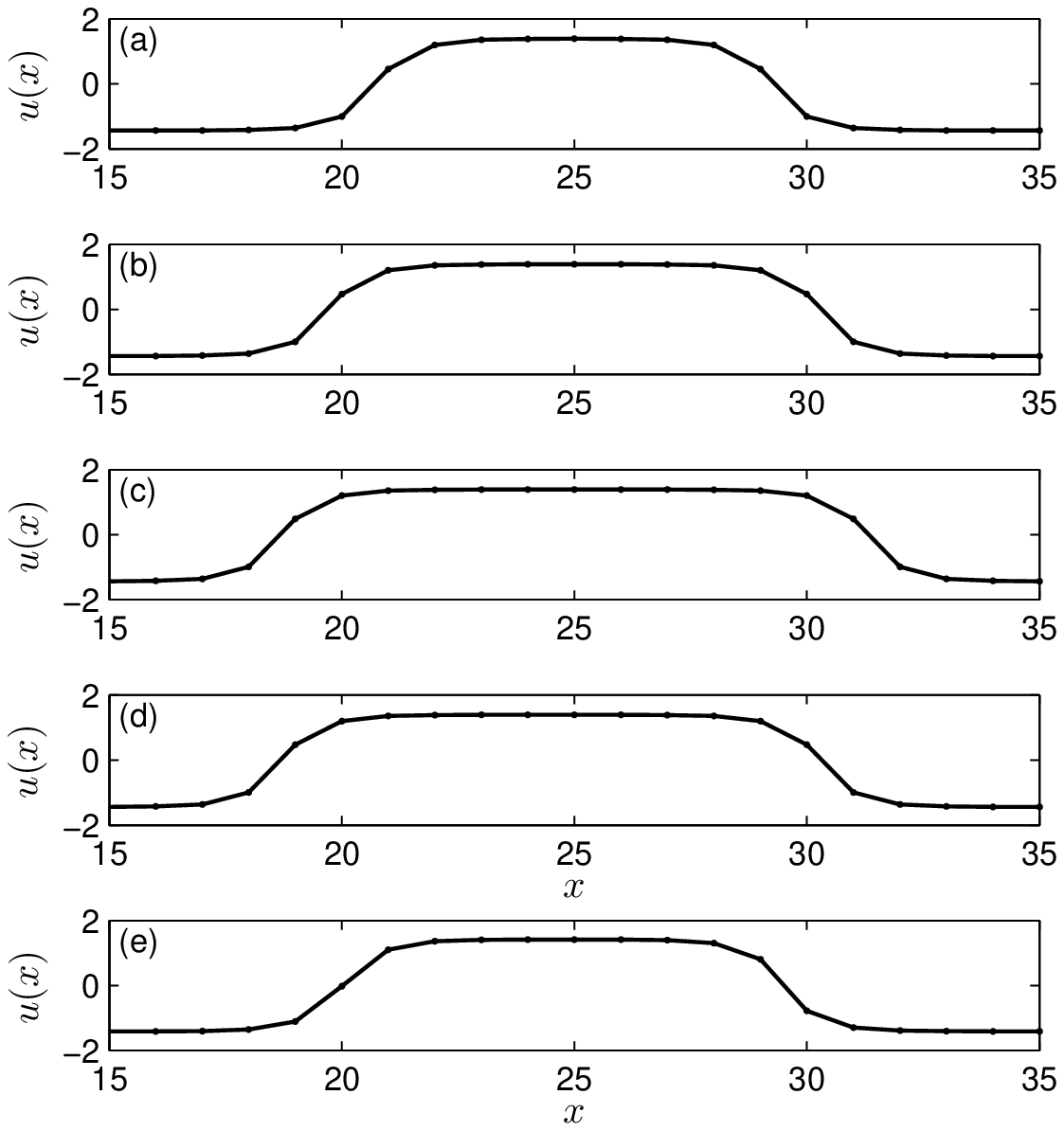}
\caption{Numerical solutions of (\ref{GenDiffEq}) with $u \equiv u(x)$ and $\epsilon^2F(u;r) = -r - 2u + u^3$ on the domain $x \in [0,d]$ with $d=50$. Left: bifurcation diagram showing the double saddle-node bifurcation of the constant solution and the snaking of the symmetric localised solutions within the region of bistability. For clarity, we omit the asymmetric `rung' solution branches in this panel. Centre: a zoomed-in view of the pinning region, rungs included. Thick (thin) lines indicate stable (unstable) solutions; we do not show stability of the snaking curves in the left-hand panel. Right: example solutions, zoomed in to the range $x \in [15,35]$. Labels indicate the position of each solution in the snaking diagram shown in the centre panel. (a)-(c) are site-centred solutions, (d) a bond-centred solution and (e) an asymmetric `rung' solution.}
\label{Pic_Const3ExampleBifDiag}
\end{figure}

\begin{figure}[th!]
\centering
\includegraphics[trim = 0 0.7cm 1cm 0, clip, width=0.333\textwidth]{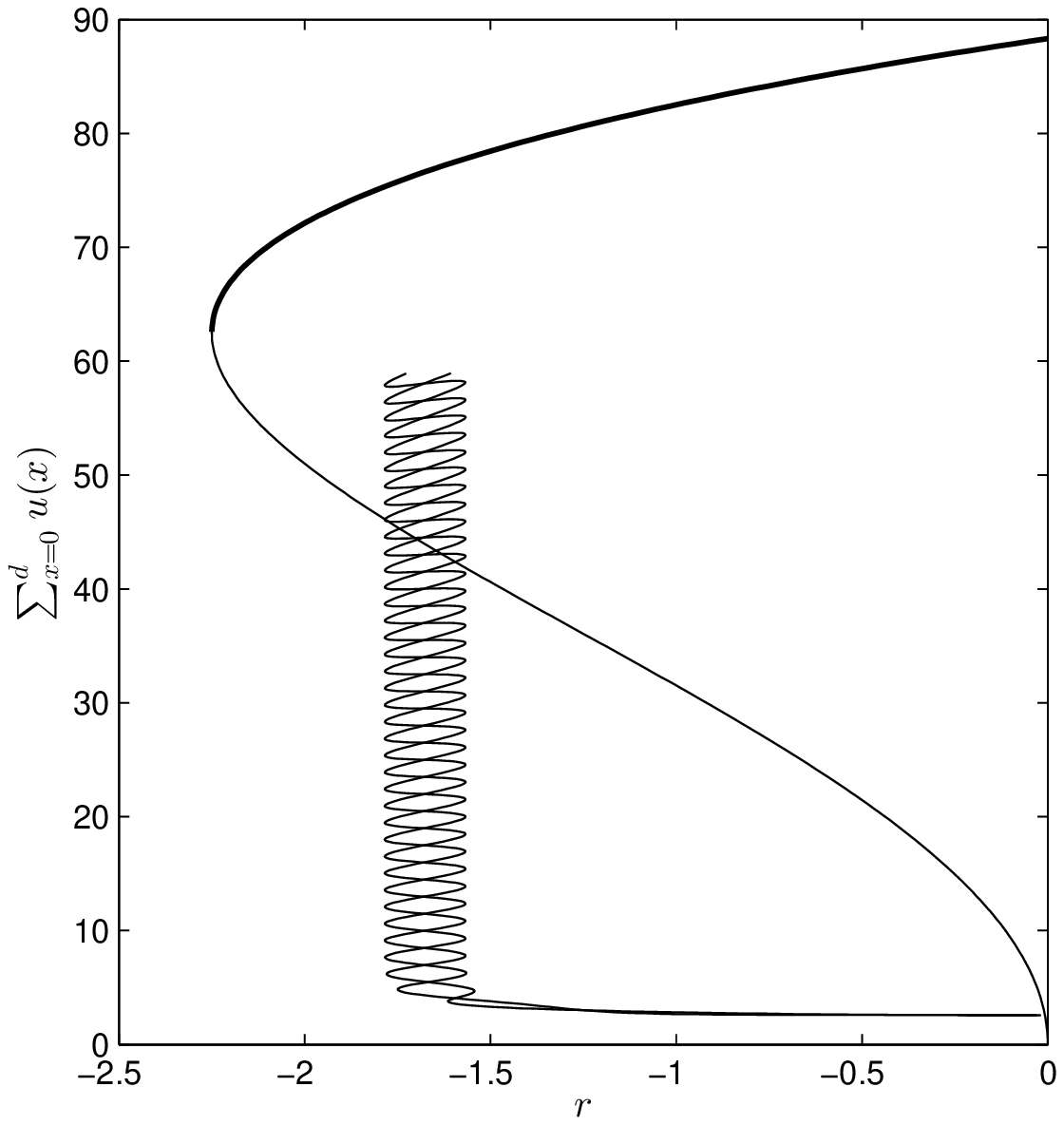}\includegraphics[trim = 0 0.7cm 1cm 0, clip, width=0.333\textwidth]{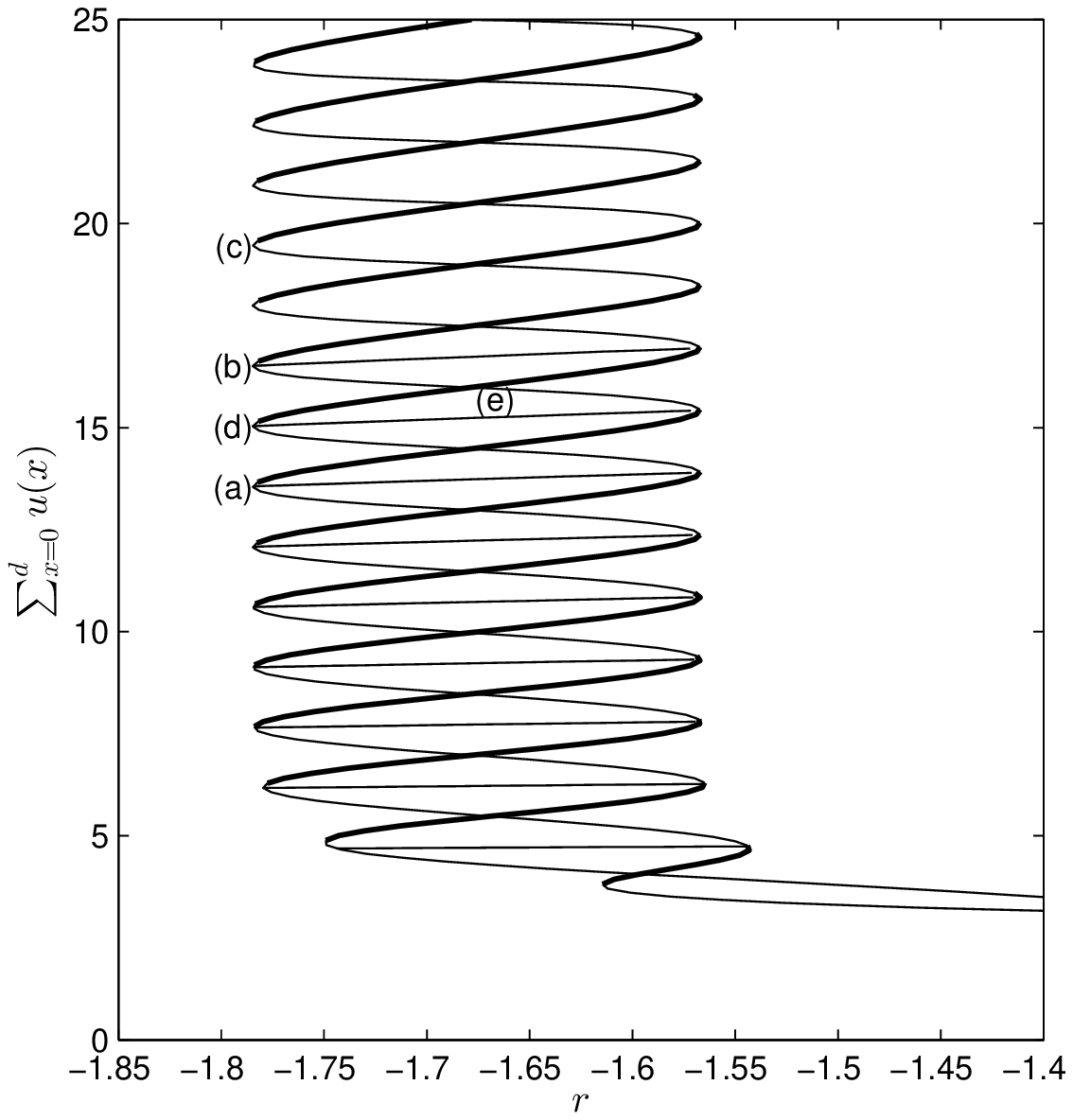}\includegraphics[trim = 0 1.1cm 1cm 0, clip, width=0.334\textwidth]{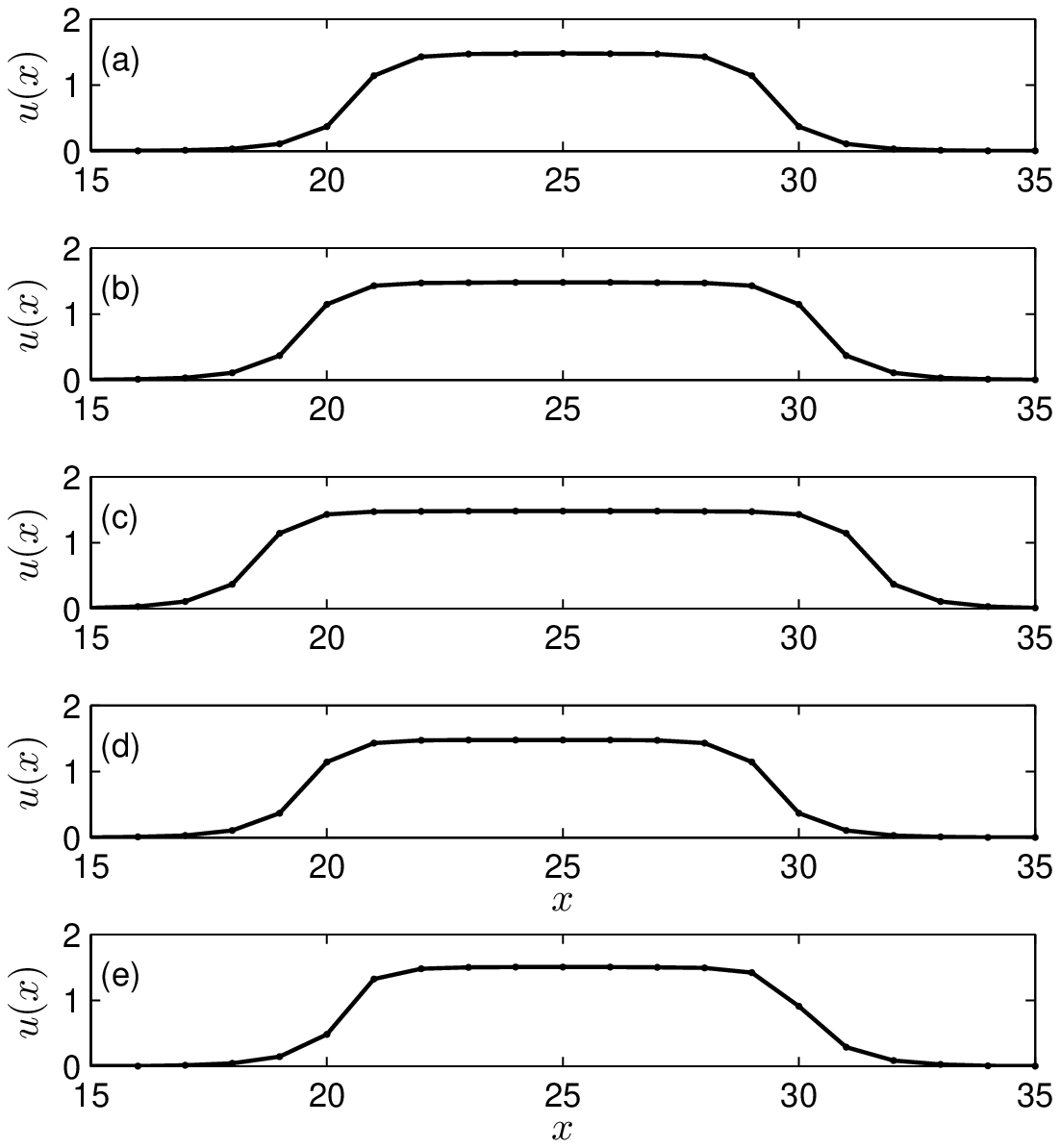}
\caption{Numerical solutions of (\ref{GenDiffEq}) with $u \equiv u(x)$ and $\epsilon^2F(u;r) = -ru - 3u^3 + u^5$ on the domain $x \in [0,d]$ with $d=50$. Left: bifurcation diagram showing the subcritical pitchfork and subsequent saddle-node bifurcation of the constant solution and the snaking of the symmetric localised solutions within the region of bistability. For clarity, we omit the asymmetric `rung' solution branches in this panel. Centre: a zoomed-in view of the pinning region, rungs included. Thick (thin) lines indicate stable (unstable) solutions; we do not show stability of the snaking curves in the left-hand panel. Right: example solutions, zoomed in to the range $x \in [15,35]$. Labels indicate the position of each solution in the snaking diagram shown in the centre panel. (a)-(c) are site-centred solutions, (d) a bond-centred solution and (e) an asymmetric `rung' solution.}
\label{Pic_35ExampleBifDiag}
\end{figure}

The formulation (\ref{GenDiffEq}) corresponds to scaling the system close to bifurcation, an approach common to nonlinear dynamical treatments of continuous systems (cf. \cite{burke2007homoclinic, chapman2009exponential, dean2011exponential} for snaking examples, and \cite{cross1993pattern} for a comprehensive review of others). This is equivalent to the limit of small mesh spacing, but we shall continue to use the language of nonlinear dynamics in order to facilitate comparison with other work on snaking. We note that a specific choice of nonlinearity may require rescaling of $u$ and $r$ before the system is in the form (\ref{GenDiffEq}), as we shall see in Section \ref{Sec_Examples}, where we apply our general results to the specific examples of figures \ref{Pic_Const3ExampleBifDiag} and \ref{Pic_35ExampleBifDiag}.

Two different snaking scenarios for (\ref{GenDiffEq}) with $u \equiv u(x)$ in the periodic domain $x \in [0,d]$ with $d=50$ are shown in figures \ref{Pic_Const3ExampleBifDiag} and \ref{Pic_35ExampleBifDiag}. Figure \ref{Pic_Const3ExampleBifDiag} is the result of setting $F = -r - 2u + u^3$; here bistability is the product of a pair of saddle-node bifurcations which together form an S-shaped solution curve. An alternative scenario is seen in figure \ref{Pic_35ExampleBifDiag}, which shows the bifurcation diagram and example solutions for $F = -ru - 3u^3 + u^5$. In this case bistability is due to a subcritical pitchfork bifurcation followed by a saddle-node, analogous to the cubic-quintic Swift-Hohenberg equation \cite{burke2007homoclinic, dean2011exponential}. The solution measure used is $\sum_{x=0}^d u(x)$; although not a norm in the strict sense, this choice is motivated by the desire that the visual representation of each loop of the snake is distinct in the first example. We remark that the second example is invariant under the reflection $u \rightarrow -u$, and so figure \ref{Pic_35ExampleBifDiag} can be reflected in the $r$-axis, while the first is invariant under the rotation in phase-space $(u,r) \rightarrow (-u,-r)$, and so there exists a second set of snaking curves emerging near the upper saddle-node in figure \ref{Pic_Const3ExampleBifDiag}. In an infinite domain the snaking curves continue indefinitely as the localised patch grows; in a finite domain, when the fronts near the domain boundaries the snaking curves leave the pinning region and reconnect to the constant solution branch \cite{bergeon2008eckhaus, dawes2009modulated}. We will apply the analytical results derived in the present work to both these examples in Section \ref{Sec_Examples}.

In the present work, we shall analyse the discrete snaking phenomenon via the method of exponential asymptotics \cite{adams2003beyond, berry1989uniform, oldedaalhuis1995stokes}. In a nutshell, this involves the calculation of exponentially small terms hidden beyond all algebraic orders in a divergent asymptotic expansion. If the expansion is truncated optimally, that is, after the least term, the resultant remainder is exponentially small. Careful analysis then indicates that the maximal change in the remainder occurs at certain lines in the complex plane emanating from singularities of the leading-order solution---Stokes lines. This rapid change manifests itself as the variation from zero to non-zero of the coefficient of a complementary function to the remainder equation; the variation is confined to a narrow layer around the Stokes line (the Stokes layer), and usually takes the form of an error function \cite{berry1989uniform, chapman1996non, oldedaalhuis1995stokes}. Although exponentially subdominant to the leading-order solution when `switched on' in this manner, the remainder is often of profound importance to the solution as a whole. Applications of the method include uniformly valid asymptotic approximations of integrals \cite{berry1989uniform, oldedaalhuis1995stokes}, travelling waves \cite{adams2003beyond}, flow past submerged bodies \cite{chapman2002exponential, chapman2006exponential, lustri2012free}, shock formation \cite{chapman2007shock} and waves formed in the wake of ships \cite{trinh2011waveless}.

In the current context, we shall find that the remainder is exponentially growing as the spatial variable tends to $\pm\infty$; eliminating such unbounded terms fixes the origin of a stationary front with respect to the lattice. This is precisely the pinning mechanism by which the snakes-and-ladders bifurcation diagram is generated. Furthermore, following previous work on snaking in the (continuous) Swift-Hohenberg equation \cite{chapman2009exponential, dean2011exponential}, the inclusion of an exponentially small deviation from the Maxwell point leads to a relation between the physical origin of the front and the distance from the Maxwell point in parameter space, which can only be satisfied within a certain exponentially small region---the pinning region. Armed with the full asymptotic expansion for a stationary front, we are able to match back-to-back fronts. The resultant matching conditions provide a set of formulae which fully describe the snakes-and-ladders bifurcation structure of the pinning region. Previous work by King and Chapman \cite{king2001asymptotics} described the pinning of fronts in a purely one-dimensional system; however, they did not study snaking explicitly, nor did they consider the rotation of solutions into the plane. In addition, we include a more general parameter-dependence in our lattice equation (\ref{GenDiffEq}). Thus the current calculation represents a significant extension of that work.

We remark that a beyond-all-orders analysis of discrete fronts has been carried out by Hwang \emph{et al.} in \cite{hwang2011solitary} using a different method than that employed here or in \cite{king2001asymptotics}, following \cite{yang1997traveling, yang1997asymmetric}; see also work by Keener \cite{keener2000homogenization}. The two approaches also differ in that discrete effects in that work were modelled by allowing a coefficient of the system to vary periodically in space, rather than through a difference operator as is the case here. We expect our method to be equally applicable to both means of modelling discreteness. An incomplete analysis \cite{clerc2011continuous} of discrete snaking phenomena has also been performed in a non-autonomous system similar to that in \cite{hwang2011solitary}; this fails to fully describe the snakes-and-ladders bifurcation because it does not consider exponentially small effects. Moreover, none of these studies considered solutions rotated into the plane; we believe the present work to be the first full asymptotic description of orientation-dependent pinning of stationary fronts, and of the resultant homoclinic snaking of localised solutions. We note that orientation-dependent pinning has been studied by several authors from a dynamical systems perspective \cite{hoffman2010universality, hupkes2011propagation, mallet2001crystallographic}, but remains an open problem. The present work complements these existing results, allowing us to observe the pinning mechanism explicitly and to derive an asymptotic relationship between the orientation of the front and its pinning region.

We begin in Section \ref{Sec_Rotation} by discussing the effects of rotating a one-dimensional localised solution into the plane, before defining some properties of the leading-order front in Section \ref{Sec_Preliminaries}, after which we perform a partial analysis of the remainder of the truncated asymptotic expansion in Section \ref{Sec_RemEqn}. We then calculate the late terms in the expansion in Section \ref{Sec_LateTerms}, allowing us to calculate the remainder in full in Section \ref{Sec_Truncation}. This leads to a formula for the snaking width in Section \ref{Sec_SnakeWidth} and an asymptotic description of the snaking bifurcation diagram in Section \ref{Sec_SnakeBifEqs}. Section \ref{Sec_Examples} sees the application of our general results to two specific examples, and comparison with numerical results. This is followed in Section \ref{Sec_Hex} by a brief discussion of how our results may be applied to a problem posed on a hexagonal lattice. We conclude in Section \ref{Sec_Conc}.


\section{Rotation into the plane}
\label{Sec_Rotation}

The solutions to (\ref{GenDiffEq}) shown in figures \ref{Pic_Const3ExampleBifDiag} and \ref{Pic_35ExampleBifDiag} are functions of $x$ only. Matters become somewhat more complicated when such one-dimensional solutions are rotated into the plane. Because the continuous analogue (\ref{RDEq}) of (\ref{GenDiffEq}) is invariant under rotations, an arbitrarily rotated solution remains a solution. This is not so in a discrete problem. Consider a one-dimensional solution of (\ref{GenDiffEq}); in order to incorporate arbitrary orientation with respect to the lattice we define
\begin{equation}
\label{z}
z = x\cos\psi + y\sin\psi, \qquad \psi\in[0,2\pi),
\end{equation}
and write $u\equiv u(z,t)$. The angle $\psi$ is measured anticlockwise from the positive $x$-axis. The difference operator (\ref{Delta}) is therefore rendered
\begin{equation}
\label{Delta_z}
\Delta u(z,t) = u(z+\cos\psi,t) + u(z-\cos\psi,t) + u(z+\sin\psi,t) + u(z-\sin\psi,t) - 4u(z,t).
\end{equation}
Thus $\psi$ retains an explicit presence in the rotated, one-dimensional version of (\ref{GenDiffEq}), in contrast to its continuous analogue (\ref{RDEq}), and solutions depend parametrically on their orientation $\psi$.

In particular, fronts cannot pin to the lattice if the tangent of the angle of orientation is irrational, for reasons we shall elucidate presently. In such a case the width of the pinning region collapses to zero. This phenomenon has been the subject of much study from a dynamical systems point of view \cite{hoffman2010universality, hupkes2011propagation, mallet2001crystallographic}. The present calculation complements the more general results derived in such work, allowing us to observe explicitly the pinning mechanism and the vanishing of the pinning region at irrational orientations, and to derive asymptotic formulae for the width of the pinning region and the resultant snaking bifurcation diagram.

\begin{figure}[th!]
\centering
\includegraphics[trim=1cm 2cm 1cm 2cm, clip, width = 0.6\textwidth]{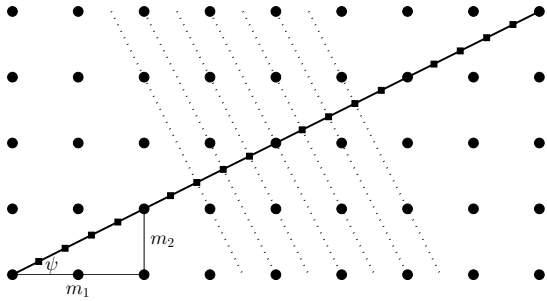}
\caption{The effective (one-dimensional) lattice $(m_1x + m_2y)(m_1^2+m_2^2)^{-1/2}$ with rational $\tan\psi = m_2/m_1$ as defined in (\ref{Rational_tanpsi}), superimposed onto the actual (two-dimensional) lattice $(x,y) \in \mathbb{Z}^2$. Actual lattice points are represented by circles; effective ones by squares. The independent variable $z$ varies in the direction of the solid line; $z$ is constant along each dotted line, which also indicate the correspondence of effective to actual lattice points.}
\label{Pic_ExplainLattice}
\end{figure}

The importance of the rationality of $\psi$ can be understood by considering the spatial domain of (\ref{Delta_z}), which is the countable set
\begin{equation}
\Psi := \lb\{ x\cos\psi + y\sin\psi\ |\ (x,y)\in\mathbb{Z}^2 \rb\}.
\end{equation}
We also define the extended set of rational numbers
\begin{equation}
\mathbb{Q}_\infty := \mathbb{Q} \cup \{\pm\infty\},
\end{equation}
assigning $\tan(\pm\fr{\pi}{2}) = \pm\infty$. If $\tan\psi\in\mathbb{Q}_\infty$, then we can set
\begin{equation}
\label{Rational_tanpsi}
\tan\psi = \fr{m_2}{m_1}, \qquad (m_1,m_2) \in \mathbb{Z}^2\backslash\{(0,0)\}, \qquad \mathrm{gcd}(|m_1|,|m_2|)=1,
\end{equation}
without loss of generality, in which case
\begin{equation}
\Psi = \lb\{ \lb. \fr{m_1x + m_2y}{\sqrt{m_1^2+m_2^2}} \ \rb| \ (x,y)\in\mathbb{Z}^2 \rb\}.
\end{equation}
Since $m_1x+m_2y$ is an integer, $\Psi$ describes a one-dimensional lattice with a well-defined lattice spacing of $(m_1^2+m_2^2)^{-1/2}$. Infinitely many points of the actual (two-dimensional) lattice are mapped to each point of this effective lattice $\Psi$, with the value of $u$ at the actual lattice point $(x,y)$ being equal to that of $u$ at the $(m_1x + m_2y)$th effective lattice point, as indicated in figure \ref{Pic_ExplainLattice}. 

In contrast, if $\tan\psi \notin \mathbb{Q}_\infty$, i.e. is irrational and finite, then $\Psi$ is a dense (and countably infinite) set. As a consequence, any point on the real line is arbitrarily close to a point in $\Psi$. Thus there is no well-defined lattice spacing for irrational $\tan\psi$, without which a front cannot pin to the lattice. We therefore expect one-dimensional snaking to occur only when $\tan\psi$ is rational or infinite, i.e. when $\tan\psi \in \mathbb{Q}_\infty$. 

\section{Setting up the beyond-all-orders calculation}
\label{Sec_Preliminaries}

Focusing now on one-dimensional solutions rotated into the plane, we write $u \equiv u(z,t)$ in (\ref{GenDiffEq}), yielding
\begin{equation}
\label{GenDiffEq_z}
\fr{\pd u(z,t)}{\pd t} = \Delta u(z,t)  - \epsilon^2 F(u(z,t);r),
\end{equation}
where $z$ is defined as in (\ref{z}) and $\Delta u(z,t)$ is given in (\ref{Delta_z}). Note (\ref{GenDiffEq_z}) has continuous analogue
\begin{equation}
\label{RDEq_Z}
\fr{\pd u}{\pd T} = \fr{\pd^2u}{\pd Z^2} - F(u;r),
\end{equation}
where
\begin{equation}
\label{Z}
Z = \epsilon \lb( z - z_0 \rb).
\end{equation}
We include the (constant) origin $z_0$ in order to enable the derivation of the pinning mechanism later on; although the continuous equation (\ref{RDEq_Z}) is invariant under translations in space, the discrete equation (\ref{GenDiffEq_z}) is invariant only under translations by integer multiples of the effective lattice spacing $(m_1^2+m_2^2)^{-1/2}$. Fixing $z_0$ therefore corresponds to the pinning of a front to the lattice. This is analogous to fixing the phase of the spatial oscillations in a continuous system \cite{burke2007homoclinic, dean2011exponential, kozyreff2006asymptotics, pomeau1986front}. We will expand upon the significance of $z_0$ and the means by which it can be determined presently. 

We define $u_c$ to be a constant solution of (\ref{GenDiffEq_z}), so that $F(u_c;r)=0$. We can investigate the stability of $u_c$ by setting $u = u_c + \hat{u} e^{\sigma t + ikz}$ in (\ref{GenDiffEq_z}), where $k \in [0,2\pi)$,  and linearising with $|\hat{u}| \ll 1$. Thus we obtain the growth rate equation
\begin{equation}
\label{GrowthRate}
\sigma \sim -2\lb[ 2 - \cos(k\cos\psi) - \cos(k\sin\psi) \rb] - \epsilon^2 F_u\lb( u_c ; r \rb),
\end{equation}
where the subscript $u$ denotes the first derivative of $F(u;r)$ with respect to $u$. Hence $u_c$ is linearly stable provided $F_u\lb(u_c;r\rb)>0$, but becomes linearly unstable to perturbations with small wavenumber $k$ as $F_u\lb(u_c;r\rb)$ becomes negative. We can therefore describe the dynamics of (\ref{GenDiffEq_z}) close to bifurcation using the double limit $\epsilon \rightarrow 0$ and $k \rightarrow 0$, under which (\ref{GrowthRate}) becomes
\begin{equation}
\sigma \sim - k^2 - \epsilon^2 F_u\lb(u_c;r\rb).
\end{equation}
This suggests that (\ref{GenDiffEq_z}) evolves with the slow scales $(Z,T) = (\epsilon (z-z_0),\epsilon^2t)$ as $\epsilon\rightarrow0$, precisely the independent variables of the continuous analogue (\ref{RDEq_Z}). 

Writing $u \equiv u(Z,T)$, we have $u(z\pm\cos\psi,t) \rightarrow u(Z\pm\epsilon\cos\psi,T)$ and $u(z\pm\sin\psi,t) \rightarrow u(Z\pm\epsilon\sin\psi,T)$. The small-$\epsilon$ limit can therefore be exploited to expand the difference operator $\Delta u$ (\ref{Delta_z}) in powers of $\epsilon$ using Taylor's theorem, rendering (\ref{GenDiffEq_z}) as
\begin{equation}
\label{TaylorExpandDiffEq}
\epsilon^2 \fr{\pd u}{\pd T} = 2 \sum_{p=1}^{\infty} \epsilon^{2p}\fr{\cos^{2p}\psi+\sin^{2p}\psi}{(2p)!} \fr{\pd^{2p}u}{\pd Z^{2p}} - \epsilon^2F(u;r).
\end{equation}
Note that only even powers of $\epsilon$ are present. The leading-order approximation to (\ref{GenDiffEq_z}) is therefore simply the continuous analogue (\ref{RDEq_Z}) of (\ref{GenDiffEq_z}).  

As we are interested in stationary solutions of (\ref{GenDiffEq_z}), we now write $u \equiv u(Z)$ and expand in powers of $\epsilon^2$ as
\begin{equation}
\label{Expandu}
u(Z) \sim \sum_{n=0}^{N-1} \epsilon^{2n}u_n(Z) + R_N(z,Z).
\end{equation}
Note that we have truncated the expansion after $N$ terms; this is  because it is divergent. As (\ref{TaylorExpandDiffEq}) is a singular perturbation problem, in that successively higher derivatives contribute at successive orders in $\epsilon$, the $n$th term in the expansion depends upon the derivatives of the previous terms. Therefore, if the leading-order solution has singularities (in the present context, these are bounded away from the real line), the resultant asymptotic expansion is divergent in the form of a factorial over a power, and must be truncated \cite{adams2003beyond, berry1989uniform, oldedaalhuis1995stokes}. In an abuse of notation, we have retained $z$-dependence in the remainder $R_N$, for reasons to be explained in Section \ref{Sec_Truncation}; we remark for now that the pinning mechanism will manifest as an interplay between the fast scale $z$ of the lattice and the slow scale $Z$ of the front. If we choose the point of truncation optimally by truncating at the point at which the expansion begins to diverge, the remainder will be exponentially small in $\epsilon$, thus allowing us to investigate exponentially small effects.

Now, the leading-order (steady) contribution to (\ref{TaylorExpandDiffEq}) is
\begin{equation}
\label{FrontEq}
0 = \fr{\od^2 u_0}{\od Z^2} - F(u_0;r),
\end{equation}
which is of course simply the steady version of (\ref{RDEq_Z}). We shall assume that $u_0(Z)$ takes the form of a stationary front, and hence impose the boundary conditions
\begin{equation}
\label{BCs}
u_0 \rightarrow u_\pm \ \mathrm{as} \ Z \rightarrow \pm\infty,
\end{equation}
where $u_\pm$ are stable, constant solutions of (\ref{GenDiffEq}) and therefore satisfy
\begin{equation}
\label{FofSteadyStateConditions}
F(u_\pm;r) = 0,\qquad F_u(u_\pm;r)>0.
\end{equation}
We shall also assume, without loss of generality, that $u_-<u_+$, as the front with opposite orientation in the plane is simply given by the rotation $\psi \rightarrow \psi+\pi$.

In order to investigate the phenomenon of homoclinic snaking, we shall restrict our attention to the class of functions $F(u;r)$ where front solutions to (\ref{FrontEq}), connecting the two constant solutions $u_\pm$, exist only at a particular value of the bifurcation parameter, $r=r_M$ say. This is the Maxwell point, the point in parameter space at which travelling waves connecting $u_-$ to $u_+$ have zero velocity. Because (\ref{FrontEq}) can be integrated once, the constant solutions $u_\pm$ must also satisfy its first integral; thus $r_M$ must satisfy
\begin{equation}
\label{rM_integral}
\int_{u_-}^{u_+} F(v;r_M) \od v = 0,
\end{equation}
as well as $F(u_\pm;r_M)=0$. Note that these conditions form a system of three algebraic equations in the three unknowns $u_\pm$ and $r_M$, providing a means of determining the Maxwell point. From this we might (erroneously) infer that stationary fronts exist only at the Maxwell point, in direct contradiction of numerical results showing homoclinic snaking within a well-defined region of parameter space centred on the Maxwell point (e.g. figures \ref{Pic_Const3ExampleBifDiag} and \ref{Pic_35ExampleBifDiag}). A standard continuum approximation cannot reconstruct such behaviour, as snaking is confined to an exponentially small distance from $r_M$; such scales are indistinguishable using techniques based solely on algebraic powers of $\epsilon$. Thus we must employ exponential asymptotics in order to capture the snaking phenomenon.

We remark that for some choices of $F$ the integral condition (\ref{rM_integral}) is satisfied without the need to impose a specific value of $r$. For example, if we choose $F = r\sin u$, $u_+=2\pi$ and $u_-=0$ then (\ref{FofSteadyStateConditions}) and (\ref{rM_integral}) hold for all $r>0$. For such an $F$ snaking does not occur, as there is no Maxwell point and stationary front solutions to the leading-order approximation may be found across an $\Or(1)$ interval of $r$-values rather than at a specific point. However, in such a case fronts still pin to the lattice by selecting an origin $z_0$, and so much of the following calculation remains relevant.

In order to incorporate exponentially small deviations from the Maxwell point into subsequent calculations, we write $r = r_M + \delta r$ and expand $F(u;r)$ around $r_M$ as
\begin{equation}
\label{FatMP}
F(u;r_M+\delta r) = F_M(u) + \delta r F_{r,M}(u) + \Or(\delta r^2),
\end{equation}
where we define
\begin{equation}
F_M(u) := F(u;r_M), \qquad F_{r,M}(u) := \lb. \fr{\pd F}{\pd r} \rb|_{r=r_M}.
\end{equation}
We assume that $F_{r,M} \neq 0$ for simplicity, but note that the present work may in principle be extended to choices of $F$ whose first non-zero derivative with respect to $r$ at the Maxwell point is of higher order. $\delta r$ is thus the bifurcation parameter we shall use to describe the snaking bifurcations; it will turn out to be exponentially small. In principle, one should also include further algebraic corrections to the Maxwell point by writing $r = r_M + \epsilon^2 r_2 + \cdots + \epsilon^{2N-2} r_{2N-2} + \delta r$; each of the $r_j$ can be fixed by successive solvability conditions at successive orders in $\epsilon^2$ (cf. the derivation of higher-order corrections to the Maxwell point in the Swift-Hohenberg equation in \cite{chapman2009exponential, dean2011exponential}). However, only the leading-order term $r_M$ and the exponentially small remainder $\delta r$ are important to the present calculation, so we shall not discuss such algebraic terms further.

We remark that the instability of $u_c$ to modes with small wavenumber is in contrast to the equivalent situation in, for example, the Swift-Hohenberg equation, in which the zero solution loses stability to modes with wavenumber $\pm1$ \cite{burke2007homoclinic}; such an instability is pattern-forming and produces a spatial structure to which fronts may pin. No such pattern-forming mechanism is present in the second-order equation (\ref{RDEq_Z}). Hence there is no spatial structure, and nothing for a front to pin to. Therefore the leading-order continuum approximation (\ref{RDEq_Z}) of (\ref{GenDiffEq_z}) does not exhibit snaking. However, snaking persists in numerical computations of (\ref{GenDiffEq_z}) even very close to bifurcation, indicating that the continuum approximation (\ref{RDEq_Z}) does not tell the whole story. This discrepancy can be resolved by incorporating higher-order effects in the asymptotic solution to (\ref{GenDiffEq_z}), in particular those which are exponentially small \cite{chapman2009exponential, dean2011exponential, king2001asymptotics}.

The present calculation is in some respects simpler than analogous work in the Swift-Hohenberg equation \cite{chapman2009exponential, dean2011exponential}. For instance, the appropriate method for studying continuous pattern formation near onset is that of multiple scales, rather than the relatively simpler continuum approximation employed in discrete problems. Moreover, the nonlinearities present in the Swift-Hohenberg equation lead to an ever-increasing number of harmonics $e^{kix}$ at successive orders in $\epsilon$, with the obvious consequence of an ever-increasing number of equations determining their coefficients. That said, the Taylor expansion of slow differences results in what is in effect an infinite-order differential equation, with successively higher derivatives contributing at successive orders in $\epsilon$, and so the current calculation is not without its own complexities. 


\section{The remainder equation}
\label{Sec_RemEqn}

Although we are not yet in a position to solve for the remainder, we are able to determine much information about it. The leading-order equation for $R_N$ is
\begin{equation}
\label{RNEqn_unknownForcing}
\Delta R_N(z,Z) - \epsilon^2 F_M'(u_0(Z))R_N(z,Z) \sim \epsilon^2 \delta r F_{r,M}(u_0(Z)) + \ \mathrm{forcing\ due\ to\ truncation},
\end{equation}
where
\begin{align}
\Delta R_N(z,Z) = &\ R_N(z+\cos\psi, Z+\epsilon\cos\psi) + R_N(z-\cos\psi, Z-\epsilon\cos\psi)
\nonumber\\
& {}  + R_N(z+\sin\psi, Z+\epsilon\sin\psi)+ R_N(z-\sin\psi, Z-\epsilon\sin\psi) - 4R_N(z,Z),
\end{align}
and the exact scalings of $R_N$ and $\delta r$, although exponentially small, are yet to be determined. The left-hand side of (\ref{RNEqn_unknownForcing}) is simply the linearisation of the steady version of (\ref{GenDiffEq_z}) around $u_0$; the first term on the right-hand side is due to the linearisation of $F(u;r)$ about the Maxwell point as in (\ref{FatMP}), while the second is the result of the truncation of the asymptotic series after $N$ terms in (\ref{Expandu}). The forcing due to truncation is at present unknown, since we do not yet possess an expression for the large-$n$ terms in (\ref{Expandu}). However, we are able at this point to derive the complementary functions of (\ref{RNEqn_unknownForcing}). It is these which will be switched on as Stokes lines are crossed. Furthermore, we can determine the forcing due to the deviation $\delta r$ from the Maxwell point, and see how the combination of this integral and the complementary functions lead to a solvability condition on the leading-order front. The derivation of the precise solvability condition requires the large-$n$ terms; these we calculate in Section \ref{Sec_LateTerms}, allowing us to evaluate the at present undetermined forcing in (\ref{RNEqn_unknownForcing}) and carry out the full beyond-all-orders calculation for $R_N$ in Section \ref{Sec_Truncation}.

As (\ref{RNEqn_unknownForcing}) is linear, and autonomous with regard to the fast scale $z$, we can look for a solution to the homogeneous equation of the form
\begin{equation}
R_N(z,Z) = e^{i\kappa z} S_N(Z) + \cc,
\end{equation}
for some eigenvalue $\kappa \in \mathbb{C}$, and Taylor expand the slow-scale differences in powers of $\epsilon$. After cancellation of the common factor $e^{i\kappa z}$, this results in
\begin{align}
\label{ExpandRemEq}
& 2 \lb[ \cos(\kappa\cos\psi) + \cos(\kappa\sin\psi) - 2 \rb] S_N + 2 i \epsilon \lb[ \cos\psi\sin(\kappa\cos\psi) + \sin\psi \sin(\kappa\sin\psi) \rb] S_N'
\nonumber\\
& {} + \epsilon^2 \lb[ \cos^2\psi\cos(\kappa\cos\psi) + \sin^2\psi\cos(\kappa\sin\psi) \rb]S_N'' - \epsilon^2 F_M'(u_0(Z)) S_N = \Or\lb(\epsilon^3S_N\rb).
\end{align}
Expanding $S_N$ as
\begin{equation}
S_N(Z) = S_{N,0}(Z) + \epsilon S_{N,1}(Z) + \epsilon^2 S_{N,2}(Z) + \cdots,
\end{equation}
then, if $S_{N,0}$ is to be non-zero, we obtain at $\Or(S_N)$ the condition
\begin{equation}
\label{kappa_cond_O1}
\cos(\kappa\cos\psi) + \cos(\kappa\sin\psi) - 2 = 0.
\end{equation}
Real solutions to (\ref{kappa_cond_O1}) are given by $\kappa\cos\psi=2M_1\pi$ and $\kappa\sin\psi = 2M_2\pi$, for any $(M_1,M_2)\in\mathbb{Z}^2$; these exist only when $\tan\psi \in \mathbb{Q}_\infty$. Hence there are no real, non-zero solutions to (\ref{kappa_cond_O1}) for irrational $\tan\psi$; however, in general there exist complex solutions to (\ref{kappa_cond_O1}). If $\tan\psi \in \mathbb{Q}_\infty$ and $\kappa \in \mathbb{R}$, we may therefore set
\begin{equation}
\label{tanphi_realkappa}
\begin{split}
\begin{array}{clcl}
& \cos\psi = \dfrac{m_1}{\sqrt{m_1^2+m_2^2}}, & \qquad & \sin\psi = \dfrac{m_2}{\sqrt{m_1^2+m_2^2}}, 
\\
\vspace{-1ex}
\\
& (m_1,m_2) \in \mathbb{Z}^2\backslash\{(0,0)\}, & \qquad & \mathrm{gcd}(|m_1|,|m_2|) = 1,
\end{array}
\end{split}
\end{equation}
without loss of generality. This then gives
\begin{equation}
\label{realkappa}
\kappa = 2M\pi\sqrt{m_1^2+m_2^2}, \qquad M\in\mathbb{Z}.
\end{equation}
Of particular note are the axes and principal diagonals, given by $\psi = \fr{k\pi}{4}$, $k \in \{0,1,\ldots,7\}$. These correspond to either $\cos\psi$ having unit modulus and $\sin\psi$ vanishing, or vice versa, or both $\cos\psi$ and $\sin\psi$ having modulus $1/\sqrt{2}$. In each of these eight instances, (\ref{realkappa}) describes all solutions to (\ref{kappa_cond_O1}), i.e. (\ref{kappa_cond_O1}) has no complex solutions when $\psi = \fr{k\pi}{4}$, $k \in \{0,1,\ldots,7\}$.

If $\psi$ and $\kappa$ respectively satisfy (\ref{tanphi_realkappa}) and (\ref{realkappa}), we find that $\Or(\epsilon S_N)$ terms also vanish in (\ref{ExpandRemEq}). Proceeding to $\Or(\epsilon^2 S_N)$, we then obtain
\begin{equation}
\label{SN0_eqn}
S_{N,0}''
 - F_M'(u_0(Z)) S_{N,0} = 0.
\end{equation}
As $u_0$ satisfies (\ref{FrontEq}), the complementary functions of (\ref{SN0_eqn}) are
\begin{equation}
\label{CompFunctions}
g(Z) := u_0'\lb( Z \rb),\qquad G(Z;\zeta) := u_0'\lb( Z \rb) \int_\zeta^{Z} \fr{1}{u_0'(t)^2} \od t,
\end{equation}
where $g(Z)$ can be found by noting that (\ref{SN0_eqn}) with $S_{N,0} = u_0'$ is simply the first derivative of (\ref{FrontEq}), after which $G(Z;\zeta)$ can readily be found using the method of reduction of order. The parameter $\zeta$ is a (complex) singularity of $u_0(Z)$, included to simplify subsequent calculations. Thus each real $\kappa$ provides a contribution to $R_N$ of the form $e^{i\kappa z} \lb( a_\kappa g + A_\kappa G \rb)$, for some constants $a_\kappa$ and $A_\kappa$.

We now turn our attention to complex, with $\Im(\kappa) \neq 0$, solutions of (\ref{kappa_cond_O1}), noting that there are no non-zero, purely imaginary solutions to (\ref{kappa_cond_O1}). Requiring that $\Or(\epsilon S_N)$ terms in (\ref{ExpandRemEq}) vanish, we must have either
\begin{equation}
\label{kappa_cond_Oeps}
\lb[ \cos\psi\sin(\kappa\cos\psi) + \sin\psi \sin(\kappa\sin\psi) \rb] = 0
\end{equation}
or
\begin{equation}
S_{N,0}' = 0.
\end{equation}
It can be shown that if $\kappa$ is complex and satisfies (\ref{kappa_cond_O1}) then it does not satisfy (\ref{kappa_cond_Oeps}); we defer this calculation to the appendix. Thus, if $\kappa$ is complex, we have $S_{N,0} = B_\kappa$, for some constant $B_\kappa$.

Finally, we may seek the particular integral of (\ref{RNEqn_unknownForcing}) due to the term involving $\delta r$. Setting $R_N(z,Z) = \delta r P(Z)$, we obtain at leading order
\begin{equation}
P'' - F_M'(u_0(Z))P = F_{r,M}(u_0(Z)),
\end{equation}
which can be solved using the method of variation of parameters to give
\begin{equation}
\label{PI_deltar}
P(Z) = u_0'\lb(Z\rb) \int^{Z} \fr{1}{u_0'(t)^2} \lb[ \int_{u_-}^{u_0(t)} F_{r,M}(v) \od v \rb] \od t.
\end{equation}
Combining the contributions for real and complex $\kappa$ and the particular integral, the leading-order solution to (\ref{RNEqn_unknownForcing}) for each $\zeta$ is thus
\begin{equation}
\label{RN_soln}
R_N(z,Z) \sim \delta r P(Z) + \sum_{\kappa\in\mathbb{R}} e^{i\kappa z} \lb[ a_\kappa g(Z) + A_\kappa G(Z;\zeta) \rb] + \sum_{\kappa \notin \mathbb{R}} e^{i\kappa z} B_\kappa,
\end{equation}
for arbitrary constants $a_\kappa$, $A_\kappa$ and $B_\kappa$. We emphasize that this solution does not account for the forcing due to truncation in (\ref{RNEqn_unknownForcing}); we remedy this in Section \ref{Sec_Truncation}. Crucially, the form of the solution (\ref{RN_soln}) is dependent upon the rationality of $\tan\psi$. If $\tan\psi \in \mathbb{Q}_\infty$, we can define $\psi$ as in (\ref{tanphi_realkappa}), in which case real $\kappa$ are given by (\ref{realkappa}). Recall that if $\psi = \fr{k\pi}{4}$ with $k \in \{0,1,\ldots,7\}$ then all solutions are purely real and the second summation in (\ref{RN_soln}) does not contribute; this is not the case for $\tan\psi \in \mathbb{Q}_\infty$ in general. On the other hand, if $\tan\psi \notin \mathbb{Q}_\infty$, then the only real solution to (\ref{kappa_cond_O1}) is $\kappa=0$, and so the first summation comprises only the contribution from this one value of $\kappa$. 

We also note that if $\tan\psi$ is rational or infinite then $\kappa z = 2M\pi( m_1x + m_2y)$; hence $e^{i\kappa z} = 1$ on lattice points. However, writing $R_N$ in the form (\ref{RN_soln}) will prove to be useful later on, when we come to evaluate the effects of the as yet unknown forcing in (\ref{RNEqn_unknownForcing}), and so we shall continue to write $e^{i\kappa z}$ even when $\kappa \in \mathbb{R}$.

\subsection{The form of the solvability condition}
\label{Sec_FormSolvabilityCondition}

We are now able to deduce the source of the beyond-all-orders solvability condition which determines the origin of the leading-order front. Linearising (\ref{FrontEq}) around the constant solutions $u_0(Z)\equiv u_\pm$, we can find expressions for $u_0$ in the far-fields, namely 
\begin{equation}
\label{FarFieldu0}
u_0 \sim u_\pm \mp D_\pm e^{\mp \alpha_\pm Z} \ \mathrm{as} \ Z \rightarrow \pm \infty,
\end{equation}
where
\begin{equation}
\label{alpha}
\alpha_\pm := +\sqrt{ F_M'(u_\pm) } > 0
\end{equation}
and $D_\pm > 0$. Therefore the complementary functions defined in (\ref{CompFunctions}) are given by
\begin{align}
g & \sim \alpha_\pm D_\pm e^{\mp \alpha_\pm Z}, 
\\
\label{FarFieldG}
G & \sim \pm \fr{1}{2\alpha_\pm^2 D_\pm} e^{\pm\alpha_\pm Z}
\end{align}
as $Z \rightarrow \pm\infty$, and the particular integral by
\begin{equation}
\label{FarFieldP}
P \sim \lb\{ \begin{array}{lcl}
\lb( 2 \alpha_+^2D_+ \rb) ^{-1} \lb( \int_{u_-}^{u_+} F_{r,M}(v) \od v \rb) e^{\alpha_+ Z}, & & Z \rightarrow \infty,
\\
 - \alpha_-^{-2} F_{r,M}(u_-), & & Z \rightarrow -\infty.
\end{array} \rb. 
\end{equation}
Thus there appears exponential growth in $G$ as $Z \rightarrow \pm\infty$, and in $P$ as $Z\rightarrow\infty$. We shall see that the as yet undetermined forcing in (\ref{RNEqn_unknownForcing}), due to truncation of the divergent series (\ref{Expandu}), results in a non-zero multiple of $G$ being present as $Z \rightarrow \infty$, but not as $Z \rightarrow -\infty$. This is an example of Stokes' phenomenon \cite{berry1989uniform, oldedaalhuis1995stokes}, in which $G$ is switched on as Stokes lines are crossed; we shall show this explicitly in Section \ref{Sec_Truncation}. Thus terms which grow exponentially as $Z \rightarrow \infty$ appear in the remainder from two sources: a complementary function switched on according to Stokes' phenomenon, and the particular integral due to the deviation from the Maxwell point. However, we shall show in Section \ref{Sec_SnakeWidth} that for certain values of the origin $z_0$ of $u_0$, defined in terms of $\delta r$, the coefficient of these unbounded terms vanishes, and this is possible only for $\delta r$ within an exponentially small range of values---the pinning region.


\section{Calculating late terms in the expansion}
\label{Sec_LateTerms}

In order to determine the forcing in the remainder equation (\ref{RNEqn_unknownForcing}), we require a formula for the $n$th term in the expansion (\ref{Expandu}). In light of the Taylor expansion (\ref{TaylorExpandDiffEq}) of the difference equation (\ref{GenDiffEq_z}) in the continuum limit, we can see that the $n$th term is given by differentiating the $(n-1)$th term four times and integrating twice, and so on. Therefore, if $u_0$ is singular at some point(s) in the complex plane, the expansion (\ref{Expandu}) is divergent in the form of a factorial over a power \cite{adams2003beyond, berry1989uniform, king2001asymptotics}. Hence we propose the ansatz
\begin{equation}
\label{unAnsatz}
u_n \sim (-1)^n\fr{\Gamma( 2n + \beta )}{W(Z)^{2n+\beta}} \lb( f_0(Z) + \fr{1}{2n} f_1(Z) + \fr{1}{(2n)^2} f_2(Z) + \cdots \rb)
\end{equation}
as $n\rightarrow\infty$, in which all dependence on $n$ and $Z$ has been written down explicitly and the large-$n$ limit has been exploited in order to write $u_n$ as an asymptotic series in inverse powers of $n$. Therefore the equation for $u_n$ as $n\rightarrow\infty$, $n \leq N-1$, given by equating terms in (\ref{GenDiffEq_z}) at $\Or(\epsilon^{2n+2})$, is
\begin{equation}
\label{unEq}
0 = 2\sum_{p=1}^{n+1} \fr{\cos^{2p}\psi + \sin^{2p}\psi}{(2p)!} \fr{\od^{2p}u_{n-p+1}}{\od Z^{2p}} - F_M'\lb(u_0\rb)u_n + \cdots,
\end{equation}
where the neglected terms contribute at higher order in $1/n$. 

In light of (\ref{unAnsatz}), the derivative terms in (\ref{unEq}) are $\Or(\Gamma(2n+2+\beta))$, whereas the terms arising due to the expansion of $F_M(u)$ around $u_0$ are merely $\Or(\Gamma(2n+\beta))$. As a result, the leading-order contribution to (\ref{unEq}) is
\begin{equation}
\label{unEq_LeadingOrder}
0 = 2 (-1)^{n+1}\fr{\Gamma(2n+2+\beta)}{W^{2n+\beta+2}} \sum_{p=1}^{n+1}(-1)^p \fr{\cos^{2p}\psi + \sin^{2p}\psi}{(2p)!} (W')^{2p} f_0
\end{equation}
The summation is dominated by $p=\Or(1)$, and so we can replace the upper limit with infinity to give
\begin{equation}
0 = \sum_{p=0}^{\infty}(-1)^p \fr{(W'\cos\psi)^{2p} + (W'\sin\psi)^{2p}}{(2p)!} - 2.
\end{equation}
Evaluating the summation, we therefore have
\begin{equation}
\label{W'Eqn}
0 = \cos(W'\cos\psi) + \cos(W'\sin\psi) - 2.
\end{equation}
This is precisely the eigenvalue equation (\ref{kappa_cond_O1}) derived in Section \ref{Sec_RemEqn}. Hence we set $W'=\kappa$, where $\kappa$ is a (possibly complex) \emph{non-zero} solution of (\ref{kappa_cond_O1}). Recall the existence of real, non-zero solutions is dependent upon the rationality of $\tan\psi$, as discussed in Section \ref{Sec_RemEqn}. We note that both the eigenvalue equation (\ref{kappa_cond_O1}) and the $\Or(\epsilon S_N)$ condition (\ref{kappa_cond_Oeps}) are invariant under $\kappa \rightarrow -\kappa$ and $\kappa \rightarrow \overline{\kappa}$; furthermore, (\ref{kappa_cond_O1}) admits no non-zero, purely imaginary solutions. Therefore we can restrict $W'=\kappa$ to the right half-plane $\Re(\kappa)>0$ without loss of generality. Hence we have
\begin{equation}
\label{WSoln}
W = \kappa \lb( Z - \zeta \rb),
\end{equation}
where $Z=\zeta$ is a singularity of $u_0$, and therefore also of each $u_n$. Since (\ref{unEq}) is linear in $u_n$, the full solution consists of the sum of the contributions from each singularity $\zeta$, which in turn is the sum of the contributions for each eigenvalue $\kappa$, with $\Re(\kappa)>0$.

We now proceed to higher orders in $1/n$ in (\ref{unEq}) in order to determine $\beta$ and $f_0$. As $W$ is a linear function of $Z$, we have
\begin{align}
\label{2pthDerivative}
\fr{\od^{2p}u_{n-p+1}}{\od Z^{2p}} = &\ (-1)^{n+p+1}  \fr{\Gamma(2n+2+\beta)}{W^{2n+\beta+2}} \kappa^{2p} \lb( f_0 + \fr{1}{2n}f_1 + \fr{(2p-2)f_1 + f_2}{(2n)^2} \rb)
\nonumber\\
& {} + (-1)^{n+p} \fr{\Gamma(2n+1+\beta)}{W^{2n+\beta+1}} \kappa^{2p-1} 2p \lb( f_0' + \fr{f_1'}{2n} \rb)
\nonumber\\
& {} + (-1)^{n+p-1} \fr{\Gamma(2n+\beta)}{2W^{2n+\beta}} \kappa^{2p-2} 2p(2p-1) f_0'' + \Or(\Gamma(2n-1+\beta)).
\end{align}
Substitution of (\ref{2pthDerivative}) into (\ref{unEq}) and division by $\Gamma(2n+2+\beta)$ yields
\begin{align}
0 = &\ 2 \lb\{ \fr{(-1)^{n+1}}{W^{2n+\beta+2}} \lb[ \sum_{p=1}^\infty (-1)^p \fr{(\kappa\cos\psi)^{2p} + (\kappa\sin\psi)^{2p}}{(2p)!} \lb( f_0 + \fr{f_1}{2n} + \fr{f_2-2f_1}{(2n)^2} \rb) \rb. \rb.
\nonumber\\
& \lb. \lb. {} - \kappa\lb( \cos\psi \sum_{p=1}^\infty (-1)^{p-1} \fr{(\kappa\cos\psi)^{2p-1}}{(2p-1)!} + \sin\psi \sum_{p=1}^\infty (-1)^{p-1} \fr{(\kappa\sin\psi)^{2p-1}}{(2p-1)!} \rb) \fr{f_1}{(2n)^2} \rb] \rb.
\nonumber\\
& \lb. {} + \fr{(-1)^{n+1}}{W^{2n+\beta+1}} \lb[ \cos\psi \sum_{p=1}^\infty (-1)^{p-1} \fr{(\kappa\cos\psi)^{2p-1}}{(2p-1)!} + \sin\psi \sum_{p=1}^\infty (-1)^{p-1} \fr{(\kappa\sin\psi)^{2p-1}}{(2p-1)!} \rb] \rb.
\nonumber\\
& {} \times \lb. \lb( \fr{1}{2n} - \fr{1+\beta}{(2n)^2} \rb) \lb( f_0' + \fr{f_1'}{2n} \rb) + \fr{(-1)^n}{2W^{2n+\beta}} \lb[ \cos^2\psi \sum_{p=1}^\infty (-1)^{p-1} \fr{(\kappa\cos\psi)^{2p-2}}{(2p-2)!} \rb. \rb.
\nonumber\\
& \lb. \lb. {} + \sin^2\psi \sum_{p=1}^\infty (-1)^{p-1} \fr{(\kappa\sin\psi)^{2p-2}}{(2p-2)!} \rb] \fr{f_0''}{(2n)^2}  \rb\} - F_M'(u_0)  \fr{(-1)^n}{W(Z)^{2n+\beta}} \fr{f_0}{(2n)^2} + \cdots
\end{align}
as $n\rightarrow\infty$. Each of these summations may be evaluated explicitly, yielding
\begin{align}
\label{Expanded_unEq}
0 = &\ 2\lb\{ \fr{(-1)^{n+1}}{W^{2n+\beta+2}} \lb[ \lb( \cos(\kappa\cos\psi) + \cos(\kappa\sin\psi) - 2 \rb) \lb( f_0 + \fr{f_1}{2n} + \fr{f_2-2f_1}{(2n)^2} \rb) \rb. \rb.
\nonumber\\
 & \lb. \lb. {} - \lb( \kappa\cos\psi \sin(\kappa\cos\psi) + \kappa\sin\psi \sin(\kappa\sin\psi) \rb) \fr{f_1}{(2n)^2} \rb] \rb.
\nonumber\\
& \lb. {} + \fr{(-1)^{n+1}}{W^{2n+\beta+1}} \lb( \cos\psi\sin(\kappa\cos\psi) + \sin\psi \sin(\kappa\sin\psi) \rb) \lb( \fr{1}{2n} - \fr{1+\beta}{(2n)^2} \rb) \lb( f_0' + \fr{f_1'}{2n} \rb) \rb.
\nonumber\\
& \lb. {} + \fr{(-1)^n}{2W^{2n+\beta}} \lb( \cos^2\psi\cos(\kappa\cos\psi) + \sin^2\psi \cos(\kappa\sin\psi) \rb) \fr{f_0''}{(2n)^2} \rb\}
\nonumber\\
& {} - F_M'(u_0) (-1)^n \fr{1}{W(Z)^{2n+\beta}} \fr{f_0}{(2n)^2}  + \cdots.
\end{align}
As $\kappa$ satisfies (\ref{kappa_cond_O1}), the first line on the right-hand side of (\ref{Expanded_unEq}) vanishes. Because this includes all $\Or(1)$ terms, we proceed to $\Or(1/n)$ and find that
\begin{equation}
\label{O1/2n_Eqn}
0 = \lb[ \cos\psi\sin(\kappa\cos\psi) + \sin\psi \sin(\kappa\sin\psi) \rb] f_0'.
\end{equation}
If $\kappa$ is real then (\ref{O1/2n_Eqn}) is automatically satisfied (cf. Section \ref{Sec_RemEqn}) and we must continue to $\Or(1/n^2)$. The zero eigenvalue does not contribute to $u_n$, as this would lead to division by zero, so we must have $\tan\psi \in \mathbb{Q}_\infty$. Thus we can define $\tan\psi = m_2/m_1$ as in (\ref{tanphi_realkappa}), in which case $\kappa = 2M\pi(m_1^2+m_2^2)^{1/2}$ as in (\ref{realkappa}), albeit with $M>0$ as we have fixed $\Re(\kappa) > 0$. Consequently, the first three lines of (\ref{Expanded_unEq}) vanish and we are left with
\begin{equation}
0 = f_0'' - F_M'\lb(u_0(Z)\rb) f_0.
\end{equation}
This we have already solved; the complementary functions $g(Z)$ and $G(Z;\zeta)$ are defined in (\ref{CompFunctions}). Hence if $\kappa$ is real then either $f_0=\lambda_{M,\psi} g$ or $f_0=\Lambda_{M,\psi} G$, for some $\psi$-dependent constants $\lambda_{M,\psi}$ and $\Lambda_{M,\psi}$. On the other hand, if $\kappa$ has non-zero imaginary part then (\ref{O1/2n_Eqn}) can be satisfied only if $f_0'=0$ (cf. Section \ref{Sec_RemEqn}), and we therefore set $f_0 = \Omega_{\kappa,\psi}$ in this case, for some $\psi$-dependent constant $\Omega_{\kappa,\psi}$.

It remains to evaluate $\beta$; this is readily achieved by consideration of the singularities of $u_0$. We shall assume that these singularities are all either of the form
\begin{equation}
\label{AlgebraicSingularity}
u_0 = \Or\lb( \lb(Z-\zeta\rb)^{-\gamma} \rb) \ \mathrm{as} \  Z\rightarrow \zeta, \quad -\gamma \notin \mathbb{N}\cup\{0\},
\end{equation}
or
\begin{equation} 
\label{LogarithmicSingularity}
u_0 = \Or\lb( \lb(Z-\zeta\rb)^{-\gamma} h\lb(\ln\lb(Z-\zeta\rb)\rb) \rb) \ \mathrm{as} \ Z\rightarrow \zeta, \quad \gamma\in\mathbb{R},
\end{equation}
for some function $h$ satisfying $h(\ln(t)) \neq At^\alpha$ for any constants $(A,\alpha)\in\mathbb{C}^2$. We shall henceforth refer to the constant $\gamma$ (which we take to be real for simplicity; results are similar for complex $\gamma$) as the \emph{strength} of the singularity at $\zeta$. The systems giving rise to figures \ref{Pic_Const3ExampleBifDiag} and \ref{Pic_35ExampleBifDiag} are both examples of an algebraic singularity (\ref{AlgebraicSingularity}); the system studied in \cite{yulin2010discrete} has a logarithmic singularity (\ref{LogarithmicSingularity}). By inspection of (\ref{unEq}), we can see that if $u_0$ has a singularity of strength $\gamma$ then $u_n$ must have one of strength $2n+\gamma$, as $u_n$ is found by differentiating $u_{n-1}$ four times and integrating twice. Considering the three possible solutions for $f_0$, $g$ has a singularity of strength $\gamma+1$ and $G$ has one of strength $-\gamma-2$, whereas the constant function has none at all. Thus, substituting (\ref{WSoln}) for $W$ into the factorial-over-power ansatz (\ref{unAnsatz}) and expanding $u_n$ near the singularity $\zeta$ for each $f_0$ in turn provides the following:
\begin{equation}
\label{BetaCalculation}
\begin{array}{lclcl}
f_0 = \lambda_{M,\psi} g & \Rightarrow\ \ & 2n + \gamma =  2n + \beta + \gamma + 1 & \Rightarrow\ \ & \beta = -1
\\
f_0 = \Lambda_{M,\psi} G & \Rightarrow\ \ & 2n + \gamma = 2n + \beta - \gamma - 2 & \Rightarrow\ \ & \beta = 2\gamma+2
\\ 
f_0 = \Omega_{\kappa,\psi} & \Rightarrow\ \ & 2n + \gamma = 2n + \beta & \Rightarrow\ \ & \beta = \gamma
\end{array}
\end{equation}
Therefore the contribution to $u_n$ from each singularity $\zeta$ is
\begin{align}
\label{unSoln}
u_n(Z) \sim & \sum_{M=1}^\infty \lb[ \fr{ (-1)^n \Gamma(2n-1) \lambda_{M,\psi} g(Z) }{ \lb[ 2M\pi\sqrt{m_1^2+m_2^2} (Z-\zeta) \rb]^{2n-1}} + \fr{ (-1)^n \Gamma(2n+2\gamma+2) \Lambda_{M,\psi} G(Z;\zeta) }{ \lb[ 2M\pi\sqrt{m_1^2+m_2^2} (Z-\zeta) \rb]^{2n+2\gamma+2} } \rb] 
\nonumber\\
& {} + \sum_{\substack{ \kappa\notin \mathbb{R} \\ \Re(\kappa)>0 }} \fr{(-1)^n \Gamma(2n+\gamma) \Omega_{\kappa,\psi}}{\lb[ \kappa (Z-\zeta) \rb] ^{2n+\gamma}}.
\end{align}
We note that if $\tan\psi \notin \mathbb{Q}_\infty$ then $\lambda_{M.\psi} = \Lambda_{M,\psi} = 0$ in (\ref{unSoln}) for all $M$, there being no real, non-zero eigenvalues. If, on the other hand, $\tan\psi \in \mathbb{Q}_\infty$ then we define $\psi$ as in (\ref{tanphi_realkappa}) and both summations contribute, unless $\psi = \fr{k\pi}{4}$, $k \in \{0,1,\ldots,7\}$, in which case there are no complex eigenvalues so the second summation vanishes. 

Clearly, it is the eigenvalues of smallest modulus which are dominant as $n \rightarrow \infty$. For eigenvalues of equal size, dominance is then determined by comparing the offsets within the $\Gamma$-functions. When $\kappa \in \mathbb{R}$, the dominant eigenvalue is given by $M = 1$. Therefore, if the modulus of the smallest complex eigenvalue is less than $2\pi(m_1^2+m_2^2)^{1/2}$ then the third term is dominant over the other two. Otherwise, the third term is subdominant to the first two, in which case the second dominates the first provided $\gamma > -\fr{3}{2}$. The question of dominance plays no role in determining the remainder, as each contribution can be considered separately by making use of the superposition principle of linear equations. However, it does become significant when calculating the constant $\Lambda_{1,\psi}$, a prerequisite for accurate comparison with numerical results. This will be discussed in detail in the context of a cubic nonlinearity with constant forcing in Section \ref{Sec_Const3FindLambda}.


\section{Optimal truncation and Stokes lines}
\label{Sec_Truncation}

We can now turn our attention to the unknown forcing in the remainder equation (\ref{RNEqn_unknownForcing}). Before we can evaluate it, we must first determine the point of truncation $n = N-1$, desiring to truncate the expansion optimally so the resultant remainder is exponentially small. To this end, we shall treat the contribution to the expansion from each singularity $\zeta$ and each eigenvalue $\kappa$ separately, as each contribution has a different least term. This we are free to do, since both the large-$n$ equation (\ref{unEq}) and the remainder equation (\ref{RNEqn_unknownForcing}) are linear. We shall therefore for the moment work in terms of a general solution pair $(f_0,\beta)$, rather than one of the three specific solutions derived in the previous section. In light of the large-$n$ solution (\ref{unSoln}), each contribution to $u_n$ is minimal with respect to $n$ when
\begin{equation}
\fr{\od}{\od n} \lb| \fr{\epsilon^{2n} \Gamma(2n+\beta) }{ \lb[ \kappa (Z-\zeta) \rb]^{2n+\beta} } \rb| = 0,
\end{equation} 
where $\beta$ is determined by the choice of $f_0$ under consideration. The limit $n\rightarrow\infty$ can be exploited in order to approximate this using Stirling's formula, yielding
\begin{equation}
\fr{\od}{\od n} \lb( \fr{\epsilon^{2n} \sqrt{2\pi} (2n+\beta)^{2n+\beta-1/2} e^{-2n-\beta} }{ \lb|\kappa (Z-\zeta)\rb| ^{2n+\beta} } \rb) = 0,
\end{equation}
which solves to give
\begin{equation}
N \sim \fr{|\kappa(Z-\zeta)|}{2\epsilon} + \nu,
\end{equation}
where $\nu=\Or(1)$ is added to ensure $N$ is an integer.

We are now able to evaluate the forcing due to truncation of (\ref{Expandu}) appearing in (\ref{RNEqn_unknownForcing}). This consists of those terms not accounted for by equating coefficients at $\Or( \epsilon^{2n+2} )$ in (\ref{TaylorExpandDiffEq}) for $n = 0, 1, \ldots, N-1$. Considering the $u_n$ equation (\ref{unEq}), it follows that this forcing, denoted henceforth by $\mathrm{RHS}$, is given by the double summation
\begin{equation}
\label{TruncationForcing}
\mathrm{RHS} \sim - 2\sum_{m=N}^\infty \epsilon^{2m+2} \sum_{p=m-N+2}^{m+1} \fr{\cos^{2p}\psi + \sin^{2p}\psi}{(2p)!} \fr{\od^{2p}u_{m-p+1}}{\od Z^{2p}} + \cdots,
\end{equation}
where the lower limit of summation in $p$ arises because the asymptotic expansion has been truncated after the $N$th term and neglected terms contribute at higher order in $\epsilon$. Since $m\gg1$ and the range $p=\Or(1)$ is dominant, we can make use of (\ref{2pthDerivative}) and Stirling's formula to give
\begin{align}
\mathrm{RHS} \sim & - 2\sqrt{2\pi} \sum_{m=N}^\infty \sum_{p=m-N+2}^{\infty} \lb[ \epsilon^{2m+2} (-1)^{m+1} \fr{( 2m+2+\beta )^{2m+3/2+\beta} e^{-(2m+2+\beta)}}{[\kappa(Z-\zeta)]^{2m+2+\beta}} \rb]
\nonumber\\
& {} \times \lb[ (-1)^p \fr{(\kappa\cos\psi)^{2p} + (\kappa\sin\psi)^{2p}}{(2p)!} \rb] f_0 + \cdots.
\end{align}
After writing $m=N+\hat{m}$, we find that this is dominated by the range $\hat{m}=\Or(1)$ and, because $(t+a)^{t+a-1/2}e^{-t-a} = \exp[ (t+a-1/2)\ln(t+a)-t-a] \sim t^{t+a-1/2}e^{-t}$ as $t\rightarrow\infty$, can be written
\begin{align}
\mathrm{RHS} \sim & - 2\sqrt{2\pi} (-1)^{N+1} \fr{ \epsilon^{2N+2} (2N)^{2N+3/2+\beta} e^{-2N} }{[\kappa(Z-\zeta)]^{2N+2+\beta}} \sum_{\hat{m}=0}^\infty \sum_{p=\hat{m}+2}^{\infty} \lb[ (-1)^{\hat{m}} \fr{\epsilon^{2\hat{m}} (2N)^{2\hat{m}} }{[\kappa(Z-\zeta)]^{2\hat{m}}} \rb]
\nonumber\\
& {} \times \lb[ (-1)^p \fr{(\kappa\cos\psi)^{2p} + (\kappa\sin\psi)^{2p}}{(2p)!} \rb] f_0 + \cdots
\end{align}
as $N \rightarrow \infty$. Reversing the order of summation, we now have
\begin{align}
\mathrm{RHS} \sim & - 2\sqrt{2\pi} (-1)^{N+1} \fr{ \epsilon^{2N+2} (2N)^{2N+3/2+\beta} e^{-2N} }{[\kappa(Z-\zeta)]^{2N+2+\beta}}
\nonumber\\ 
& {} \times \sum_{p=2}^\infty \lb[ (-1)^p \fr{(\kappa\cos\psi)^{2p} + (\kappa\sin\psi)^{2p}}{(2p)!} \sum_{\hat{m}=0}^{p-2} (-1)^{\hat{m}} \lb( \fr{2\epsilon N}{\kappa(Z-\zeta)} \rb)^{2\hat{m}} \rb] f_0 + \cdots.
\end{align}
This we can evaluate, since
\begin{align}
\sum_{p=2}^\infty \lb[ (-1)^p \fr{v^{2p}}{(2p)!} \sum_{m=0}^{p-2} \lb( -w^2 \rb)^m \rb] & = \sum_{p=2}^\infty  (-1)^p \fr{v^{2p}}{(2p)!} \fr{1 + (-1)^p w^{2p-2}}{1+w^2}
\nonumber\\
& = \fr{1}{1+w^2} \lb[ \cos v - 1 + \fr{1}{w^2} \lb( \cosh(vw) - 1 \rb) \rb].
\end{align}
Therefore, since $\kappa$ satisfies (\ref{kappa_cond_O1}), the leading-order forcing due to truncation can be written as
\begin{align}
\label{EvaluatedRHS}
\mathrm{RHS} \sim &\ 2\sqrt{2\pi} (-1)^N \fr{ \epsilon^{2N} (2N)^{2N-1/2+\beta} e^{-2N} }{[\kappa(Z-\zeta)]^{2N-2+\beta} ( \kappa^2(Z-\zeta)^2 + 4\epsilon^2N^2 )} \lb[ \cosh\lb( \fr{2\epsilon N\cos\psi}{Z-\zeta} \rb) \rb.
\nonumber\\
& \lb. {} + \cosh\lb( \fr{2\epsilon N \sin\psi}{Z-\zeta} \rb) - 2 \rb] f_0 + \cdots.
\end{align}

In order to simplify subsequent calculations we now define
\begin{equation}
\kappa(Z-\zeta) = \rho e^{i\theta},
\end{equation}
which gives $N \sim \rho/(2\epsilon) + \nu$. Therefore we can write the prefactor of (\ref{EvaluatedRHS}) as
\begin{align}
& \fr{ \epsilon^{2N} (2N)^{2N-1/2+\beta} e^{-2N} }{[\kappa(Z-\zeta)]^{2N-2+\beta} ( \kappa^2(Z-\zeta)^2 + 4\epsilon^2N^2 )}
\nonumber\\
& \hspace{0.3\textwidth} \sim \fr{\epsilon^{1/2-\beta}}{\sqrt{\rho}} \fr{\lb( 2\epsilon N/\rho \rb)^{2N-1/2+\beta} e^{-i\theta(2N-2+\beta)} e^{-2N}}{e^{2i\theta} + 4\epsilon^2N^2/\rho^2} 
\nonumber\\
& \hspace{0.3\textwidth} = \fr{\epsilon^{1/2-\beta}}{\sqrt{\rho}} \fr{\lb( 1 + 2\epsilon\nu/\rho \rb)^{2N-1/2+\beta} e^{-i\theta(2N-2+\beta)} e^{-2N}}{e^{2i\theta} + 1 + 4\epsilon\nu/\rho + 4\epsilon^2\nu^2/\rho^2}
\nonumber\\
& \hspace{0.3\textwidth} \sim \fr{\epsilon^{1/2-\beta}}{\sqrt{\rho}} \fr{e^{2\epsilon \nu ( \rho/\epsilon + 2\nu - 1/2 + \beta)/\rho} e^{-i\theta(2N-2+\beta)} e^{-\rho/\epsilon - 2\nu}}{e^{2i\theta} + 1}
\nonumber\\
& \hspace{0.3\textwidth} \sim \fr{\epsilon^{1/2-\beta}}{\sqrt{\rho}} \fr{ e^{-i\theta(2N-2+\beta)} e^{-\rho/\epsilon}}{e^{2i\theta} + 1},
\end{align}
and obtain
\begin{align}
\label{RHS_expsmall}
\mathrm{RHS} \sim &\ 2\sqrt{2\pi} (-1)^N \fr{\epsilon^{1/2-\beta}}{\sqrt{\rho}} \fr{ e^{-i\theta(2N-2+\beta)} e^{-\rho/\epsilon}}{e^{2i\theta} + 1} 
\nonumber\\
& \times \lb[ \cosh\lb( \kappa \cos\psi e^{-i\theta} \rb)  + \cosh\lb( \kappa \sin\psi e^{-i\theta} \rb) - 2 \rb] f_0 + \cdots.
\end{align}
Thus we see from the factor $e^{-\rho/\epsilon}$ that RHS is exponentially small.

We now substitute the RHS (\ref{RHS_expsmall}) into the remainder equation (\ref{RNEqn_unknownForcing}) and seek a particular integral. Following Section \ref{Sec_RemEqn}, we write $R_N = e^{\mp i\kappa z} S_N(Z)$ in (\ref{RNEqn_unknownForcing}), where $\kappa$ is as usual a solution of (\ref{kappa_cond_O1}) with $\Re(\kappa)>0$, and Taylor expand the differences in $Z$. This gives
\begin{align}
\label{RemEqnRHS}
& 2\lb[ \cos(\kappa\cos\psi) + \cos(\kappa\sin\psi) - 2 \rb] S_N + 2i\epsilon \lb[ \cos\psi\sin(\kappa\cos\psi) + \sin\psi \sin(\kappa\sin\psi) \rb] S_N'
\nonumber\\
& {} + \epsilon^2 \lb[ \cos^2\psi\cos(\kappa\cos\psi) + \sin^2\psi \cos(\kappa\sin\psi) \rb] S_N'' - \epsilon^2 F_M'(u_0(Z))S_N + \cdots
\nonumber\\
& \hspace{0.7\textwidth} = e^{\pm i\kappa z} \mathrm{RHS} + \cdots.
\end{align}
Note that on the left-hand side, the $\Or(S_N)$ contribution vanishes because $\kappa$ satisfies (\ref{kappa_cond_O1}). The $\Or(\epsilon S_N)$ terms on the left-hand side also vanish if $\kappa$ is real, in which case the leading-order contribution is $\Or(\epsilon^2 S_N)$; otherwise it is $\Or(\epsilon S_N)$.

Now,
\begin{equation}
\exp(\pm i\kappa z - \rho/\epsilon) = \exp\lb[ \pm i \kappa z_0 + \lb( \pm i \kappa \zeta  \pm i \rho e^{i\theta} - \rho \rb) /\epsilon \rb].
\end{equation}
Therefore we can see that, although it remains exponentially small on the real line, $e^{\pm i\kappa z}\mathrm{RHS}$ is maximal with respect to $\theta$ at $\theta = \mp\fr{\pi}{2}$. These values of $\theta$ define the Stokes lines, two emanating from each singularity, at which the main change in $S_N$ will occur. As we are concerned with the behaviour of the solution for real $z$, the Stokes lines of importance are those which cross the real line. Focusing on those singularities in the upper half-plane, so that $\Im(\zeta)>0$, the Stokes line of interest is $\theta=-\fr{\pi}{2}$. We therefore concentrate on solutions $R_N=e^{-i\kappa z}S_N$. Symmetry considerations then allow the contribution from the conjugate singularity at $\overline{\zeta}$ to be recovered simply by taking the complex conjugate $\overline{R}_N = e^{i\overline{\kappa}z}\overline{S}_N$.

In order to capture the effects of maximal forcing, we rescale $\theta$ in the vicinity of the Stokes line as $\theta = -\fr{\pi}{2} + \eta(\epsilon)\hat{\theta}$, where the scaling $\eta(\epsilon)$ is to be determined. The region $\hat{\theta} = \Or(1)$ thus defines the Stokes layer, in which the remainder changes rapidly as coefficients of complementary functions to (\ref{RNEqn_unknownForcing}) vary from zero to non-zero. This gives
\begin{align}
e^{+ i\kappa z} (-1)^N e^{-i\theta(2N-2+\beta)} e^{-\rho/\epsilon} \hspace{-0.27\textwidth} &
\nonumber\\
 \sim &\ \exp\bigg[ i\kappa z_0 + \fr{i}{\epsilon} \lb(\rho e^{-i\pi/2 + i\eta\hat{\theta}} + \kappa\zeta \rb) + iN\pi - i \lb( -\fr{\pi}{2} + \eta\hat{\theta} \rb) \lb( 2N - 2 + \beta \rb) - \fr{\rho}{\epsilon} \bigg]
\nonumber\\
\sim &\ \exp \bigg[ i\kappa \lb( z_0 + \fr{\zeta}{\epsilon} \rb) + \fr{1}{\epsilon} \lb( \rho + i\rho\eta\hat{\theta} - \fr{1}{2}\rho\eta^2\hat{\theta}^2 \rb)
\nonumber\\
& \hspace{0.2\textwidth} {} - i\eta\hat{\theta}\lb( \fr{\rho}{\epsilon} + 2\nu \rb) +\ i\lb( \beta - 2 \rb)\lb(\fr{\pi}{2} - \eta\hat{\theta} \rb) - \fr{\rho}{\epsilon} \bigg]
\nonumber\\
\sim & - e^{i\beta\pi/2} e^{i\kappa( z_0 + \zeta/\epsilon )} e^{-\rho\eta^2\hat{\theta}^2/(2\epsilon)},
\end{align}
which suggests the scaling $\eta = \sqrt{\epsilon}$. We therefore also have
\begin{align}
\fr{\cosh\lb( \kappa\cos\psi e^{-i\theta} \rb) + \cosh\lb( \kappa \sin\psi e^{-i\theta} \rb) - 2 }{e^{2i\theta} + 1} \hspace{-0.2\textwidth} &  
\nonumber\\
\sim &\ \kappa \lb[ \cos\psi \sin(\kappa\cos\psi) + \sin\psi \sin(\kappa\sin\psi) \rb] \lb( -\fr{1}{2} + \fr{3i}{4} \sqrt{\epsilon}\hat{\theta} \rb) 
\nonumber\\
& {} + \fr{i}{4}\kappa^2 \lb[ \cos^2\psi\cos(\kappa\cos\psi) + \sin^2\psi\cos(\kappa\sin\psi) \rb]\sqrt{\epsilon}\hat{\theta} + \cdots,
\end{align}
where we have made use of the fact that $\kappa$ satisfies (\ref{kappa_cond_O1}) in order to eliminate terms; note that if $\kappa$ is real the first contribution to the right-hand side also vanishes and the second is simply equal to $i\kappa^2\sqrt{\epsilon}\hat{\theta}/4$. We shall now consider the two types of eigenvalue in turn, $\kappa \in \mathbb{R}$ and $\kappa \notin \mathbb{R}$, in order to elucidate precisely what contribution to the remainder each makes. 


\subsection{Contribution to $R_N$ from $\kappa\in\mathbb{R}$}

As discussed in Section \ref{Sec_RemEqn}, $\kappa$ can be real and non-zero only if $\tan\psi \in \mathbb{Q}_\infty$, in which case we define $\tan\psi = m_2/m_1$ as in (\ref{tanphi_realkappa}). This gives $\kappa = 2M\pi(m_1^2+m_2^2)^{1/2}$ as in (\ref{realkappa}), with $M>0$ due to our restriction that $\Re(\kappa) > 0$. Thus if $\kappa$ is real then the leading-order balance in (\ref{RemEqnRHS}) is
\begin{equation}
\epsilon^2 S_N'' - \epsilon^2 F_M'(u_0(Z))S_N =  - i\sqrt{\pi/2} e^{i\beta\pi/2} \fr{\epsilon^{1-\beta}}{\sqrt{\rho}} \kappa^2 e^{i \kappa ( z_0 + \zeta/\epsilon )} \hat{\theta} e^{-\rho\hat{\theta}^2/2} f_0(Z) + \cdots,
\end{equation}
where either $f_0 = \lambda_{M,\psi} g(Z)$ and $\beta = -1$ or $f_0 = \Lambda_{M,\psi} G(Z;\zeta)$ and $\beta = 2\gamma+2$. Writing
\begin{equation}
S_N(Z) = \epsilon^{-\beta} e^{i \kappa ( z_0 + \zeta/\epsilon )} f_0(Z) \hat{S}_N(\hat{\theta}),
\end{equation}
we have
\begin{equation}
\epsilon^2 S_N''(Z) = \epsilon^{1-\beta} e^{i \kappa ( z_0 + \zeta/\epsilon )} f_0(Z) \fr{\kappa^2}{\rho^2} \fr{\od^2 \hat{S}_N}{\od \hat{\theta}^2} + \Or\lb( \epsilon^{3/2-\beta} e^{i\kappa( z_0 + \zeta/\epsilon )} \rb).
\end{equation}
Thus
\begin{equation}
\fr{\od^2 \hat{S}_N}{\od \hat{\theta}^2} \sim - i\sqrt{\pi/2} e^{i\beta\pi/2} \rho^{3/2} \hat{\theta} e^{-\rho\hat{\theta}^2/2} + \cdots.
\end{equation}
Imposing the boundary condition $\hat{S}_N\rightarrow0$ as $\hat{\theta}\rightarrow-\infty$, i.e. that the particular integral due to truncation and that due to the deviation from the Maxwell point (\ref{PI_deltar}) are bounded in the same far-field, this has leading-order solution
\begin{equation}
\hat{S}_N(\hat{\theta}) \sim \fr{i\pi}{2} e^{i\beta\pi/2} \erfc\lb(-\hat{\theta}\sqrt{\fr{\rho}{2}}\rb),
\end{equation}
where $\erfc(t)$ is the complementary error function
\begin{equation}
\erfc(t) := \fr{2}{\sqrt{\pi}} \int_t^\infty e^{-s^2} \od s.
\end{equation}
Therefore the exponentially small terms
\begin{equation}
\label{SwitchedOn_realkappa}
R_N \sim \sum_{M=1}^\infty i\pi e^{i\beta\pi/2} \epsilon^{-\beta} e^{-2M\pi i \sqrt{m_1^2+m_2^2}( z - z_0 - \zeta/\epsilon )}  f_0(Z),
\end{equation}
for each singularity $\zeta$ in the upper half-plane, are present to the right of the Stokes layer. By symmetry, the contribution from the conjugate singularity at $Z=\overline{\zeta}$ is simply the complex conjugate of (\ref{SwitchedOn_realkappa}). Note that here $e^{-i\kappa z} = \exp[-2M \pi i (m_1^2+m_2^2)^{1/2} z] = 1$ on lattice points, as $\tan\psi \in \mathbb{Q}_\infty$ for real $\kappa$. 


\subsection{Contribution to $R_N$ from $\kappa\notin\mathbb{R}$}

We now consider the forcing of the remainder equation (\ref{RemEqnRHS}) due to those eigenvalues with $\Im(\kappa)\neq 0$ (recall that we have set $\Re(\kappa)>0$). We know that for such $\kappa$
\begin{equation}
\cos\psi\sin(\kappa\cos\psi) + \sin\psi \sin(\kappa\sin\psi) \neq 0,
\end{equation}
(cf. Section \ref{Sec_RemEqn}) and so the leading-order balance in (\ref{RemEqnRHS}) is
\begin{equation}
2\epsilon i S_N'(Z) =  \sqrt{2\pi} e^{i\beta\pi/2} \fr{\epsilon^{1/2-\beta}}{\sqrt{\rho}} e^{i\kappa( z_0 + \zeta/\epsilon )} e^{-\rho\hat{\theta}^2/2} f_0 + \cdots,
\end{equation}
where $f_0 = \Omega_{\kappa,\psi}$, a constant, and $\beta=\gamma$. Writing
\begin{equation}
S_N(Z) = \epsilon^{-\beta} e^{i\kappa( z_0 + \zeta/\epsilon )} f_0 \hat{S}_N(\hat{\theta}),
\end{equation}
we have
\begin{equation}
\epsilon S_N'(Z) = \epsilon^{1/2-\beta} e^{i\kappa( z_0 + \zeta/\epsilon )} f_0 \fr{\kappa}{\rho} \fr{\od \hat{S}_N}{\od \hat{\theta}} + \Or\lb( \epsilon^{1-\beta} e^{i\kappa( z_0 + \zeta/\epsilon )} \rb).
\end{equation}
Thus
\begin{equation}
\label{SNEqn_complexkappa}
\fr{\od \hat{S}_N}{\od \hat{\theta}} \sim -i\sqrt{\pi/2} e^{i\beta\pi/2} \sqrt{\rho} e^{-\rho\hat{\theta}^2/2} + \cdots.
\end{equation}
Since $R_N = e^{-i\kappa z}S_N$ and $\Im(\kappa)\neq0$, when $\Im(\pm\kappa)>0$ we have $R_N \rightarrow 0$ as $z \rightarrow \mp\infty$. Although these contributions to the remainder are bounded in the pertinent limit, in the opposite limit we have $e^{-i\kappa z} \rightarrow \infty$ as $z \rightarrow \pm\infty$ when $\Im(\pm\kappa)>0$. However, the resultant unbounded growth may be prevented by choosing appropriately the constant of integration when integrating (\ref{SNEqn_complexkappa}). Doing this, we have
\begin{equation}
\hat{S}_N(\hat{\theta}) \sim \lb\{ \begin{array}{lcl}
\fr{i\pi}{2} e^{i\beta\pi/2} \erfc\lb( \hat{\theta}\sqrt{\fr{\rho}{2}} \rb), & & \Im(\kappa)>0,
\\
-\fr{i\pi}{2} e^{i\beta\pi/2} \erfc\lb( -\hat{\theta}\sqrt{\fr{\rho}{2}} \rb), & & \Im(\kappa)<0.
\end{array} \rb. 
\end{equation}
Therefore the exponentially small terms
\begin{equation}
\label{SwitchedOn_Imkappapositive}
R_N \sim \sum_{\substack{ \Re(\kappa)>0, \\ \Im(\kappa)>0 }} i\pi e^{i\beta\pi/2} \epsilon^{-\beta} e^{i\kappa( z_0 + \zeta/\epsilon )} e^{-i\kappa z} f_0
\end{equation}
are present to the left of the Stokes layer, and the exponentially small terms
\begin{equation}
\label{SwitchedOn_Imkappanegative}
R_N \sim \sum_{\substack{ \Re(\kappa)>0, \\ \Im(\kappa)<0 }} -i\pi e^{i\beta\pi/2} \epsilon^{-\beta} e^{i\kappa( z_0 + \zeta/\epsilon )} e^{-i\kappa z} f_0
\end{equation}
to the right, for each singularity $\zeta$ in the upper half-plane. By symmetry, the contributions from the conjugate singularity at $Z=\overline{\zeta}$ are the complex conjugates of (\ref{SwitchedOn_Imkappapositive}) and (\ref{SwitchedOn_Imkappanegative}). Note that, due to our selection of the constants of integration, the Stokes lines relevant to complex $\kappa$ do not switch on any exponentially growing terms as they are crossed; in fact, the terms which are switched on decay exponentially in the fast scale $z$. Thus contributions from $\kappa \notin \mathbb{R}$ remain exponentially small in the far-fields and play no role in selecting the leading-order solution.


\section{The width of the pinning region}
\label{Sec_SnakeWidth}

Although we have verified that complex $\kappa$ do not produce any unbounded terms in the remainder, we have yet to deal with the exponentially growing contributions from real $\kappa$. Because $G(Z;\zeta)$ has coefficient zero to the left of the Stokes lines and the particular integral $P(Z)$ is bounded as $Z \rightarrow -\infty$ (cf. (\ref{FarFieldP})), the remainder is bounded to the left of the Stokes layer. On the other hand, $G$ has non-zero coefficient to the right of the Stokes layer, and both $G$ and $P$ experience exponential growth as $Z \rightarrow \infty$ (cf. (\ref{FarFieldG}) and (\ref{FarFieldP})). We must eliminate these unbounded terms if the asymptotic expansion is to remain uniform. Note that we shall now evaluate our solution on the lattice points, and so have $e^{-i\kappa z} \equiv 1$ on lattice points in (\ref{SwitchedOn_realkappa}), as $\tan\psi \in \mathbb{Q}_\infty$ for real $\kappa$.

$G$ and $P$ are given in the far-field by (\ref{FarFieldG}) and (\ref{FarFieldP}), respectively. In light of (\ref{SwitchedOn_realkappa}), the dominant terms which are switched on are given by those singularities closest to, and equidistant from, the real line, with $M=1$. For the sake of simplicity, we shall assume henceforth that there are only two such singularities. In this instance, focusing on the exponentially growing complementary function $G$, the leading-order contribution which is switched on as the Stokes lines are crossed is
\begin{equation}
\sim -i\pi e^{i\gamma\pi} \epsilon^{-2\gamma-2} e^{2\pi \sqrt{m_1^2+m_2^2} i( z_0 + \zeta/\epsilon )} \Lambda_{1,\psi} G(Z;\zeta) + \cc,
\end{equation}
where we have written $\kappa = 2M\pi(m_1^2+m_2^2)^{1/2}$ with $M=1$. Note that, as we are focusing solely on $f_0 = \Lambda_{1,\psi}G$, we have $\beta = 2\gamma+2$. Including the particular integral $P(Z)$ due to the forcing $\epsilon^2\delta r F_{r,M}(u_0)$ in (\ref{RNEqn_unknownForcing}), the remainder in the far-field $Z\rightarrow\infty$ is therefore
\begin{align}
\label{FarFieldRN}
R_N \sim &\ \lb\{ \fr{\pi |\Lambda_{1,\psi}| e^{-2\pi \sqrt{m_1^2+m_2^2} \Im(\zeta)/\epsilon}}{\epsilon^{2\gamma+2} \alpha_+^2D_+} \cos\lb[ 2\pi z_0 \sqrt{m_1^2+m_2^2} + \chi \rb] \rb.
\nonumber\\
\\ & \lb. {} + \fr{\delta r \int_{u_-}^{u_+} F_{r,M}(v) \od v }{ 2 \alpha_+^2D_+} \rb\} e^{\alpha_+Z},
\end{align}
where
\begin{equation}
\label{chi_eq}
\chi = - \fr{\pi}{2} + \gamma \pi + \fr{2\pi}{\epsilon} \Re(\zeta) \sqrt{m_1^2+m_2^2} + \Arg\lb(\Lambda_{1,\psi}\rb),
\end{equation}
and we have made use of the far-field representations (\ref{FarFieldG}) and (\ref{FarFieldP}) of $G$ and $P$. For the expansion to remain uniform as $Z\rightarrow\infty$, we require the coefficient of these unbounded terms to be zero. This is true if
\begin{equation}
\label{PinningEqn}
\delta r = -\fr{2\pi |\Lambda_{1,\psi}| e^{-2\pi \sqrt{m_1^2+m_2^2} \Im(\zeta)/\epsilon}}{\epsilon^{2\gamma+2} \int_{u_-}^{u_+} F_{r,M}(v) \od v} \cos \lb[ 2\pi z_0 \sqrt{m_1^2+m_2^2} + \chi \rb],
\end{equation}
thus fixing the origin of the front $z_0$ to be one of two values modulo the effective lattice spacing $(m_1^2+m_2^2)^{-1/2}$. Furthermore, real solutions exist only if
\begin{equation}
\label{DiscreteSnakingWidth}
|\delta r| \leq \fr{2\pi|\Lambda_{1,\psi}| e^{-2\pi \sqrt{m_1^2+m_2^2} \Im(\zeta)/\epsilon}}{\epsilon^{2\gamma+2} |\int_{u_-}^{u_+} F_{r,M}(v) \od v|}.
\end{equation}
i.e. stationary fronts exist only for $\delta r$ within this (exponentially small) region. (\ref{DiscreteSnakingWidth}) is the width of the pinning region in which one-dimensional front solutions to (\ref{GenDiffEq}) pin to the underlying lattice. Furthermore, as localised solutions are constructed from back-to-back stationary fronts, (\ref{DiscreteSnakingWidth}) provides a formula for the width of the pinning region; we believe this to be the first time such a result has been reported in full. (\ref{DiscreteSnakingWidth}) should be compared with (17) in \cite{kozyreff2013analytical}, in which an analogous calculation for fronts oriented with respect to a hexagonal pattern is carried out. The quantities $|\Delta k|$ and $\pi/\lambda$ in \cite{kozyreff2013analytical} correspond respectively to $2\pi\sqrt{m_1^2+m_2^2}$ and $\Im(\zeta)$ in (\ref{DiscreteSnakingWidth}). Due to the general nature of the analysis in \cite{kozyreff2013analytical}, only the exponential part of the snaking width is derived; because we have studied the specific problem (\ref{GenDiffEq}), we are able to derive a complete formula, including the algebraic scaling. Continuing the calculation to its conclusion in this manner thus confirms the general results of \cite{kozyreff2013analytical}. Note the constant $\Lambda_{1,\psi}$ is at present undetermined; in fact, it cannot be determined analytically due to the linear nature of the large-$n$ equation (\ref{unEq}). However, the leading-order contribution to (\ref{unEq}) as $Z \rightarrow \zeta$ yields a recurrence relation which can in principle be used to obtain a good approximation to $\Lambda_{1,\psi}$ \cite{adams2003beyond, dean2011exponential, king2001asymptotics}. As this must be done on a case-by-case basis for each choice of $F(u;r)$, we defer further discussion of the calculation of $\Lambda_{1,\psi}$ to Section \ref{Sec_Examples}, in which we shall consider the two specific examples presented in figures \ref{Pic_Const3ExampleBifDiag} and \ref{Pic_35ExampleBifDiag}.


\section{The snakes-and-ladders bifurcation equations}
\label{Sec_SnakeBifEqs}

Armed with the full asymptotic expansion for a stationary front, including exponentially growing terms in the remainder, we are now able to construct localised solutions, or spatially homoclinic connections to the constant solution $u_-$ via $u_+$, by means of matching two stationary fronts back-to-back. Such a solution consists of an up-front $u(\epsilon z - \epsilon z_0)$ matched to a distant down-front $u( -\epsilon z + \epsilon z_0 + L/\epsilon)$, where $L>0$ is an $\Or(1)$ constant. Note that the down-front is produced by applying the rotation $(\psi,z_0) \rightarrow (\psi + \pi,-z_0)$ to $u(Z)$, followed by the translation $\epsilon z_0 \rightarrow \epsilon z_0 + L/\epsilon$. Therefore the origin of the up-front is at $z = z_0$, as before, whereas the translation of the down-front to the right shifts its origin to $-z = -z_0 - L/\epsilon^2$. The scaling of the front separation $L/\epsilon$ is motivated by the fact that the exponentially growing contribution to the remainder (\ref{FarFieldRN}) is no longer exponentially small when $Z = \Or( 1/\epsilon )$ and is positive. This allows us to observe the interplay between three exponentially small effects: the locking of fronts to the lattice, the deviation from the Maxwell point and the front matching error. The first two are responsible for the existence of the pinning region, as already shown in Section \ref{Sec_SnakeWidth}; we shall see now that the third is responsible for the way the solution curves are skewed to the right of the pinning region when the localised patch is small enough, e.g. figures \ref{Pic_Const3ExampleBifDiag} and \ref{Pic_35ExampleBifDiag}. 

From the far-field expansions (\ref{FarFieldu0}) and (\ref{FarFieldRN}), we see that an up-front \linebreak \mbox{$u \sim u_0 + \cdots + R_N$} is given by
\begin{align}
\label{upFront}
u \sim &\ u_+ - D_+ e^{-\alpha_+ Z} + \lb\{ \fr{\pi |\Lambda_{1,\psi}| e^{-2\pi \sqrt{m_1^2+m_2^2} \Im(\zeta)/\epsilon}}{\epsilon^{2\gamma+2} \alpha_+^2D_+} \cos\lb[ 2\pi z_0 \sqrt{m_1^2+m_2^2} + \chi \rb] \rb.
\nonumber\\
& \lb. {} + \fr{\delta r \int_{u_-}^{u_+} F_{r,M}(v) \od v }{ 2 \alpha_+^2D_+} \rb\} e^{\alpha_+ Z} 
\end{align}
as $Z \rightarrow \infty$. By symmetry, the down-front is given within the matching region by (\ref{upFront}) under the combined rotation and translation $(\psi, z_0) \rightarrow (\psi + \pi, - z_0 - L/\epsilon^2 )$. Thus $Z \rightarrow -Z + L/\epsilon$ and we have
\begin{align}
\label{DownFront}
u \sim & \ u_+ - D_+ e^{\alpha_+( Z - L/\epsilon )} + \lb\{ \fr{\pi |\Lambda_{1,\psi}| e^{-2\pi \sqrt{m_1^2+m_2^2} \Im(\zeta)/\epsilon}}{\epsilon^{2\gamma+2} \alpha_+^2D_+} \rb.
\nonumber\\
& \lb. {} \times \cos\lb[ 2\pi \lb( -z_0 - \fr{L}{\epsilon^2} \rb) \sqrt{m_1^2+m_2^2} + \chi \rb]+ \fr{\delta r \int_{u_-}^{u_+} F_{r,M}(v) \od v }{ 2 \alpha_+^2D_+} \rb\} e^{-\alpha_+( Z - L/\epsilon )} 
\end{align}
as $( - Z + L/\epsilon ) \rightarrow \infty$. Note that we have not yet eliminated exponentially growing terms; these are necessary in order to match exponentially growing and decaying terms between fronts. Unbounded terms are then removed by adding the up-front and down-front together and subtracting matched parts, following the usual method of matched asymptotic expansions. 

Matching growing and decaying exponential terms in the matching region, we obtain
\begin{align}
\label{MatchGrowExp}
- D_+ e^{-\alpha_+L/\epsilon} = &\ \fr{\pi |\Lambda_{1,\psi}| e^{-2\pi \sqrt{m_1^2+m_2^2} \Im(\zeta)/\epsilon}}{\epsilon^{2\gamma+2} \alpha_+^2D_+} \cos\lb[ 2\pi z_0 \sqrt{m_1^2+m_2^2} + \chi \rb]
\nonumber\\
& {} + \fr{\delta r \int_{u_-}^{u_+} F_{r,M}(v) \od v }{ 2 \alpha_+^2D_+},
\\
\label{MatchDecayExp}
- D_+ e^{-\alpha_+L/\epsilon} = &\ \fr{\pi |\Lambda_{1,\psi}| e^{-2\pi \sqrt{m_1^2+m_2^2} \Im(\zeta)/\epsilon}}{\epsilon^{2\gamma+2} \alpha_+^2D_+} \cos\lb[ -2\pi \lb( z_0 + \fr{L}{\epsilon^2} \rb) \sqrt{m_1^2+m_2^2} + \chi \rb]
\nonumber\\
& + \fr{\delta r \int_{u_-}^{u_+} F_{r,M}(v) \od v }{ 2\alpha_+^2D_+}.
\end{align}
We therefore have
\begin{equation}
\label{EquateCosines}
\cos\lb[ 2\pi z_0 \sqrt{m_1^2+m_2^2} + \chi \rb] = \cos\lb[ - 2\pi \lb( z_0 + \fr{L}{\epsilon^2} \rb) \sqrt{m_1^2+m_2^2} + \chi \rb].
\end{equation}
Solving (\ref{EquateCosines}) provides two cases to consider: firstly
\begin{equation}
\label{SnakesCondition}
z_0 = - \fr{L}{2\epsilon^2} + \fr{k}{2\sqrt{m_1^2+m_2^2}},
\end{equation}
and secondly
\begin{equation}
\label{LaddersCondition}
\fr{L}{\epsilon} = \lb( \fr{\chi}{\pi} + k \rb) \fr{\epsilon}{\sqrt{m_1^2+m_2^2}},
\end{equation}
where $k$ is some integer, chosen so that $L/\epsilon \gg 1$.

\subsection{The snakes}

Suppose first that (\ref{SnakesCondition}) holds. Substituting for $z_0$ in (\ref{MatchGrowExp}) and rearranging, we gain the bifurcation equation
\begin{align}
\label{GeneralSnakeEq}
\delta r = & - \fr{2}{\int_{u_-}^{u_+} F_{r,M}(v) \od v } \bigg\{ \fr{\pi |\Lambda_{1,\psi}|}{\epsilon^{2\gamma+2}} e^{-2\pi\sqrt{m_1^2+m_2^2} \Im(\zeta)/\epsilon} \cos\lb[ \fr{\pi L}{\epsilon^2} \sqrt{m_1^2+m_2^2} k\pi - \chi \rb]
\nonumber\\
& {} + \alpha_+^2D_+^2 e^{-\alpha_+ L/\epsilon} \bigg\},
\end{align}
from which the front separation $L/\epsilon$ may be determined. As (\ref{GeneralSnakeEq}) is 2-periodic in $k$, only the parity of $k$ is of importance when determining $L$; thus (\ref{GeneralSnakeEq}) describes two distinct snaking solution curves with phases that differ by $\pi$. Each solution is unique up to translations in $Z = \epsilon (z-z_0)$ by integer multiples of the effective (slow-scale) lattice spacing $\epsilon(m_1^2+m_2^2)^{-1/2}$. Furthermore, inspection of (\ref{SnakesCondition}) indicates that the localised solution is site-centred if $k$ is even and bond-centred if $k$ is odd. The second term on the right-hand side of (\ref{GeneralSnakeEq}), corresponding to the front matching error, skews the solution curves to the right of the snaking region for small enough $L$. However, as $L$ increases this term rapidly becomes negligible, in which case the snaking curves are confined to the exponentially small parameter range defined in (\ref{DiscreteSnakingWidth})---the pinning region. $L$ is free to increase without bound, resulting in an infinite multiplicity of localised solutions within this range.

\subsection{The ladders}

Now suppose that (\ref{LaddersCondition}) holds. Since $k$ is arbitrary, in this case the front separation $L/\epsilon$ may take one of a discrete set of values, provided the constraints $L>0$ and $k = \Or(1/\epsilon^2)$ (because $L=\Or(1)$) are satisfied. The origin $z_0$ of the up-front may then be found by solving (\ref{MatchGrowExp}), rewritten here as
\begin{align}
\label{GeneralLadderEq}
\delta r = & - \fr{2}{\int_{u_-}^{u_+} F_{r,M}(v) \od v } \bigg\{ \fr{\pi |\Lambda_{1,\psi}|}{\epsilon^{2\gamma+2}} e^{-2\pi\sqrt{m_1^2+m_2^2}\Im(\zeta)/\epsilon} \cos\lb[ 2\pi z_0 \sqrt{m_1^2+m_2^2} + \chi \rb]
\nonumber\\
& {} + \alpha_+^2D_+^2 e^{-\alpha_+ L/\epsilon} \bigg\}.
\end{align}
This equation therefore describe the `ladders' of the snakes-and-ladders bifurcation diagram. Each $k$ corresponds to a single rung of the ladder, which may be parametrised by $z_0$ in the range $[0,(m_1^2+m_2^2)^{-1/2})$. The deviation $\delta r$ from the Maxwell point for each $z_0$ is then provided by (\ref{GeneralLadderEq}), which has solutions in the same range of values of $\delta r$ as (\ref{GeneralSnakeEq}), as expected. Furthermore, each rung in fact represents two solution curves, corresponding to the two solutions of (\ref{GeneralLadderEq}) in the range $z_0 \in [0,(m_1^2+m_2^2)^{-1/2})$. These two solutions coincide at the maximum and minimum of the cosine, representing the bifurcation points at which the rungs meet the snakes. Note that each rung originates on one snake at $z_0 = 0$ and terminates on the other at $z_0 = (m_1^2+m_2^2)^{-1/2}/2$, linking the two snaking solution curves. 


\section{Examples}
\label{Sec_Examples}

We shall now demonstrate the application of the general results (\ref{DiscreteSnakingWidth}), (\ref{GeneralSnakeEq}) and (\ref{GeneralLadderEq}) to two specific choices of $F(u;r)$ in (\ref{GenDiffEq}). Furthermore, we shall describe how the constant $\Lambda_{1,\psi}$ can be calculated on the axes and principal diagonals, and discuss the difficulties presented by other orientations of $z$. Note that we are now interested only in those orientations for which the width of the pinning region is non-zero, and so assume that $\tan\psi \in \mathbb{Q}_\infty$ throughout the present section.


\subsection{Cubic nonlinearity with constant forcing}
\label{Sec_Constant3Example}

Our first example, arguably the simplest form of (\ref{GenDiffEq}) to exhibit snaking, is
\begin{equation}
\label{Const3DiffEq}
\fr{\pd \hat{u}}{\pd t} = \Delta \hat{u} + \hat{r} + \hat{s}\hat{u} - \hat{u}^3.
\end{equation}
An example bifurcation diagram and solutions for the bistable parameter range $\hat{s} > 0$ are shown in figure \ref{Pic_Const3ExampleBifDiag}, with $\psi = 0$; the unhatted variables in those figures correspond to hatted ones here. The two constant, stable solutions (both non-zero) are connected via an unstable branch, thus forming an S-shaped solution curve in parameter space. This results in a region of bistability, within which is the pinning region. We note that a system similar to (\ref{Const3DiffEq}) was studied in \cite{clerc2011continuous}, in which discreteness was incorporated by replacing a constant coefficient in a partial differential equation with a spatially periodic function, rather than through a difference operator as is the case here. However, that work presents an incomplete description of the snaking phenomenon, due to its failure to incorporate exponentially small terms. 

(\ref{Const3DiffEq}) is not in the form (\ref{GenDiffEq}); we remedy this by performing the rescaling $(\hat{u},\hat{r},\hat{s}) = (\epsilon u, \epsilon^3 r, \epsilon^2 s)$, yielding
\begin{equation}
\label{Constant3DiffEq_scaled}
\fr{\pd u}{\pd t} = \Delta u - \epsilon^2 \lb( -r - su + u^3 \rb),
\end{equation}
describing (\ref{Const3DiffEq}) close to the transition between monostability and bistability. Note we have $F(u;r) = -r - su + u^3$. Defining $Z$ as in (\ref{Z}) and setting $u \equiv u(Z)$, the leading-order continuum approximation is
\begin{equation}
\label{Constant3DiffEq_LO}
0 = \fr{\od^2 u_0}{\od Z^2} + r + su_0 - u_0^3.
\end{equation}
Imposing the condition (\ref{rM_integral}), this can be integrated to give the leading-order solution
\begin{equation}
\label{u0_Const3}
u_0 = \sqrt{s}\tanh\lb( \sqrt{\fr{s}{2}} Z \rb)
\end{equation}
and the value of the Maxwell point $r_M = 0$. Hence $u \rightarrow u_\pm = \pm\sqrt{s}$ as $Z \rightarrow \pm\infty$. The sign of the square root in (\ref{u0_Const3}) has been chosen in order that $u_+ > u_-$; the front of opposite orientation can be obtained by exploiting the reversibility of (\ref{Const3DiffEq}). 

From (\ref{u0_Const3}), we see that the singularities $\zeta$ of $u_0$ are
\begin{equation}
\zeta = \zeta_m = \sqrt{\fr{1}{2s}} (2m+1) \pi i
\end{equation}
each of which has strength $\gamma = 1$. Thus (\ref{chi_eq}) yields $\chi = \fr{\pi}{2} + \Arg\lb(\Lambda_{1,\psi}\rb)$. The dominant singularities are those nearest (and equidistant from) the real line, namely $\zeta_0$ and $\zeta_{-1} = \overline{\zeta}_0$. In addition, because
\begin{equation}
u_0 \sim \sqrt{s} \lb( 1 - 2 e^{-\sqrt{2s}Z} \rb)
\end{equation}
as $Z \rightarrow \infty$ and $F_{r,M}(u) \equiv -1$ we have
\begin{equation}
\alpha_+ = \sqrt{2s}, \qquad D_+ = 2\sqrt{s}, \qquad \int_{u_-}^{u_+} F_{r,M}(v) \od v = -2\sqrt{s}.
\end{equation}
Note that, although $\alpha_+$ is defined by (\ref{alpha}), it is simpler in practice to simply read it off from the leading-order exponential in the far-field. 

We are now almost in possession of the requisite detail to write down the bifurcation equations (\ref{GeneralSnakeEq})-(\ref{GeneralLadderEq}) in terms of the parameters of the scaled equation (\ref{Constant3DiffEq_scaled}). The only parameters as yet undetermined are the $\Lambda_{1,\psi}$; we discuss their calculation in detail in Section \ref{Sec_Const3FindLambda}. From (\ref{DiscreteSnakingWidth}), we see that the snaking width for (\ref{Constant3DiffEq_scaled}) is
\begin{equation}
|\delta r| \leq \fr{\pi|\Lambda_{1,\psi}| e^{-\pi^2\sqrt{2(m_1^2+m_2^2)/\epsilon^2s}}}{\epsilon^4 \sqrt{s}}.
\end{equation}
It is instructive to write this in terms of the original, unscaled variables of (\ref{Const3DiffEq}), in which case $\hat{s}$ provides the small variable. Reabsorbing the scalings in $\epsilon$, we obtain the unscaled snaking width as
\begin{equation}
\label{Const3_SnakeWidth_Unscaled}
|\hat{r}| \leq \fr{\pi|\Lambda_{1,\psi}| e^{-\pi^2\sqrt{2(m_1^2+m_2^2)/\hat{s}}}}{\sqrt{\hat{s}}},
\end{equation}
since the Maxwell point in this case is zero. For the sake of brevity, we omit to write out the snaking bifurcation equations (\ref{GeneralSnakeEq}) and (\ref{GeneralLadderEq}) for the present system.


\subsubsection{Finding $\Lambda_{1,\psi}$}
\label{Sec_Const3FindLambda}

All that remains for a comprehensive comparison between numerical computations of (\ref{Const3DiffEq}) and our analytical predictions is the evaluation of the constants $\Lambda_{1,\psi}$. Unfortunately, this cannot be done analytically due to the linear nature of the large-$n$ equation (\ref{unEq}). However, they may in principle be calculated directly through the iteration of a certain recurrence relation arising from the behaviour of the solution near the singularity $\zeta_m$. 

As the singularity in the leading-order front $u_0$ (\ref{u0_Const3}) has strength $\gamma = 1$, we have (cf. the discussion around (\ref{AlgebraicSingularity})-(\ref{BetaCalculation}))
\begin{equation}
\label{un_NearSing}
u_n \sim \fr{U_n}{(Z-\zeta_m)^{2n+1}},
\end{equation}
as $Z\rightarrow\zeta_m$, for some sequence of constants $U_n$. Upon substitution of this ansatz into (\ref{TaylorExpandDiffEq}), taking the leading-order terms in $(Z-\zeta_m)^{-1}$ leads to
\begin{equation}
\label{RecurrenceRelation}
0 = 2 \sum_{p=1}^{n+1} \fr{\lb( \cos^{2p}\psi + \sin^{2p}\psi \rb) \Gamma\lb(2n+3\rb)}{(2p)!\Gamma\lb(2n-2p+3\rb)}  U_{n-p+1} - \sum_{p_1=0}^n \sum_{p_2=0}^{n-p_1} U_{p_1} U_{p_2} U_{n-p_1-p_2}.
\end{equation}
Iteration of this recurrence relation therefore yields the sequence $U_n$. In principle, we may then compare (\ref{un_NearSing}) with the analytical formula (\ref{unSoln}) for $u_n$ as $n \rightarrow \infty$ in order to find $\Lambda_{1,\psi}$. 

There are three types of contribution to (\ref{unSoln}), two arising from real and one from complex eigenvalues (where each eigenvalue $\kappa$ is a solution of (\ref{kappa_cond_O1}); recall that we have set $\Re(\kappa) > 0$ without loss of generality in Section \ref{Sec_LateTerms}). The real eigenvalues are characterised by the integers $M$, and it is clear that the dominant one is given by $M=1$. Furthermore, in the present example the term involving $G(Z;\zeta_m)$ dominates the one involving $g(Z)$ (cf. the discussion after (\ref{unSoln})). As we are not in general able to determine complex eigenvalues analytically, we shall for now merely denote by $K$ the eigenvalue $\kappa \notin \mathbb{R}$ of smallest modulus in the quadrant $\Re(K) > 0$, $\Im(K) > 0$. Hence there is only one other complex eigenvalue with modulus equal to that of $K$, and this is simply $\overline{K}$, as solutions of (\ref{kappa_cond_O1}) occur in complex conjugate pairs. Therefore, considering in turn the contributions to (\ref{unSoln}) from real and complex $\kappa$ and retaining only the dominant part of each, we have
\begin{equation}
\label{Const3_un_DomTerms}
u_n(Z) \sim \fr{ (-1)^n \Gamma(2n+4) \Lambda_{1,\psi} G(Z;\zeta_m) }{ \lb[ 2\pi\sqrt{m_1^2+m_2^2} (Z-\zeta_m) \rb]^{2n+4} } + \fr{(-1)^n \Gamma(2n+1)}{ (Z-\zeta_m)^{2n+1} } \lb( \fr{\Omega_{K,\psi}}{K^{2n+1}} + \fr{\Omega_{\overline{K},\psi}}{{\overline{K}}^{2n+1}} \rb).
\end{equation}
Note that (\ref{Const3_un_DomTerms}) is not meant to represent a two-term asymptotic series, as there may be many more terms in (\ref{unSoln}) which are much smaller than one of those on the right-hand side of (\ref{Const3_un_DomTerms}), but much larger than the other. Now, from (\ref{u0_Const3}) we have
\begin{equation}
\label{Const3_NearSing}
u_0 \sim  \fr{\sqrt{2}}{Z-\zeta_m}
\end{equation}
as $Z\rightarrow\zeta_m$; thus $U_0 = \sqrt{2}$, and
\begin{equation}
G(Z;\zeta_m) \sim -\fr{1}{8} \sqrt{2} \lb(Z-\zeta_m\rb)^{3}
\end{equation}
in the same limit. Because $U_0$ is real, inspection of (\ref{RecurrenceRelation}) indicates that $U_n$ is real for all $n$; hence $\Omega_{\overline{K},\psi} = \overline{\Omega}_{K,\psi}$. Comparing (\ref{Const3_un_DomTerms}) with (\ref{Const3_NearSing}), the dominant contributions to $U_n$ from real and complex $\kappa$ are
\begin{align}
\label{Const3_Un_DomTerms}
U_n \sim &\ \fr{ (-1)^{n+1} \sqrt{2} \Gamma(2n+4) \Lambda_{1,\psi} }{ 8 \lb( 2\pi\sqrt{m_1^2+m_2^2} \rb)^{2n+4} } 
\nonumber\\
& {} + \fr{2 (-1)^n \Gamma(2n + 1) |\Omega_{K,\psi}|}{|K|^{2n+1}} \cos \lb[ \Arg\lb(\Omega_{K,\psi}\rb) - \lb( 2n + 1 \rb) \Arg(K) \rb]
\end{align}
as $n\rightarrow \infty$. As discussed after (\ref{unSoln}), if $|K| < 2\pi (m_1^2+m_2^2)^{1/2}$ then the second term dominates; otherwise, the first does. Immediately we see a difficulty in obtaining $\Lambda_{1,\psi}$. If the second term is dominant, rearranging (\ref{Const3_Un_DomTerms}) provides an expression for $\Omega_{K,\psi}$, whereas we require $\Lambda_{1,\psi}$. To obtain $\Lambda_{1,\psi}$ in this way, the first term must be the dominant one. 

\begin{figure}[th!]
\centering
\includegraphics[width=\textwidth]{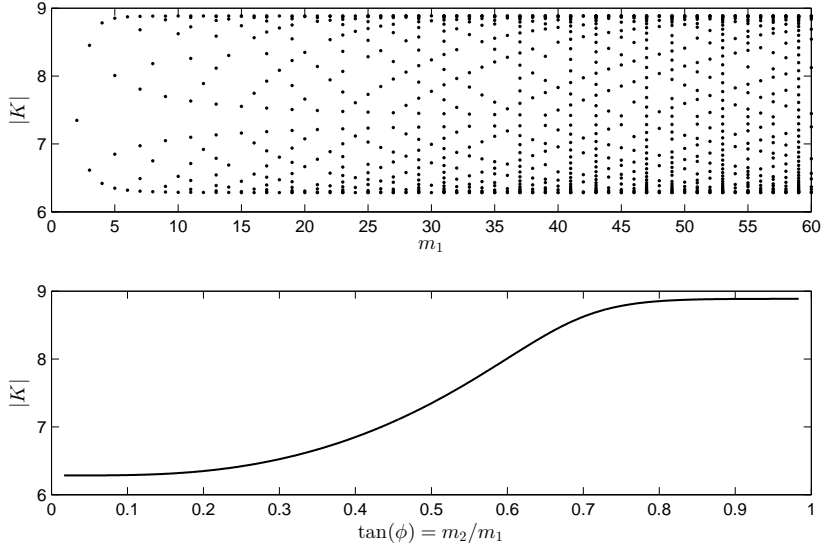}
\caption{Top: the complex solutions of (\ref{kappa_cond_O1}) having smallest modulus, for $m_1 = 2,\ldots,60$ and all values of $m_2$ satisfying $m_2 \leq m_1$, $\mathrm{gcd}(|m_1|,|m_2|) = 1$, plotted against $m_1$. Values of $m_2$ are not indicated. Bottom: data from the top figure, plotted against $\tan\psi$. Although the data are discrete, we employ a line plot for clarity.}
\label{Pic_ModComplexEigs}
\end{figure}

Unfortunately, it seems that if (\ref{kappa_cond_O1}) admits complex solutions, then \linebreak \mbox{$|K| < 2\pi (m_1^2+m_2^2)^{1/2}$} no matter the choice of $m_1, m_2$. Although we are unable to prove this, two strands of investigation provide evidence that this is indeed the case. Without loss of generality, we focus on the sector $\psi \in [0,\fr{\pi}{4}]$, so that $\tan\psi \in [0,1]$, equivalent to $m_2 \leq m_1$ with $m_1 \geq 1$ and $m_2 \geq 0$. Our results can be applied to the rest of the plane via the invariant rotations and reflections of (\ref{GenDiffEq}). First, we seek asymptotic solutions near $\psi = 0$ and $\psi = \fr{\pi}{4}$. The former is given by the limit $m_2 \ll m_1$, in which the complex solution to (\ref{kappa_cond_O1}) with smallest modulus is $K \sim 2\pi (1+m_2^2/m_1^2)^{1/2} ( 1 + i m_2/m_1)$; the latter corresponds to $m_2 \sim m_1$, in which $K \sim \pi (1+m_2^2/m_1^2)^{1/2} [ 2 + (1+i)(1-m_2/m_1)]$. In both limits we have $|K| < 2\pi(m_1^2+m_2^2)^{1/2}$, the real eigenvalue with smallest modulus. Second, solving (\ref{kappa_cond_O1}) numerically for $m_1 = 2,\ldots,60$ and all relevant values of $m_2$ in $\psi \in [0,\fr{\pi}{4}]$ (recall $\mathrm{gcd}(|m_1|,|m_2|) = 1$) has not produced a counter-example, as shown in figure \ref{Pic_ModComplexEigs}. Note that $|K|$ approaches $2\pi$ as $\tan\psi \rightarrow 0$, and $2\sqrt{2}\pi$ as $\tan\psi \rightarrow 1$, as predicted by our asymptotic solutions, and $K$ lies between these two limiting values for all choices of $\tan\psi$ for which solutions have been calculated. Thus it would appear that we cannot calculate $\Lambda_{1,\psi}$ using the above method if (\ref{kappa_cond_O1}) admits complex solutions. Furthermore, as the eigenvalues $\kappa$ are independent of the choice of $F(u;r)$, this is so for all problems of the form (\ref{GenDiffEq}).

There are, however, special cases with no complex eigenvalues at all; the axial and diagonal alignments $\psi = \fr{k\pi}{4}$, $k \in \{0,1,\ldots,7\}$. For such alignments (\ref{Const3_Un_DomTerms}) contains only the term in which $\Lambda_{1,\psi}$ appears. Rearranging, we therefore see that in such a case
\begin{equation}
\label{LambdaEqn}
\Lambda_{1,\psi} \sim \lim_{n\rightarrow\infty} \fr{(12)^{1/4} (-1)^{n+1} \lb(2\pi\sqrt{m_1^2+m_2^2}\rb)^{2n+4}}{\Gamma(2n+4)}U_n, 
\end{equation}
yielding a good approximation for $\Lambda_{1,\psi}$ if $U_n$ can be calculated for large enough $n$. Now, (\ref{RecurrenceRelation}) must in general be iterated separately for each $\psi$. However, as (\ref{Const3DiffEq}) is invariant under rotations $\psi \rightarrow \psi + \fr{\pi}{2}$, it suffices to iterate (\ref{RecurrenceRelation}) only for $\psi = 0$ and $\psi = \fr{\pi}{4}$, as the other six alignments can be recovered using said invariance. Doing so, we calculate $\Lambda_{1,0} \approx -2535$ and $\Lambda_{1,\pi/4} \approx -10141$. Thus we may carry out a quantitative comparison between the analytical formula (\ref{Const3_SnakeWidth_Unscaled}) and numerical computations for these values of $\psi$. We note that inspection of the ratio $\Lambda_{1,\pi/4}/\Lambda_{1,0} \approx 4$ and the equivalent ratio in the next example suggests the simple relationship $\Lambda_{1,\psi} = \Lambda_{1,0}(m_1^2+m_2^2)^{1+\gamma}$ (recall that $\gamma$ is the strength of the leading order singularity; cf. (\ref{AlgebraicSingularity})-(\ref{BetaCalculation}) and the surrounding paragraph). However, this is found to drastically underestimate the width of the pinning region for $\psi \neq \fr{k\pi}{4}$, $k \in \{0,1,\ldots,7\}$.


\subsubsection{Comparison of analytical and numerical results}

We have solved the one-dimensional, steady version of (\ref{Const3DiffEq}) for $\psi=0, \fr{\pi}{4}$ numerically, using pseudo-arclength continuation to compute the bifurcation diagram. The domain size is chosen to be large enough that boundaries have negligible effect on the width of the pinning region. In order to preserve this independence, the domain must be increased as $\hat{s} = \epsilon s$ decreases to counterbalance the spreading out of fronts; for example, we used three hundred points for $\hat{s}=1$, but seven hundred for $\hat{s}=0.2$. We imposed symmetric boundary conditions and sought stationary front solutions to (\ref{Const3DiffEq}); such solutions are equivalent to site-centred solutions on a domain of twice the size. Exploiting the symmetry of the solution to use only half the lattice points in this manner has the dual benefit of faster computation times and a significantly decreased chance of skipping between solution branches, which may otherwise occur all too readily given the high density of solutions within such a narrow parameter range \cite{burke2009multipulse}. We have chosen to focus here only on site-centred solutions; similar results are easy to find for the bond-centred solution branch. Of course, there is no symmetry to exploit when computing the ladders and so these must be found on the full domain, hence requiring great care at small values of $\hat{s}$.

\begin{figure}[th!]
\centering
\includegraphics[trim = 0cm 0 1cm 0, clip, width=0.5\textwidth]{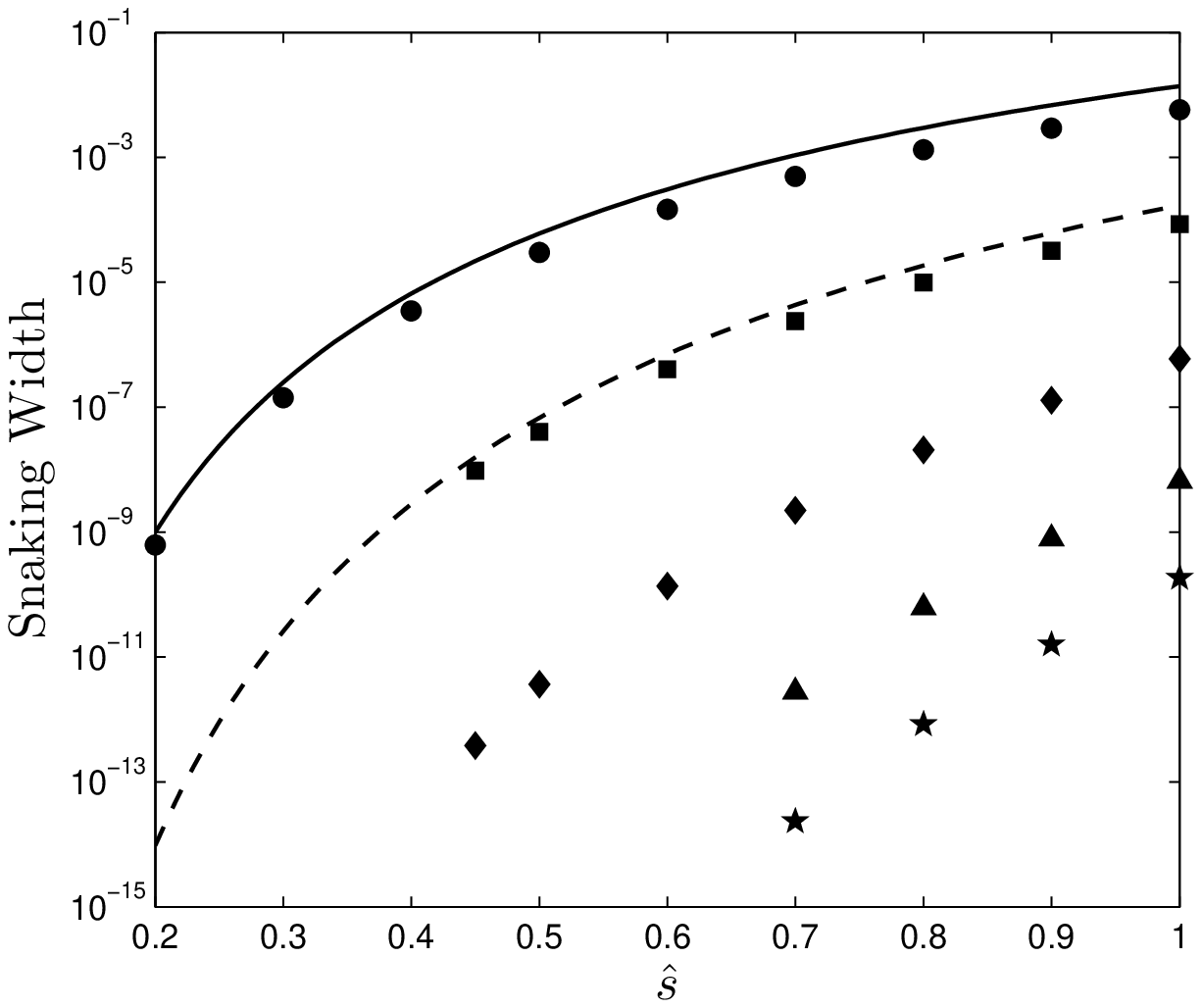}\includegraphics[trim = 0cm 0 1cm 0, clip, width=0.5\textwidth]{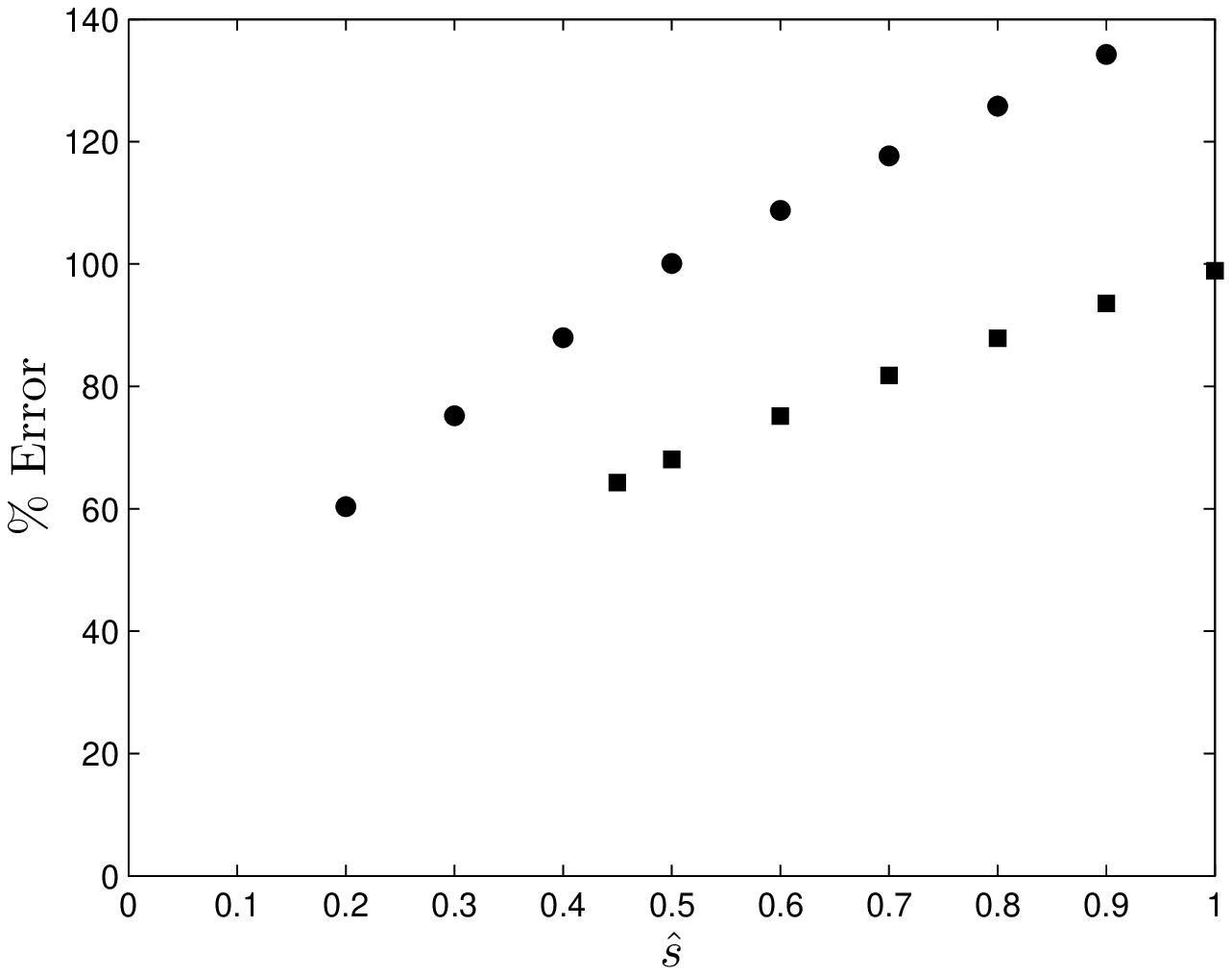}
\caption{Left: analytical (lines) and numerical (points) snaking widths for one-dimensional solutions to (\ref{Const3DiffEq}) at various orientations $\psi = \arctan (m_2/m_1)$. The solid line and circles represent $(m_1,m_2) = (1,0)$; the dashed line and squares represent $(m_1,m_2) = (1,1)$; diamonds represent $(m_1,m_2) = (2,1)$; triangles represent $(m_1,m_2) = (3,1)$; stars represent $(m_1,m_2) = (3,2)$. Note that an analytical formula is only available for the first two choices of $\psi$. Right: percentage error in analytical formula (\ref{Const3_SnakeWidth_Unscaled}) for (\ref{Const3DiffEq}) with $(m_1,m_2) = (1,0)$ represented by circles and $(m_1,m_2) = (1,1)$ by squares.}
\label{Pic_Const3CompareNumerics}
\end{figure}

\begin{figure}[th!]
\centering
\includegraphics[trim = 0cm 0 1cm 0, clip, width=0.5\textwidth]{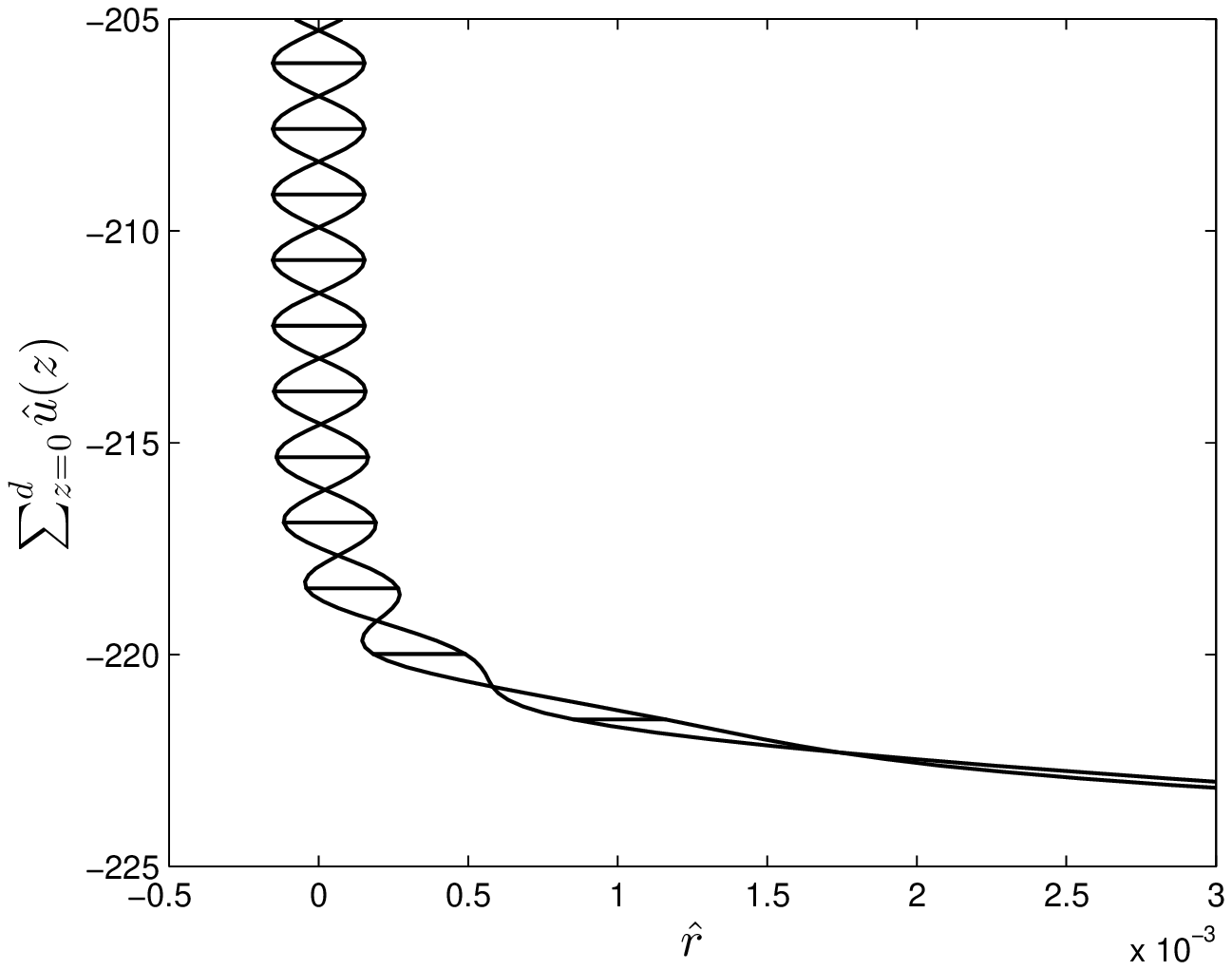}\includegraphics[trim = 0cm 0 1cm 0, clip, width=0.5\textwidth]{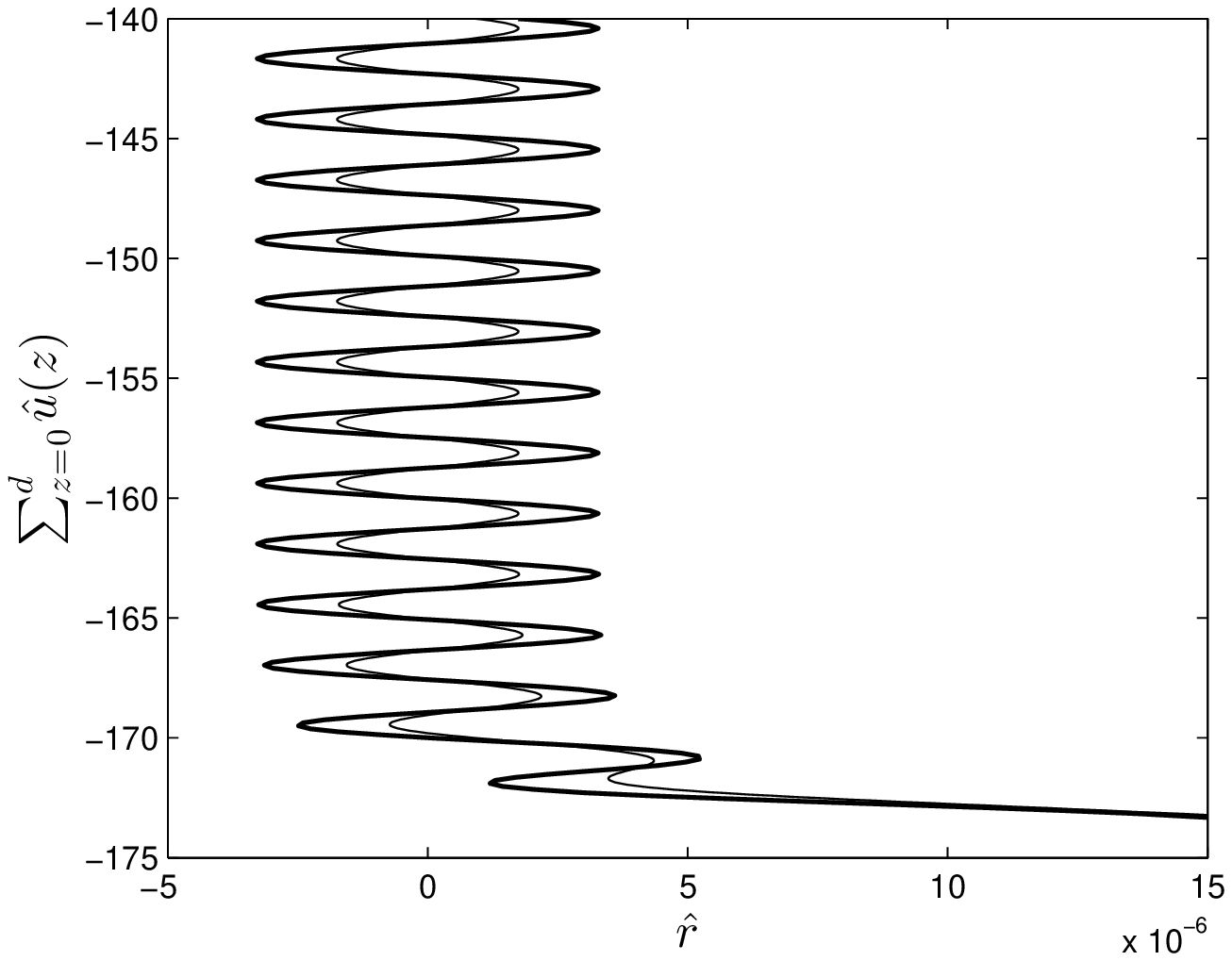}
\caption{Left: snakes-and-ladders bifurcation diagram for (\ref{Const3DiffEq}) with $\hat{s} = 0.6$ and $\psi=0$, drawn using the analytical formulae (\ref{SnakesCondition})-(\ref{GeneralLadderEq}). Right: comparison of analytical (thick line) and numerical (thin line) site-centred snaking curve for $\hat{s}=0.4$ and $\psi=0$.}
\label{Pic_Const3BifDiags}
\end{figure}

Numerical results are compared to (\ref{Const3_SnakeWidth_Unscaled}) in figures \ref{Pic_Const3CompareNumerics}, with good agreement. Although an analytical formula is unavailable for $\psi \neq \fr{k\pi}{4}$, $k \in \{0,1,\ldots,7\}$ as $\Lambda_{1,\psi}$ remains undetermined in these cases, we see that the snaking width appears to scale with $\hat{s}$ as predicted by (\ref{Const3_SnakeWidth_Unscaled}) for all values of $\psi$ shown. Note that numerical results are unavailable for very small $\hat{s}$, and that the smallest value of $\hat{s}$ for which numerical results are available increases with $m_1^2+m_2^2$; this is due to the snaking width approaching values in which machine error is significant. The full analytical bifurcation diagram (\ref{GeneralSnakeEq})-(\ref{GeneralLadderEq}) for (\ref{Const3DiffEq}) with $\psi = 0$ is drawn in figure \ref{Pic_Const3BifDiags}a, and a comparison between an analytical and a numerical snaking solution curve shown in figure \ref{Pic_Const3BifDiags}b, again with good agreement.

\subsection{A cubic-quintic nonlinearity}
\label{Sec_35Example}

Our second example is
\begin{equation}
\label{35DiffEq}
\fr{\pd \hat{u}}{\pd t} = \Delta \hat{u} + \hat{r}\hat{u} + \hat{s}\hat{u}^3 - \hat{u}^5,
\end{equation}
which is bistable when $\hat{s} > 0$. (\ref{35DiffEq}) was the subject of a numerical investigation of snaking of fully two-dimensional localisations in \cite{taylor2010snaking}, and is analogous to the Swift-Hohenberg equation with cubic and quintic nonlinear terms \cite{burke2007snakes, dean2011exponential}. The bifurcation diagram and example solutions for (\ref{35DiffEq}) with $\psi=0$ can be seen in figure \ref{Pic_35ExampleBifDiag}; unhatted variables in that figure correspond to hatted in (\ref{35DiffEq}). The system is bistable due to a subcritical pitchfork bifurcation at $\hat{r}=0$ and a subsequent saddle-node bifurcation at some $\hat{r}<0$, at which point the nontrivial solution curve turns over to form a region of bistability, containing the pinning region. 

Rescaling in the weakly subcritical limit, we define $(\hat{u}, \hat{r}, \hat{s}) = (\sqrt{\epsilon} u, \epsilon^2 r, \epsilon s)$, yielding
\begin{equation}
\label{35DiffEq_Scaled}
\fr{\pd u}{\pd t} = \Delta u - \epsilon^2 \lb( -ru - su^3 + u^5 \rb).
\end{equation}
Hence $F(u;r) = -ru-su^3+u^5$. The leading-order, one-dimensional solution is therefore given by
\begin{equation}
\label{35DiffEq_LOApprox}
0 = \fr{\od^2u_0}{\od Z^2} + ru_0 + su_0^3 - u_0^5.
\end{equation}
Imposing (\ref{rM_integral}), this exhibits the front solution
\begin{equation}
\label{35DiffEq_Front}
u_0(Z) = \fr{1}{2} \lb( \fr{3s}{1 + e^{-\sqrt{3} s Z/2}} \rb)^{1/2}.
\end{equation}
at the Maxwell point $r_M = -3s^2/16$. Thus $u \rightarrow u_\pm $ as $Z \rightarrow \pm\infty$, where $u_+ = \sqrt{3s}/2$ and $u_- = 0$. We have chosen the positive square root in order that $u_+>u_-$; the front of opposite orientation may be recovered by exploiting the reversibility of (\ref{35DiffEq}).

From (\ref{35DiffEq_Front}), we can see that the singularities $\zeta$ of $u_0$ are
\begin{equation}
\zeta = \zeta_m := (2m+1)2\pi i / \sqrt{3}s, \qquad m\in\mathbb{Z},
\end{equation}
each of which has strength $\gamma = \fr{1}{2}$. Thus (\ref{chi_eq}) gives $\chi = \Arg\lb(\Lambda_{1,\psi}\rb)$. The dominant singularities are those nearest (and equidistant from) the real line, namely $\zeta_0$ and $\zeta_{-1}=\overline{\zeta}_0$. Also, since
\begin{equation}
u_0 \sim \fr{1}{2} \sqrt{3s} \lb( 1 - \fr{1}{2} e^{-\sqrt{3} s Z/2} \rb)
\end{equation}
as $Z \rightarrow \infty$ and $F_{r,M}(u) \equiv -u$, we have
\begin{equation}
\alpha_+ = \fr{1}{2}\sqrt{3}s, \qquad D_+ = \fr{1}{4} \sqrt{3s}, \qquad \int_{u_-}^{u_+} F_{r,M}(v) \od v = - \fr{3}{8}s.
\end{equation}
Thus the bifurcation equations (\ref{GeneralSnakeEq}) and (\ref{GeneralLadderEq}) can now be written in terms of the parameters of the scaled equation (\ref{35DiffEq_Scaled}). Again, we shall not write these out in full; instead we simply write down the width of the pinning region from (\ref{DiscreteSnakingWidth}), which now reads
\begin{equation}
|\delta r| \leq \fr{16\pi \lb| \Lambda_{1,\psi} \rb| e^{-4\pi^2 \sqrt{3(m_1^2+m_2^2)}/ 3\epsilon s}}{3\epsilon^3s}.
\end{equation}
Absorbing the scalings in $\epsilon$ and writing this in terms of the original, hatted variables appearing in (\ref{35DiffEq}) then yields
\begin{equation}
\label{35_SnakeWidth_Unscaled}
\lb| \hat{r} - \hat{r}_M \rb| \leq \fr{16\pi \lb| \Lambda_{1,\psi} \rb| e^{-4\pi^2 \sqrt{3(m_1^2+m_2^2)}/ 3 \hat{s}}}{3\hat{s}},
\end{equation}
where $\hat{r}_M = - 3\hat{s}^2/16 + \Or(\hat{s}^4)$ is the unscaled Maxwell point and $\hat{s}$ provides the small variable. This formula corresponds to that derived in \cite{matthews2011variational} using variational approximations (equation (50) in that work); however, the method presented here yields a complete formula, whereas that in \cite{matthews2011variational} is unable to determine the constant factor \linebreak $16\pi |\Lambda_{1,\psi}|/3$. We also note that the functional dependence of (\ref{35_SnakeWidth_Unscaled}) on $\hat{s}$ when $(m_1,m_2) = (1,0)$ is identical to that of the corresponding formula derived in \cite{dean2011exponential} for the cubic-quintic Swift-Hohenberg equation (equation (8.7) in that work). However, the snaking width is much smaller in the present case, as $e^{-1/\hat{s}}$ is raised to the power $\approx 22.8$  in (\ref{35_SnakeWidth_Unscaled}) when $\psi = 0$, and only to the power $\approx 15.3$ in the equivalent formula for the Swift-Hohenberg equation. 

All that remains is to derive the constants $\Lambda_{1,\psi}$. In a similar manner as in Section (\ref{Sec_Const3FindLambda}), we have
\begin{equation}
u_n \sim \fr{U_n}{(Z-\zeta_m)^{2n+1/2}},
\end{equation}
as $Z\rightarrow\zeta_m$, for some sequence of constants $U_n$. These can be found by iteration of the recurrence relation
\begin{align}
\label{RecurrenceRelation35}
0 = &\ 2 \sum_{p=1}^{n+1} \fr{\lb( \cos^{2p}\psi + \sin^{2p}\psi \rb)\Gamma\lb( 2n+\tfrac{5}{2} \rb)}{(2p)!\Gamma\lb( 2n-2p+\tfrac{5}{2} \rb)}  U_{n-p+1} 
\nonumber\\
& {} - \sum_{p_1=0}^n \sum_{p_2=0}^{n-p_1} \sum_{p_3=0}^{n-p_1-p_2} \sum_{p_4=0}^{n-p_1-p_2-p_3} U_{p_1} U_{p_2} U_{p_3} U_{p_4} U_{n-p_1-p_2-p_3-p_4},
\end{align}
where $U_0 = \lb( \fr{3}{4} \rb)^{1/4}$. Again, due to the dominant contribution from complex eigenvalues at other values of $\psi$, we are able to calculate $\Lambda_{1,\psi}$ only when $\psi = \fr{k\pi}{4}$, $k \in \{0,1,\ldots,7\}$. For such orientations
\begin{equation}
\Lambda_{1,\psi} \sim \lim_{n\rightarrow\infty} \fr{(12)^{1/4} (-1)^{n+1} \lb(2\pi\sqrt{m_1^2+m_2^2}\rb)^{2n+3}}{\Gamma(2n+3)}U_n.
\end{equation}
By iteration of (\ref{RecurrenceRelation35}), we are therefore able to calculate $\Lambda_{1,0} \approx -89$ and $\Lambda_{1,\pi/4} \approx -252$; $\Lambda_{1,\psi}$ for $\psi = \fr{k\pi}{4}$, $k \in \{2,3,\ldots,7\}$ then follow using the invariance of (\ref{35DiffEq}) under rotations $\psi \rightarrow \psi + \fr{\pi}{2}$.

\begin{figure}[th!]
\centering
\includegraphics[trim = 0cm 0 1cm 0, clip, width=0.5\textwidth]{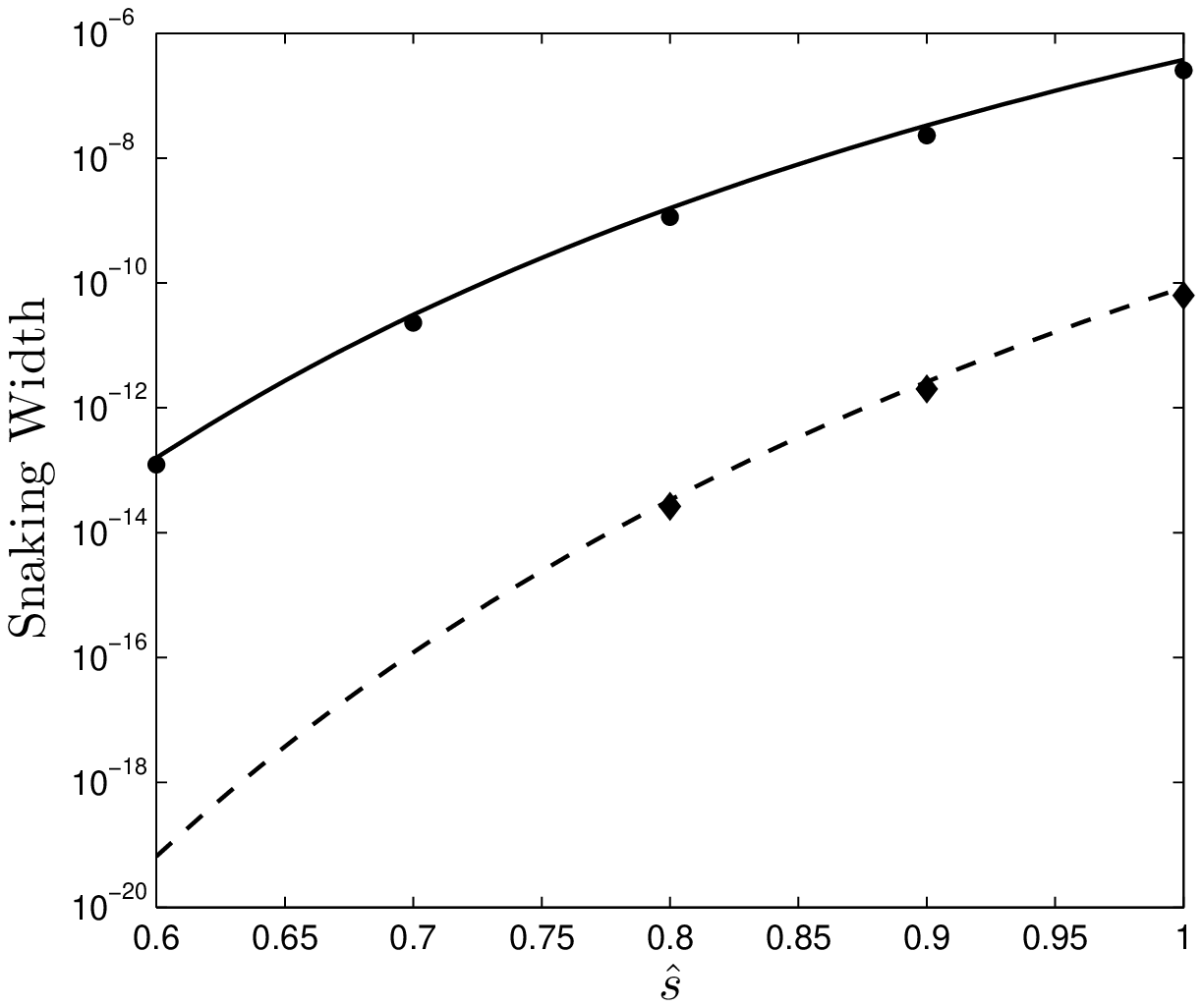}\includegraphics[trim = 0cm 0 1cm 0, clip, width=0.5\textwidth]{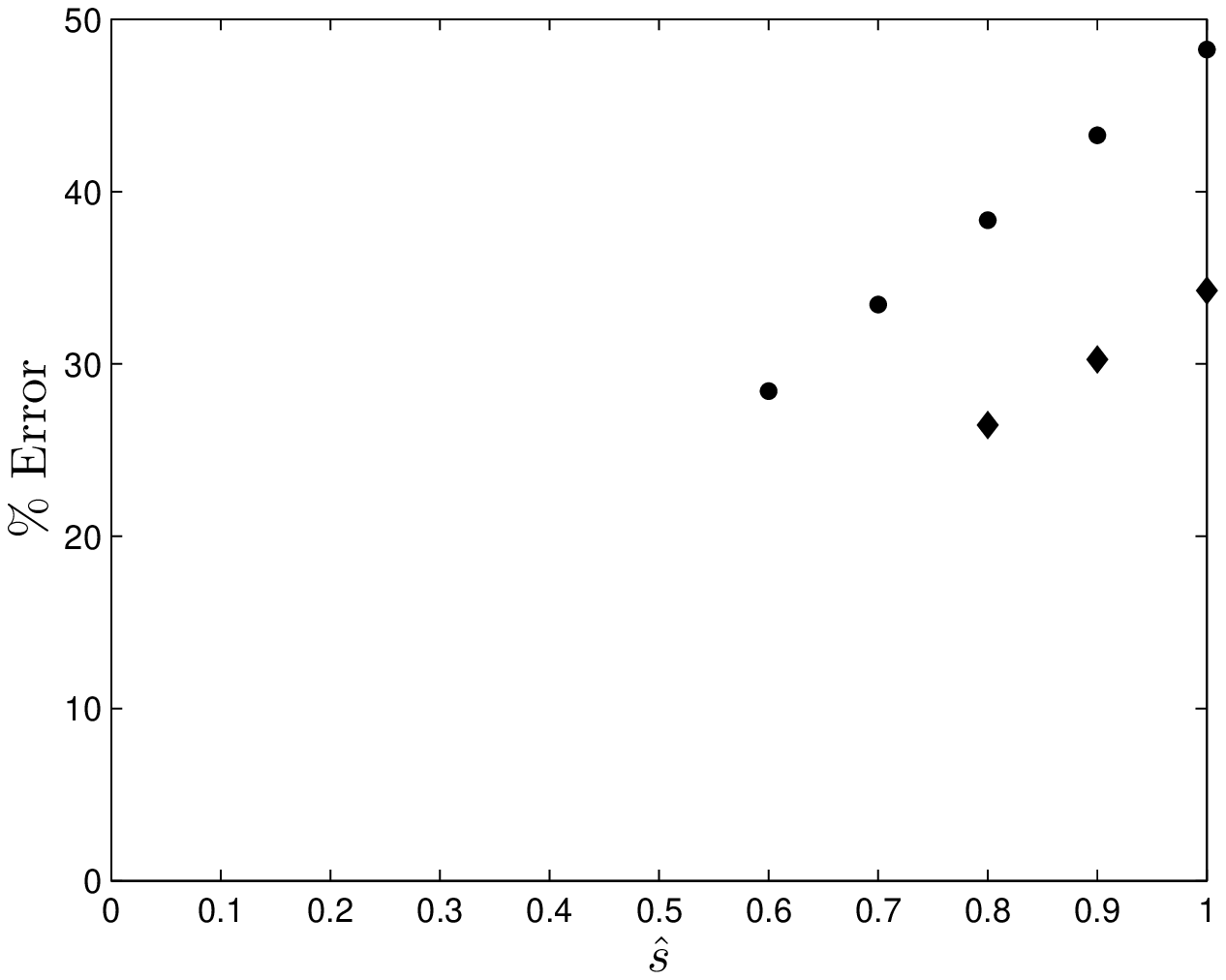}
\caption{Left: comparison of analytical (lines) and numerical (data points) snaking widths for (\ref{35DiffEq}). Right: percentage error in analytical formula compared to numerical results. Solid lines or circles correspond to $\psi=0$, and dashed lines or diamonds to $\psi=\fr{\pi}{4}$.}
\label{Pic_CompareSnakeWidths35}
\end{figure}

\begin{figure}[th!]
\centering
\includegraphics[trim = 0cm 0 1cm 0, clip, width=0.5\textwidth]{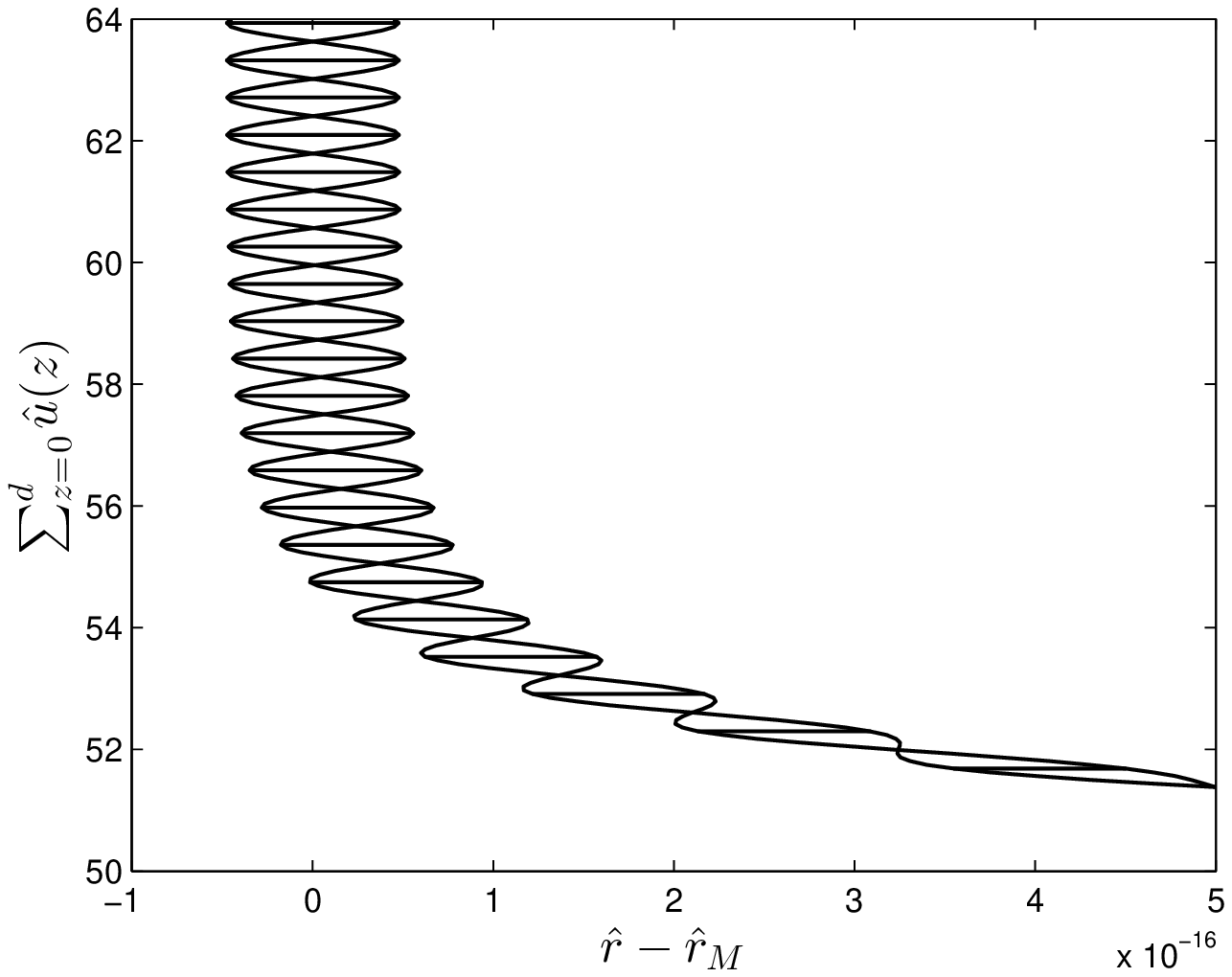}\includegraphics[trim = 0cm 0 1cm 0, clip, width=0.5\textwidth]{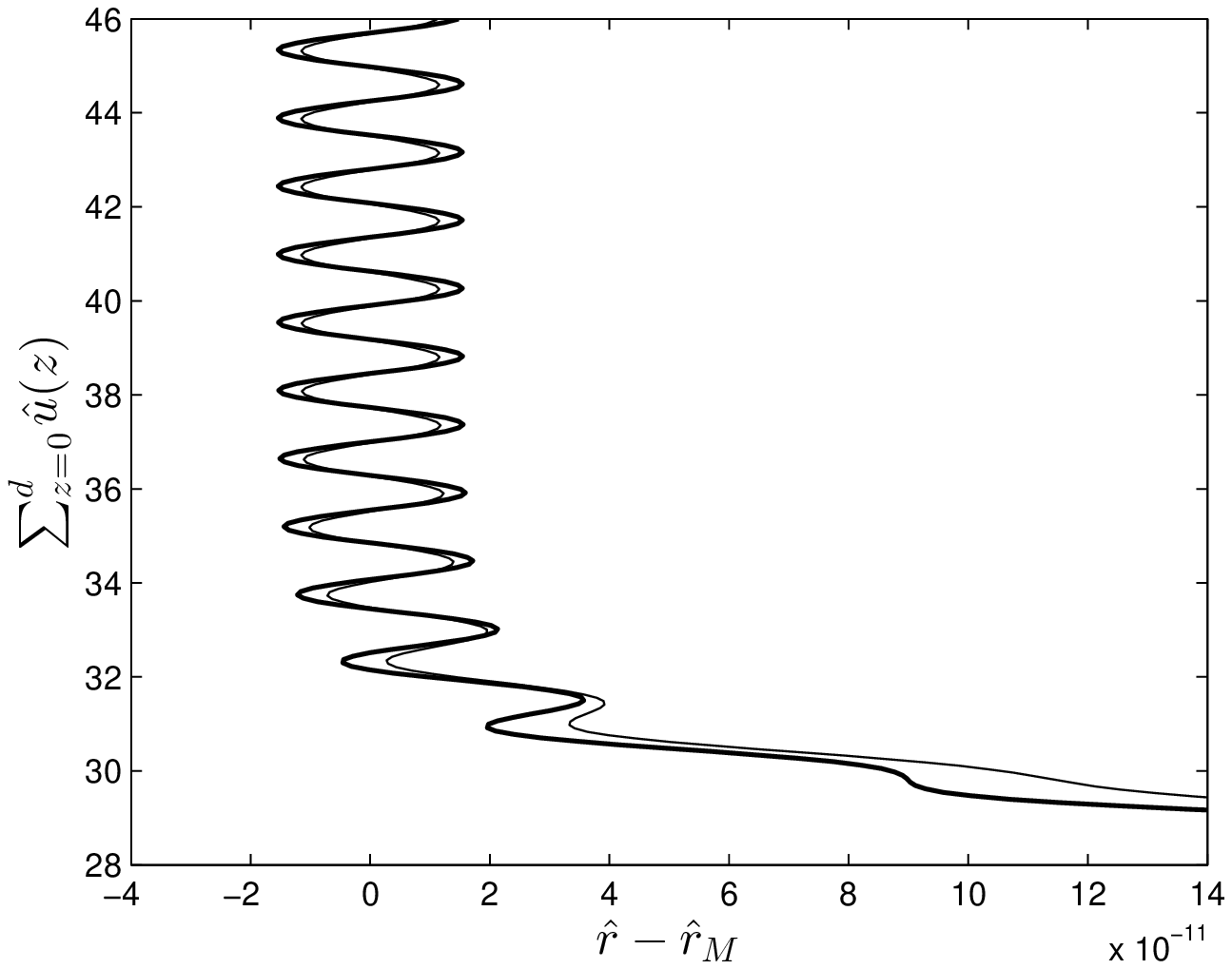}
\caption{Left: snakes-and-ladders bifurcation diagram for (\ref{35DiffEq}) with $\hat{s} = 0.5$ and $\psi=0$, drawn using the analytical formulae (\ref{SnakesCondition})-(\ref{GeneralLadderEq}). Right: comparison of analytical (thick line) and numerical (thin line) site-centred snaking curve for $\hat{s}=0.7$ and $\psi=0$.}
\label{Pic_BifDiags35}
\end{figure}

Numerical computations for $\psi=0,\fr{\pi}{4}$ are compared to (\ref{35_SnakeWidth_Unscaled}) in figure \ref{Pic_CompareSnakeWidths35}, with good agreement. Note that machine error becomes significant at much larger values of $\hat{s}$ than in the example of Section \ref{Sec_Constant3Example}; this is because the exponent in (\ref{35_SnakeWidth_Unscaled}) is more negative than that in (\ref{Const3_SnakeWidth_Unscaled}). The full analytical bifurcation diagram is drawn in figure \ref{Pic_BifDiags35}a using the value of $\Lambda_{1,0}$ calculated from the recurrence relation (\ref{RecurrenceRelation}), and a comparison between an analytical and a numerical snaking solution curve shown in figure \ref{Pic_BifDiags35}b.

\section{Snaking on a hexagonal lattice}
\label{Sec_Hex}

The methods of Sections \ref{Sec_Preliminaries}-\ref{Sec_SnakeBifEqs} can readily be followed through for an alternative choice of difference operator. As a useful example, we shall now consider (\ref{GenDiffEq}) posed on a hexagonal lattice, in which case we replace the square operator $\Delta$ (\ref{Delta}) with the hexagonal operator
\begin{align}
\label{Delta_hex}
\Delta_{\mathrm{hex}} u(x,y,t) := &\ \fr{2}{3} \lb[ u(x+1,y,t) + u(x-1,y,t) + u\lb( x+\fr{1}{2}, y+\fr{\sqrt{3}}{2},t \rb) \rb.
\nonumber\\
& \lb. {} + u\lb( x-\fr{1}{2}, y-\fr{\sqrt{3}}{2},t \rb) + u\lb( x+\fr{1}{2}, y-\fr{\sqrt{3}}{2},t \rb) \rb.
\nonumber\\
& \lb. {} + u\lb( x-\fr{1}{2}, y+\fr{\sqrt{3}}{2},t \rb) - 6u(x,y,t) \rb],
\end{align}
which generates the hexagonal lattice
\begin{equation}
\label{HexagonalLattice}
H := \lb\{ \lb. (x,y) = \lb(n_1+\fr{n_2}{2},\fr{\sqrt{3}n_2}{2} \rb) \ \rb| \ (n_1,n_2) \in \mathbb{Z}^2 \rb\}.
\end{equation}
As is the case for (\ref{GenDiffEq}) on a square lattice, (\ref{GenDiffEq}) on $H$ is the discrete analogue of the reaction-diffusion equation (\ref{RDEq}) (i.e. both $\Delta$ and $\Delta_{\mathrm{hex}}$ are discrete analogues of $\pd_X^2+\pd_Y^2$). Note that this means the leading-order front solution and corresponding singularities are unaffected by replacing $\Delta$ with $\Delta_{\mathrm{hex}}$, i.e. $u_0$ is determined by (\ref{FrontEq}) for both choices of difference operator. In fact, following through the calculations in Sections \ref{Sec_Preliminaries}-\ref{Sec_SnakeBifEqs}, we can see that choosing a hexagonal lattice only affects the orientations $\psi$ for which snaking occurs, and the values of the corresponding eigenvalues $\kappa$ (cf. Section \ref{Sec_RemEqn}). Therefore, in order to apply our results to (\ref{GenDiffEq}) on $H$, we simply need to solve the appropriate eigenvalue equation for $\Delta_{\mathrm{hex}}$ and replace $\kappa$ and $\psi$ in (\ref{DiscreteSnakingWidth}) and (\ref{SnakesCondition})-(\ref{LaddersCondition}) with their corresponding values for a hexagonal lattice.

Referring to Section \ref{Sec_RemEqn}, in particular (\ref{kappa_cond_O1}), we see that the eigenvalues $\kappa$ are given by $\Delta_{\mathrm{hex}} e^{i\kappa z} \equiv 0$. This yields
\begin{equation}
\label{kappa_eqn_hex}
\cos\lb( \kappa\cos\psi \rb) + \cos\lb( \fr{\kappa}{2}\cos\psi + \fr{\sqrt{3}\kappa}{2}\sin\psi \rb) + \cos\lb( \fr{\kappa}{2}\cos\psi - \fr{\sqrt{3}\kappa}{2}\sin\psi \rb) - 3 = 0.
\end{equation}
In the same way as for the square lattice (cf. Sections \ref{Sec_Rotation} and \ref{Sec_SnakeWidth}), (\ref{GenDiffEq}) on $H$ exhibits snaking when $\psi$ is such that (\ref{kappa_eqn_hex}) admits real solutions, i.e. $\kappa \in \mathbb{R}$. This requires $\kappa\cos\psi = 2M_1\pi$, $\kappa\lb( \cos\psi + \sqrt{3}\sin\psi \rb) = 4M_2\pi$ and $\kappa\lb( \cos\psi - \sqrt{3}\sin\psi \rb) = 4M_3\pi$, where the $M_j$ are arbitrary integers. After some manipulation, we arrive at the solution
\begin{equation}
\label{tanphi_hex}
\begin{split}
\begin{array}{clcl}
& \cos\psi = \dfrac{\sqrt{3}}{2}\dfrac{m_1}{\sqrt{m_1^2+m_1m_2+m_2^2}}, & \qquad & \sin\psi = \dfrac{1}{2} \dfrac{m_1+2m_2}{\sqrt{m_1^2+m_1m_2+m_2^2}}, 
\\
\vspace{-1ex}
\\
& (m_1,m_2) \in \mathbb{Z}^2\backslash\{(0,0)\}, & \qquad & \mathrm{gcd}(|m_1|,|m_2|) = 1,
\end{array}
\end{split}
\end{equation}
and
\begin{equation}
\label{kappa_hex}
\kappa = \fr{4M\pi}{\sqrt{3}} \sqrt{m_1^2+m_1m_2+m_2^2}, \qquad M \in\mathbb{Z}.
\end{equation}
Hence $\tan\psi = (m_1+2m_2)/\sqrt{3}m_1$, i.e. $(\sqrt{3}\tan\psi-1)/2 \in \mathbb{Q}_\infty$. In light of (\ref{HexagonalLattice}), when $\psi$ is given by (\ref{tanphi_hex}) and $\kappa$ by (\ref{kappa_hex}), the one-dimensional lattice spanned by $z$ (cf. Section \ref{Sec_Rotation}) is
\begin{equation}
\Psi = \lb\{ \lb. \fr{\sqrt{3}}{2} \fr{m_1n_1+m_1n_2+m_2n_2}{\sqrt{m_1^2+m_1m_2+m_2^2}} \ \rb| \ (n_1,n_1) \in \mathbb{Z}^2 \rb\},
\end{equation}
and therefore has a well-defined lattice spacing of $(\sqrt{3}/2)(m_1^2+m_1m_2+m_2^2)^{-1/2}$; note that $\kappa$ is given by $2M\pi$ divided by the effective lattice spacing, as is the case on a square lattice (cf. (\ref{realkappa})).

We note that the lattice generated by $\kappa(\cos\psi,\sin\psi)$ is also hexagonal, but oriented at $\pi/6$ radians to $H$. This is because $\kappa(\cos\psi,\sin\psi)$ is the reciprocal lattice of $H$, defined by $\exp[\pi\kappa(\cos\psi,\sin\psi)\cdot(x,y)] \equiv e^{i\kappa z} \equiv 1$ for $(x,y) \in H$. This is, of course, analogous to the condition (\ref{tanphi_realkappa}) for $\kappa \in \mathbb{R}$ on a square lattice derived in Section \ref{Sec_RemEqn}. However, the importance of the reciprocal lattice is not immediately apparent in Section \ref{Sec_RemEqn}, as in that case it is also a square lattice with the same orientation as that of $\Delta$. Thus we can make make the more general statement that snaking will occur only if $\kappa(\cos\psi,\sin\psi)$ is a vector of the reciprocal lattice to that the problem is posed on. 

Following the calculation through to its conclusion, we eventually arrive at a formula for the width of the snaking region, namely
\begin{equation}
\label{DiscreteSnakingWidth_hex}
|\delta r| \leq \fr{2\pi|\Lambda_{1,\psi}| e^{-4\pi \sqrt{m_1^2++m_1m_2+m_2^2} \Im(\zeta)/\sqrt{3}\epsilon}}{\epsilon^{2\gamma+2} |\int_{u_-}^{u_+} F_{r,M}(v) \od v|}.
\end{equation}
(\ref{DiscreteSnakingWidth_hex}) is simply the snaking width for a square lattice (\ref{DiscreteSnakingWidth}) with $2\pi\sqrt{m_1^2+m_2^2}$ in the exponent replaced by the smallest real eigenvalue for the hexagonal lattice, (\ref{kappa_hex}) with $M=1$. Similarly, the full equations for the snakes-and-ladders bifurcation diagram are given by (\ref{SnakesCondition})-(\ref{LaddersCondition}) with all instances of $2\pi\sqrt{m_1^2+m_2^2}$ replaced by (\ref{kappa_hex}) with $M=1$. The constant $\Lambda_{1,\psi}$ can be determined by a process analogous to that of (\ref{Sec_Const3FindLambda}), but will differ in value to that for the square lattice since, for example, (\ref{un_NearSing}) depends on the choice of lattice. As in Section \ref{Sec_SnakeWidth}, we note that (\ref{DiscreteSnakingWidth_hex}) should be compared with the equivalent result in \cite{kozyreff2013analytical}. The quantity $|\Delta k|$ in \cite{kozyreff2013analytical} now corresponds to $4\pi\sqrt{m_1^2+m_1m_2+m_2^2}/\sqrt{3}$; $\pi/\lambda$ again corresponds to $\Im(\zeta)$ since the singularity $\zeta$ is the same for both choices of lattice.

We conclude this section by noting that the relatively simple means by which our results for the square lattice have been applied to the hexagonal lattice is due to the fact that both $\Delta$ (\ref{Delta}) and $\Delta_{\mathrm{hex}}$ are discrete analogues of the same differential operator, $\pd_X^2+\pd_Y^2$. Of course, both square and hexagonal lattices may be generated by other choices of difference operator, some of which give rise to leading-order continuum approximations containing differential operators other than $\pd_X^2+\pd_Y^2$. Thus more work is required in order to apply our results to such systems than has been necessary in the present section. However, our method is applicable to a wide class of difference operators, and in many cases our results can be adapted without any great difficulty.


\section{Conclusion}
\label{Sec_Conc}

We have applied the method of exponential asymptotics to the study of one-dimensional heteroclinic and homoclinic connections in the class of differential-difference equations given by (\ref{GenDiffEq}). By studying slowly varying solutions near bifurcation (equivalent to the continuum approximation in the limit of small mesh spacing) and truncating the asymptotic expansion after its least term, we have been able to elucidate the role played by the exponentially small remainder. Rescaling near Stokes lines, at which the remainder equation (\ref{RNEqn_unknownForcing}) is maximally forced, we have observed explicitly how the coefficient of an exponentially growing complementary function varies smoothly from zero to non-zero via an error-function. These Stokes lines emanate from complex singularities of the leading-order front. Furthermore, exponentially small deviations $\delta r$ from the Maxwell point $r_M$ also produce exponentially small, but exponentially growing, particular integrals. We saw that unbounded terms vanish only for particular values of the origin of the front, defined in terms of $\delta r$. This results in an exponentially small parameter range (\ref{DiscreteSnakingWidth}) in which stationary fronts, and hence localised solutions, exist.

Armed with the full asymptotic expansion of the front solution, localised solutions were then constructed by means of matching exponentially growing and decaying terms in two back-to-back fronts. Matching conditions yielded a set of formulae (\ref{SnakesCondition})-(\ref{GeneralLadderEq}) which describe the full snakes-and-ladders bifurcation diagram associated with such solutions.

Of particular interest is the result that the snaking width is non-zero only if $\tan\psi \in \mathbb{Q}_\infty$, and then is exponentially small in $(m_1^2+m_2^2)^{1/2}$. This is to be expected, because (when $\tan\psi \in \mathbb{Q}_\infty$) the effective lattice spacing is $(m_1^2+m_2^2)^{-1/2}$. For each value of $\epsilon$, then, the largest pinning region is for those solutions oriented along an axis $\psi = \fr{k\pi}{2}$, followed by those oriented along a primary diagonal $\psi = \fr{(2k+1)\pi}{4}$, where $k \in \{1,2,3,4\}$. The decrease in snaking width as $m_1$ and $m_2$ increase is considerable, as the effective lattice spacing appears in the exponentially small term in (\ref{DiscreteSnakingWidth}); this explains why those solutions oriented along an axis or primary diagonal are the easiest to find numerically in the small-$\epsilon$ limit (cf. figure \ref{Pic_Const3CompareNumerics}). On the other hand, when $\tan\psi$ is irrational the problem is posed on a dense set, and hence there is no periodic spatial structure for fronts to pin to and localised solutions do not snake \cite{hoffman2010universality, hupkes2011propagation, mallet2001crystallographic}.

Furthermore, the existence regions of the different fronts has implications for fully two-dimensional localisations. A numerical study \cite{taylor2010snaking} of (\ref{35DiffEq}) found that two-dimensional localised patches evolve in a rather complicated manner along the snaking curve, with saddle-nodes aligning to a number of asymptotes in phase space, in contrast to the two asymptotes of the one-dimensional case; compare the snaking diagrams in \cite{taylor2010snaking} with, for example, figures \ref{Pic_Const3BifDiags} and \ref{Pic_BifDiags35}. Inspection of the results presented in \cite{taylor2010snaking} indicates that square localised solutions with sides aligned with the axes are present in a wider parameter range than, for example, those which appear octagonal, having sides aligned with both the axes and the primary diagonals. This may be explained by interpreting two-dimensional localisations as constructed from superpositions of various one-dimensional fronts, since the fronts aligned with the axes have the widest pinning region. This apparent relationship between one-dimensional fronts and two-dimensional localisations bears further investigation.

Although the focus has been on snaking problems, the analysis is equally applicable to fronts and localised solutions which do not snake. In these instances, there is no Maxwell point and any reference to $\delta r$ is meaningless, as solutions exist within an $\Or(1)$ range of $r$ rather than an exponentially small one. However, they still pin to the lattice in the manner discussed in Section \ref{Sec_SnakeWidth}.

We remark that in \cite{king2001asymptotics}, on which the present work is in part based, the (purely one-dimensional) analysis included a description of front motion outside the pinning region; the results of that paper extend readily to fronts oriented at an arbitrary angle to the lattice. In that paper, it was found that if the origin $z_0$ varies exponentially slowly with time, the term $C\od z_0/\od T$, where $C$ is some constant, must be added to the right-hand side of (\ref{PinningEqn}). This yields an equation of the form
\begin{equation}
\label{FrontMotion}
\fr{\od \hat{z}_0}{\od \hat{T}} = \delta \hat{r} - \cos \hat{z}_0,
\end{equation}
where $\hat{z}_0$, $\hat{T}$ and $\delta\hat{r}$ are suitably rescaled versions of $z_0$, $T$ and $\delta r$, and are each $\Or(1)$. Thus, within the pinning region the constant values of $z_0$ given by (\ref{PinningEqn}) are stable solutions of (\ref{FrontMotion}), as expected---this simply describes the pinning mechanism. On the other hand, if $\delta r$ increases to just outside the pinning region while remaining exponentially small, (\ref{FrontMotion}) has no real constant solutions and so the front drifts with non-constant velocity, `clicking' through lattice points as described in Section 7 of \cite{king2001asymptotics}. When $\delta r$ increases to $\Or(1)$, $\od z_0/\od T$ is constant and determined by a travelling wave solution to the leading-order reaction-diffusion equation (\ref{RDEq_Z}). This result of \cite{king2001asymptotics}, derived for $\psi = 0$ only, also holds when $\tan\psi \in \mathbb{Q}_\infty$ (after a rescaling to account for the effective lattice spacing); if $\tan\psi$ is irrational and finite, however, there are no exponentially growing terms to eliminate from the remainder and $\od z_0/\od T$ is constant and determined at leading-order.

The methods presented in this paper are readily applicable to lattice systems with difference operators other than $\Delta$ and $\Delta_{\mathrm{hex}}$, e.g. \cite{draelants2014localised, taylor2010snaking}. The results also have implications for finite-difference approximations to differential equations, where approximation of a continuous problem by a lattice creates an artificial pinning region. Furthermore, we note that the formulae describing the bifurcation structure of the pinning region have been derived without explicit knowledge of the leading-order front, only its far-field behaviour; thus our results also apply to systems in which a leading-order front cannot be found analytically, such as that in \cite{yulin2010discrete} containing a nonlinear term of the form $u^3/(1+u^2)$. To this end, we note that in \cite{kozyreff2013analytical}, the imaginary part of the singularity was derived to be $\Im(\zeta) = \pi/\alpha_+$, where $\alpha_+$ is defined in (\ref{alpha}); thus, with the exception of $\Lambda_{1,\psi}$, the formula (\ref{DiscreteSnakingWidth}) for the snaking width is comprised entirely of information which can be determined analytically from (\ref{GenDiffEq}) and (\ref{FrontEq}). The same can not be said of the snakes-and-ladders equations (\ref{SnakesCondition})-(\ref{LaddersCondition}), however, as $D_+$ cannot be found analytically (note that $\Re(\zeta)$ can be set to zero due to the invariance of (\ref{FrontEq}) under spatial translations). Nevertheless, this property will be invaluable in higher-dimensional studies, where analytical results are scarce; we therefore expect the present work to provide a valuable stepping-stone towards the analysis of fully two-dimensional localised solutions \cite{avitabile2010snake, chong2009multistable, lloyd2008localized, lloyd2009localized, taylor2010snaking}. 


\Appendix 

\section{On complex solutions of the eigenvalue equation (\ref{kappa_cond_O1})}
\label{App_Complexkappa}

Suppose that $\kappa = a + ib$ is such that (\ref{kappa_cond_O1}) holds, i.e.
\begin{equation}
\label{A_kappa_cond_O1}
\cos(\kappa\cos\psi) + \cos(\kappa\sin\psi) - 2 = 0
\end{equation}
and $a$ and $b$ are real, non-zero constants. As $\sin^2\Theta + \cos^2\Theta \equiv 1$, (\ref{kappa_cond_Oeps}) can be rewritten to give
\begin{equation}
\label{A_kappa_cond_Oeps_cos}
\cos^2\psi\lb( 1 - \cos^2(\kappa\cos\psi) \rb) = \sin^2\psi \lb( 1 - \cos^2(\kappa\sin\psi) \rb).
\end{equation}
Thus (\ref{A_kappa_cond_O1}) and (\ref{A_kappa_cond_Oeps_cos}) taken together may be formulated as a system of two algebraic equations, treating $\cos(\kappa\cos\psi)$ and $\cos(\kappa\sin\psi)$ as two unknown constants. Of course, in actuality there is only one unknown, the eigenvalue $\kappa$; any solution must therefore provide a consistent value of $\kappa$. 

Solving this system is a simple exercise, and we find upon doing so that either \linebreak $\cos(\kappa\cos\psi) = \cos(\kappa\sin\psi) = 1$, or
\begin{equation}
\cos(\kappa\cos\psi) = \fr{3\tan^2\psi + 1}{\tan^2\psi - 1},\qquad \cos(\kappa\sin\psi) = - \fr{\tan^2\psi + 3}{\tan^2\psi - 1}.
\end{equation}
The first instance gives real $\kappa$ and is simply the solution given by (\ref{tanphi_realkappa}) and (\ref{realkappa}), which we have already discussed fully in Section \ref{Sec_RemEqn}. In the second instance, separating $\kappa$ into its real and imaginary parts, we have
\begin{align}
\cos(a\cos\psi)\cosh(b\cos\psi) - i\sin(a\cos\psi)\sinh(b\cos\psi) & = \fr{3\tan^2\psi + 1}{\tan^2\psi - 1},
\\
\cos(a\sin\psi)\cosh(b\sin\psi) - i\sin(a\sin\psi)\sinh(b\sin\psi) & = - \fr{\tan^2\psi + 3}{\tan^2\psi - 1}.
\end{align}
However, the imaginary part of both of the above equations must vanish, as the right-hand side of each is real. Therefore $\sin(a\cos\psi) = \sin(a\sin\psi) = 0$, giving $a\cos\psi = 2M_1\pi$, $a \sin\psi = 2M_2\pi$ for $(M_1,M_2) \in\mathbb{Z}^2$. However, (\ref{A_kappa_cond_O1}) now reads 
\begin{equation}
\cos(\kappa\cos\psi) + \cos(\kappa\sin\psi) = \cosh(b\cos\psi) + \cosh(b\sin\psi) = 2.
\end{equation}
This has real solutions only if $b\cos\psi = b\sin\psi = 0$, which gives $b = 0$, a contradiction as $b = \Im(\kappa) \neq 0$. Thus there are no solutions $\kappa$ to (\ref{kappa_cond_O1}) with $\Im(\kappa)\neq0$ that also satisfy (\ref{kappa_cond_Oeps}).


\section*{Acknowledgements}
The authors are grateful to the referees for their comments on an early version of this manuscript. A.D. Dean is grateful to the EPSRC for a Ph.D. studentship, in which much of this work was carried out, and the award of a Doctoral Prize Fellowship. 


\bibliography{summat}
\bibliographystyle{siam}

\end{document}